\documentclass[a4paper, 10pt, american]{amsart}

\usepackage{times,latexsym,amssymb}
\usepackage{amsmath,amsthm,bm}
\usepackage{xcolor}
\usepackage{tikz}
\usetikzlibrary{arrows.meta}
\usepackage[colorlinks,pdfpagelabels,pdfstartview = FitH,bookmarksopen= true,bookmarksnumbered = true,linkcolor = blue,plainpages = false,hypertexnames = false,citecolor = red,pagebackref=false]{hyperref}
\usepackage{amsbsy}
\usepackage{amstext}
\usepackage{amssymb}
\usepackage{stmaryrd}
\usepackage{esint}
\usepackage{stmaryrd}
\usepackage{graphicx}
\usepackage{floatflt}
\usepackage{float} 
\usepackage{appendix}
\usepackage{mathrsfs}
\usepackage{psfrag}
\usepackage{caption}
\usepackage[normalem]{ulem}

\SetSymbolFont{stmry}{bold}{U}{stmry}{m}{n}

\allowdisplaybreaks
\newtheorem{theorem}{Theorem}
\newtheorem{lemma}[theorem]{Lemma}
\newtheorem{definition}[theorem]{Definition}
\newtheorem{proposition}[theorem]{Proposition}

\newtheorem{remark}[theorem]{Remark}
\numberwithin{theorem}{section}
\numberwithin{equation}{section}

\def\N{\mathbb{N}}
\def\R{\mathbb{R}}

\renewcommand{\d}{\mathrm{d}}
\newcommand{\dx}{\mathrm{d}x}
\newcommand{\dy}{\mathrm{d}y}

\newcommand{\dt}{\mathrm{d}t}

\newcommand{\dtau}{\mathrm{d}\tau}

\renewcommand{\epsilon}{\varepsilon}

\DeclareMathOperator{\Div}{div}

\DeclareMathOperator{\dist}{dist}
\DeclareMathOperator{\loc}{loc}
\newcommand{\osc}{\operatornamewithlimits{osc}}

\renewcommand{\epsilon}{\varepsilon}

\renewcommand{\rho}{\varrho}

\def\lkk{[\![}
\def\rkk{]\!]}
\def\eqn#1$$#2$${\begin{equation}\label#1#2\end{equation}}

\newcommand{\wto}{\rightharpoondown}

\newcommand{\pl}{\partial}

\newcommand{\dl}{\delta}

\newcommand{\varep}{\varepsilon}

\newcommand{\z}{\zeta}

\newcommand{\sfB}{\boldsymbol{\mathsf{B}}}
\newcommand{\sfP}{\boldsymbol{\mathsf{P}}}
\newcommand{\sfY}{\boldsymbol{\mathsf{Y}}}
\newcommand{\sfM}{\boldsymbol{\mathsf{M}}}

\newcommand{\sfA}{\boldsymbol{\mathsf{A}}}
\newcommand{\sfC}{\boldsymbol{\mathsf{C}}}
\newcommand{\sfD}{\boldsymbol{\mathsf{D}}}
\newcommand{\sfj}{\boldsymbol{\mathsf{j}}}
\newcommand{\sfF}{\boldsymbol{\mathsf{F}}}
\newcommand{\sfJ}{\boldsymbol{\mathsf{J}}}
\newcommand{\sfL}{\boldsymbol{\mathsf{L}}}
\newcommand{\sfa}{\boldsymbol{\mathsf{a}}}
\newcommand{\sfm}{\boldsymbol{\mathsf{m}}}
\newcommand{\sfw}{\boldsymbol{\mathsf{w}}}
\newcommand{\sfR}{\boldsymbol{\mathsf{R}}}

\newcommand{\sfK}{\boldsymbol{\mathsf{K}}}
\newcommand{\sfE}{\boldsymbol{\mathsf{E}}}

\newcommand{\fprod}{%
  \mathbin{\vcenter{\hbox{\scalebox{0.50}{$\bullet$}}}}%
}
\let\TeXchi\chi
\newbox\chibox
\setbox0 \hbox{\mathsurround0pt $\TeXchi$}
\setbox\chibox \hbox{\raise\dp0 \box 0 }
\def\chi{\copy\chibox}

\def\Xint#1{\mathchoice
    {\XXint\displaystyle\textstyle{#1}}%
    {\XXint\textstyle\scriptstyle{#1}}%
    {\XXint\scriptstyle\scriptscriptstyle{#1}}%
    {\XXint\scriptscriptstyle\scriptscriptstyle{#1}}%
    \!\int}
\def\XXint#1#2#3{\setbox0=\hbox{$#1{#2#3}{\int}$}
    \vcenter{\hbox{$#2#3$}}\kern-0.5\wd0}

\def\dashint{\Xint{\raise4pt\hbox to7pt{\hrulefill}}}

\def\XXiint#1#2#3{\setbox0=\hbox{$#1{#2#3}{\iint}$}
    \vcenter{\hbox{$#2#3$}}\kern-0.5\wd0}

\newcommand{\D}{\nabla}
\DeclareMathOperator{\diver}{div}
\DeclareMathOperator{\tr}{tr}
\newcommand{\calL}{\boldsymbol{\mathcal{L}}}


\begin{document}

\title[Lipschitz bounds and gradient continuity for the Trudinger equation]
{Lipschitz bounds and gradient continuity\\ for the Trudinger equation}
\date{\today}

\author[V. B\"ogelein]{Verena B\"{o}gelein}
\address{Verena B\"ogelein,
Fachbereich Mathematik, Universit\"at Salzburg,
Hellbrunner Str. 34, 5020 Salzburg, Austria}
\email{verena.boegelein@plus.ac.at}

\author[F. Duzaar]{Frank Duzaar}
\address{Frank Duzaar,
Fachbereich Mathematik, Universit\"at Salzburg,
Hellbrunner Str. 34, 5020 Salzburg, Austria}
\email{frankjohannes.duzaar@plus.ac.at}

\author[U. Gianazza]{Ugo Gianazza}
\address{Ugo Gianazza,
Dipartimento di Matematica ``F. Casorati",
Universit\`a di Pavia,
via Ferrata 5, 27100 Pavia, Italy}
\email{ugogia04@unipv.it}

\author[N. Liao]{Naian Liao}
\address{Naian Liao,
Fachbereich Mathematik, Universit\"at Salzburg,
Hellbrunner Str. 34, 5020 Salzburg, Austria}
\email{naian.liao@plus.ac.at}

\author[C. Scheven]{Christoph Scheven}
\address{Christoph Scheven,
Fakult\"at f\"ur Mathematik, 
Universit\"at Duisburg-Essen, Thea-Leymann-Str. 9, 45127 Essen, Germany}
\email{christoph.scheven@uni-due.de}

\begin{abstract}
We establish local spatial Lipschitz and gradient H\"older estimates
for general nonnegative weak solutions to Trudinger's
equation $\partial_t u^{p-1}=\Delta_pu$, $p>1$.
The estimates are uniform across the zero set.
The proofs combine Bernstein-type estimates for
power transformations, smooth approximation,
and comparison with exponential barriers
near the zero set.
\end{abstract}

\maketitle

\setcounter{tocdepth}{1}
 \tableofcontents

\section{Introduction}

We consider nonnegative weak solutions to the homogeneous doubly nonlinear equation
\begin{equation}\label{PDE}
    \partial_t u^{p-1}
    - \Div\big(|\nabla u|^{p-2}\nabla u\big)=0
    \quad\mbox{in $\Omega_T$,}
\end{equation}
also known as Trudinger's equation. Here $p>1$, $\Omega\subset\R^N$
is a bounded open set, $T>0$, and $\Omega_T:=\Omega\times(0,T]$.

Equation~\eqref{PDE} is a benchmark in the nonlinear parabolic regularity theory: although nonlinear
in both the solution and its spatial gradient, it is invariant under
multiplication of the solution by positive constants. For $p\neq2$,
however, invariance under the addition of constants is lost, creating
a basic difficulty in the regularity analysis. While H\"older continuity
is well established
(see, e.g., \cite{BDL-holder, KuSiUr, KuLaSiUr}), a central question
is whether spatial gradient estimates remain uniform as the solution
approaches zero.

Our first main result answers this question by establishing   local Lipschitz continuity in space.
The estimates do not depend on a positive lower bound for the solution
and therefore apply also across its zero set.

\begin{theorem}[Lipschitz estimates]\label{thm:Lipschitz}
Let $u$ be a locally bounded, nonnegative, local weak solution to \eqref{PDE} in
$\Omega_T$. Then, we have 
$$
    \nabla u
    \in 
    L^{\infty}_{\rm loc}(\Omega_T,\R^N).
$$
Moreover, there is a constant $C=C(N,p)$ such that for any cylinder $Q_{2R}(z_o)\Subset\Omega_T$ there holds
\begin{equation*}
    \|\D u\|_{L^\infty(Q_R(z_o))}
    \le
    \frac{C}{R} \|u\|_{L^\infty(Q_{2R}(z_o))}.
\end{equation*}
\end{theorem}

We refer to Definition~\ref{def:weaksol} below for the precise notion of a weak solution. The boundedness assumption in Theorem~\ref{thm:Lipschitz} is supported by the known regularity theory;
see~\cite{KiKu,Tr}.

Next, we show that the spatial gradient is not only locally bounded but also locally H\"older continuous.

\begin{theorem}[Gradient H\"older estimate]\label{thm:holder}
Under the hypotheses of Theorem~\ref{thm:Lipschitz}, the continuous representative of $u$ is spatially differentiable at every point of $\Omega_T$ and the spatial gradient is locally Hölder continuous for some Hölder exponent $\alpha=\alpha(N,p)\in(0,1)$.
Moreover, there is a constant $C=C(N,p)>1$ such that
\begin{equation}\label{nn:eq:holder-main}
 |\D u(x,t)-\D u(y,s)|
 \le \frac{C}{R} \|u\|_{L^\infty(Q_{2R}(z_o))}
 \bigg(\frac{|x-y|+|t-s|^{\frac1p}}{R}\bigg)^\alpha
\end{equation}
for all $(x,t),(y,s)\in Q_R(z_o)$. In addition, we have $\D u=0$ on $\{u=0\}$.
\end{theorem}

The significance of these results lies in the uniformity of
the gradient estimates as the solution approaches zero.
For strictly positive solutions, regularity estimates may
depend on a positive lower bound for the solution and therefore
provide no uniform control when that bound vanishes.
Removing this dependence addresses a structural difficulty:
the coefficient $(p-1)u^{p-2}$ of the time derivative degenerates
or becomes singular at zero, while the lack of invariance under
the addition of constants prevents a direct reduction to the
strictly positive case.

Theorems~\ref{thm:Lipschitz} and~\ref{thm:holder} establish
quantitative first-order spatial regularity throughout the
all nonnegative solutions, with estimates that remain valid near the
zero set. They treat the homogeneous endpoint $q=p-1$
of the general equation,
\begin{equation}\label{PDE-doubly}
    \partial_t u^q=\Delta_pu,\qquad q>0,
\end{equation}
which lies outside the gradient regularity
theory developed in~\cite{BDGLS-23, p-Schauder} for $q>p-1$.
Thus the results extend the regularity theory of Trudinger's
equation from H\"older continuity of the solution to the gradient level, without a
quantitative positivity assumption.

\subsection{Previous regularity results}

The regularity theory for~\eqref{PDE} originates in the work of
Trudinger~\cite{Tr}, who established a parabolic Harnack inequality.
Alternative approaches to Harnack's inequality for nonnegative weak
solutions were subsequently developed by Gianazza and
Vespri~\cite{GiVe} and by Kinnunen and Kuusi~\cite{KiKu}.

H\"older continuity for the broader family of doubly nonlinear
parabolic equations \eqref{PDE-doubly}
was studied in various ranges of $p$ and $q$ in the early 1990s by
Chen, Ivanov, Porzio, Vespri, and
Zhou~\cite{Chen-Zhou, Ivanov1, Ivanov2, Ivanov3, PoVe, Ve}. Such theory is still fragmented.
For Trudinger's equation, Kuusi, Laleoglu, Siljander, and
Urbano~\cite{KuSiUr, KuLaSiUr} established local H\"older continuity
of nonnegative weak solutions using Harnack's inequality.
Interior and boundary H\"older continuity for sign-changing
solutions was later obtained by B\"ogelein, Duzaar, and
Liao~\cite{BDL-holder}.

The regularity of the spatial gradient presents further difficulties.
An earlier approach was proposed in~\cite{Siljander}; however,
the proof of~\cite[Lemma~3.3]{Siljander} contains an oversight in the
treatment of the time-derivative term, already in the case of
Lebesgue measure.

For the broader equation~\eqref{PDE-doubly}, B\"ogelein, Duzaar,
Gianazza, Liao, and Scheven~\cite{BDGLS-23, p-Schauder} established local
H\"older continuity of the spatial gradient in a range with
$q>p-1$. Their argument relies on a time-insensitive Harnack
inequality and does not extend to the borderline case $q=p-1$
corresponding to Trudinger's equation.

For Trudinger's equation itself, Lindqvist, Parviainen, and
Siltakoski~\cite{LPS} proved local Lipschitz regularity when solutions are strictly positive by combining the Ishii--Lions method
with the equivalence of positive weak and viscosity solutions.
The resulting Lipschitz estimate depends on a positive lower
bound for the solution and deteriorates as that bound tends
to zero; see~\cite[Corollary~3.2]{LPS}.
Lipschitz regularity under such a positivity assumption also
holds for the more general equation~\eqref{PDE-doubly};
see~\cite[Theorem~6.8]{p-Schauder}.

A related line of research concerns higher integrability of
the spatial gradient. B\"ogelein, Duzaar, Kinnunen, and
Scheven~\cite{BoDuKiSc} established a Gehring-type
self-improvement in the range
    $\max\{1,\frac{2N}{N+2}\}
    <p<\frac{2N}{(N-2)_+}$.
Their result applies to systems and uses an intrinsic geometry
depending on both the solution and its spatial gradient.
For nonnegative scalar solutions, Saari and
Schwarzacher~\cite{SaSc} subsequently removed the upper
restriction on $p$, obtaining higher integrability throughout
the range $p\ge2$.

The central issue addressed here is pointwise gradient
regularity as the solution approaches zero.
Our main results, Theorems~\ref{thm:Lipschitz}
and~\ref{thm:holder}, establish gradient estimates for
nonnegative weak solutions that remain uniform near the zero
set, without assuming a positive lower bound for the solution.

\subsection{Ideas of the proof}
The argument is developed within the weak formulation and yields quantitative estimates for the spatial gradient. The proof builds on the structure of the positivity set: by the expansion of positivity for Trudinger's equation, a nonnegative solution may vanish identically up to a time $t_\ast$, while it is strictly positive thereafter. Moreover, comparison with an exponential barrier shows that, for a short time after $t_\ast$, the solution remains exponentially small in the interior; see \S\,\ref{subsec:exp-barrier}. 

The key ingredient is a pointwise gradient estimate for a suitable power of the solution on its positivity set. Writing $u=v^k$ with $k>1$ sufficiently large, we employ a localized Bernstein argument to control $\nabla v$ independently of any positive lower bound for $u$. To give a glimpse into the machinery, consider the linear case $p=2$, $w=|\nabla v|^2$ will satisfy
\begin{equation*}
    \bigg(\partial_t - \Delta -\underbrace{2(k-1)\frac{\nabla v}{v}\cdot \nabla}_{\text{drift term}}\bigg) w
    =
    -2|\nabla^2v|^2 -\underbrace{2(k-1)\frac{w^2}{v^2}}_{\text{zeroth-order term}}.
\end{equation*}
At an interior maximum of $w$, the drift term disappears, while the negative zeroth-order term prevents $w$ from becoming too large, suggesting an estimate of the form $w\lesssim v^2$. Notice that this mechanism already requires $k>1$ for the heat equation.

The same idea extends, in a considerably more delicate form, to  Trudinger's equation. Its homogeneity is essential, since it allows the power $k$ to be chosen so that the favorable gradient contribution absorbs the terms generated by differentiating the diffusion coefficient. Indeed, for the general equation~\eqref{PDE-doubly} the transformation $u=v^k$ gives
\[
    \partial_t v
    =
    \frac{k^{p-2}}q\,v^\sigma\Delta_pv +
    \frac{k^{p-2}}q(p-1)(k-1)v^{\sigma-1}|\D v|^p,
    \qquad
    \sigma=2-p+k(p-1-q).
\]
Precisely in the Trudinger case $q=p-1$, the exponent $\sigma$ is independent of $k$.

Since the Bernstein computation involves third-order derivatives, we carry it out for a suitable regularization and then pass to the limit. For $k$ sufficiently large and on every cylinder $B_{2r}(x_o)\times(s-2h,s]$ with $0<\sfm\le u\le\sfM$, this yields the gradient bound
\begin{equation}\label{Bernstein-est}
    \|\nabla v\|_{L^\infty(B_r(x_o)\times(s-h,s))}
    \lesssim
    \sfM^{\frac1k}
    \left(\frac1r+\frac1{h^{\frac1p}}\right),
\end{equation}
with a constant independent of the lower bound $\sfm$. Returning to the original variables gives local spatial Lipschitz continuity together with quantitative decay of $\nabla u$ as $u$ approaches zero. Combined with the exponential decay near the zero set, this allows the estimates on the positivity set and the zero set to be patched together into the local gradient bound stated in Theorem~\ref{thm:Lipschitz}.

The estimate in~\eqref{Bernstein-est} is in fact stronger than the gradient bound asserted in Theorem~\ref{thm:Lipschitz}. As stated in Theorem~\ref{thm:lip-k}, it yields, since $u=v^k$ with $k>1$, the quantitative decay
\[
    |\nabla u|
    \lesssim
    u^{1-\frac1k}
\]
as $u$ approaches zero. Combining this decay near the zero set with the local gradient estimates in the positivity set provided by the Schauder theory in~\cite{p-Schauder}, we obtain the H\"older continuity of the gradient stated in Theorem~\ref{thm:holder}.

To the best of our knowledge, this approach is new. Related ideas, however, appear in earlier work. Kotschwar and Ni~\cite{KN} use a logarithmic transformation, a nonlinear Bochner identity, and the maximum principle for Trudinger's equation, but their estimates involve different quantities, in particular combinations of spatial gradients and time derivatives. 
For a related $p$-Laplace equation with gradient absorption, Bartier and Lauren\c{c}ot~\cite{BL} employ a power transformation and a Bernstein-type argument in which the Hessian terms are controlled by completing squares.

\section{Notation and preliminaries}

For a point  $z_o=(x_o,t_o)\in \R^N\times\R$ in space-time and a radius $\rho>0$ we consider parabolic cylinders of the form
\begin{equation*}
    Q_\rho(z_o)
    :=
    B_\rho(x_o)\times (t_o-\rho^p,t_o],
\end{equation*}
where $B_\rho(x_o)\subset\R^N$ denotes the open ball with radius $\rho>0$ and center $x_o$. 
Occasionally, we also use more general cylinders of the type 
\begin{equation*}
    Q_{r,s}(z_o):= B_r(x_o)\times (t_o-s,t_o]
\end{equation*}
with $r,s>0$.
The parabolic boundary of a cylindrical domain $\Omega_T:=\Omega\times(0,T]$ is given by 
$$
 \pl_\mathrm{par} \Omega_T
    :=
    \big(\overline \Omega\times\{0\}\big) \cup \big(\partial \Omega\times(0,T]\big).
$$
For two points $z_1=(x_1,t_1), z_2=(x_2,t_2)\in\R^{N+1}$ their parabolic distance is defined as 
$$
    \dist_\mathrm{par}(z_1,z_2)
    :=
    |x_1-x_2| + |t_1-t_2|^{\frac{1}{p}}.
$$

The spatial gradient of $u(x,t)$ is denoted by $\nabla u$, whereas the Hessian matrix is $\nabla^2u$. We also use some operators of matrices:
\[
A\otimes B:= (a_i b_j)_{N\times N} \quad\text{for } A=(a_1,\cdots, a_N),\,B=(b_1,\cdots, b_N)
\]
and
\[
\mathbf{A}\fprod\mathbf{B}:=\sum_{1\le i,j\le N} a_{ij} b_{ij} \quad\text{for } \mathbf{A}=(a_{ij})_{N\times N},\, \mathbf{B}=(b_{ij})_{N\times N}.
\]

\begin{definition}\label{def:weaksol}
A function
\begin{equation*}
	u\in C\big(0,T;L^p_{\loc}(E)\big)\cap L^p_{\loc}\big(0,T; W^{1,p}_{\loc}(E)\big)
\end{equation*}
is a local, weak solution to \eqref{PDE}, if the identity
\begin{align*}
	\iint_{\Omega_T}\big[u^{p-1}\varphi_t - |\nabla u|^{p-2}\nabla u\cdot \nabla\varphi\big]\dx\dt
    =
    0
\end{align*}
holds, for any  testing function $\varphi\in C_0^\infty(\Omega_T)$.
\end{definition}

\begin{lemma}[Scaling invariance]\label{lem:scaling}
Let $\rho,\sfM>0$ and $a<b$.
Suppose that $u$ is a nonnegative weak solution
to \eqref{PDE} in
\[
 B_\rho(x_o)\times
 (t_\ast+a\rho^p,t_\ast+b\rho^p]
\]
and satisfies $0\le u\le\sfM$ there.
Then
\[
 v(y,\tau)
 :=
 \frac1{\sfM}
 u(x_o+\rho y,t_\ast+\rho^p\tau)
\]
is a nonnegative weak solution to \eqref{PDE}
in $B_1\times(a,b)$ and satisfies $0\le v\le1$.
\end{lemma}

The boundedness and Hölder continuity of weak solutions was proved in~\cite{BDL-holder,KuSiUr,KuLaSiUr}.
\begin{proposition}\label{prop:holder}
Let $u$ be a nonnegative weak solution to~\eqref{PDE} in $Q_{2R}(z_o)$. Then, $u\in C^\eta_{\loc}(Q_{2R}(z_o))$ for some $\eta=\eta(N,p)\in (0,1)$, and there exists $C=C(N,p)>1$, such that
 for all $\rho\in(0,R)$,
\begin{equation*}
    \osc_{Q_\rho(z_o)} u \le C\bigg(\frac{\rho}{R}\bigg)^\eta \sup_{Q_{R}(z_o)} u .
\end{equation*}
\end{proposition}

\begin{proposition}[Strong maximum principle]
\label{prop:Harnack}
Let $\Omega\subset\R^N$ be open and connected,
and let $u$ be the continuous representative of a nonnegative weak solution
to \eqref{PDE} in $\Omega_T:=\Omega\times(0,T]$.
Then there exists $t_\ast\in[0,T]$ such that 
\[u(x,t)=0
 \qquad\text{for all }(x,t)\in
 \Omega\times(0,t_\ast],
 \]
 and
\[
 u(x,t)>0
 \qquad\text{for all }(x,t)\in
 \Omega\times(t_\ast,T].\,\,
\]
In particular, $t_\ast=0$ corresponds to
$u>0$ throughout $\Omega_T$, whereas
$t_\ast=T$ corresponds to
$u\equiv0$ throughout $\Omega_T$.
\end{proposition}

\begin{proof}
We extend $u$ up to time $2T$ by solving~\eqref{PDE} on $\Omega\times(T,2T]$ with initial datum $u(\cdot,T)$ and zero lateral boundary data.
Pick a point $x_*\in\Omega$ and let $t_*$ be the supremum of $\{0\}\cup\{t\in(0,T]:\,u(x_*,t)=0\}$.
It suffices to show that if
$u(x_\ast,t_\ast)=0$, then
\[
 u\equiv0
 \qquad\text{in }\Omega\times(0,t_\ast].
\]
Indeed, suppose that $u(\tilde x,\tilde t)=0$ for some $\tilde x\in\Omega$ with $\tilde x\neq x_*$ and some $\tilde t>t_*$. Then, by the argument below, $u(x_*,t)=0$ for all $t_*<t\le\tilde t$, contradicting the assumption that $t_*$ is the supremum of the times at which $u(x_*,\cdot)$ vanishes.
We use the following expansion of positivity, cf.~\cite[Proposition~B.1]{BDL-holder}. There exist constants $\eta\in(0,1)$, $\theta>0$, and $\Lambda>\max\{2,\theta\}$, depending only on $N$ and $p$, such that the following holds.
If $r>0$ and
\[
 \overline{B_{\Lambda r}(y)}
 \times[s,s+\Lambda r^p]
 \subset\Omega_T,
\]
and if, for some $\ell>0$,
\[
 u(x,s)\ge\ell
 \qquad\text{for every }x\in B_r(y),
\]
then
\begin{equation}\label{eq:expansion}
 u(x,s+\theta r^p)\ge\eta\ell
 \qquad\text{for every }x\in B_{2r}(y).
\end{equation}

Suppose, towards a contradiction, that
$u(y,s)>0$ for some $y\in\Omega$ and
$0<s<t_\ast$.
We will propagate this positivity along a
chain of balls to obtain $u(x_\ast,t_\ast)>0$,
contrary to our assumption.
By continuity, there exist $\ell>0$ and
$\rho>0$, with $\overline{B_\rho(y)}\subset\Omega$,
such that
\begin{equation}\label{eq:initial}
 u(x,s)\ge\ell
 \qquad\text{for every }x\in B_\rho(y).
\end{equation}
No quantitative bounds for $\ell$ or $\rho$
are needed.
Since $\Omega$ is open and connected, there
exists a polygonal path
$\gamma\colon[0,1]\to\Omega$ joining $y$ to $x_\ast$.
If $y\ne x_\ast$, we parametrize $\gamma$
proportionally to arc length, so that
\[
 \gamma(0)=y,
 \qquad
 \gamma(1)=x_\ast,
 \qquad
 |\gamma(a)-\gamma(b)|\le L|a-b|
 \quad\text{for }a,b\in[0,1],
\]
where $L$ denotes the length of the path.
If $y=x_\ast$, we take the constant path and
set $L=0$.
The image $\gamma([0,1])$ is a compact subset
of $\Omega$. Hence, we can choose $d>0$ such
that its closed $d$-neighborhood satisfies
\[
 \bigl\{
   x\in\R^N:
   \operatorname{dist}(x,\gamma([0,1]))\le d
 \bigr\}
 \subset\Omega.
\]

Set $\Delta:=t_\ast-s>0$.
For an integer $m\ge1$, to be chosen below,
define
\begin{equation*}
 r_m:=\left(\frac{\Delta}{\theta m}\right)^{\frac1p},
 \qquad
 \tau_j:=s+\frac{j\Delta}{m},
 \qquad
 z_j:=\gamma\left(\frac jm\right),
 \quad j=0,\ldots,m.
\end{equation*}
Thus, for $j=0,\ldots,m-1$,
\begin{equation*}
 \tau_{j+1}-\tau_j=\theta r_m^p,
 \qquad
 |z_{j+1}-z_j|\le\frac Lm.
\end{equation*}
Since $p>1$, we have
\[
 \frac{L/m}{r_m}
 =
 L\left(\frac{\theta}{\Delta}\right)^{\frac1p}
 m^{\frac1p-1}
 \longrightarrow0
 \qquad\text{as }m\to\infty.
\]
Moreover, also $r_m\to0$ as $m\to\infty$.
Therefore, we can choose a finite integer $m$
sufficiently large such that
\begin{equation}\label{eq:choices}
 \frac Lm\le r_m,
 \qquad
 r_m<\rho,
 \qquad
 2\Lambda r_m<d,
 \qquad
 \Lambda r_m^p<2T-t_\ast.
\end{equation}
Fix such an $m$ and write $r:=r_m$.
These choices ensure that every enlarged cylinder
\[
 \overline{B_{\Lambda r}(z_j)}
 \times[\tau_j,\tau_j+\Lambda r^p],
 \qquad j=0,\ldots,m-1,
\]
is contained in $\Omega_{2T}$.
Indeed, each $z_j$ lies on the path and
$\Lambda r<d$, so
$\overline{B_{\Lambda r}(z_j)}\subset\Omega$.
Moreover,
\[
 0<s\le\tau_j,
 \qquad
 \tau_j+\Lambda r^p
 \le t_\ast+\Lambda r^p<2T.
\]
The last application of the expansion of positivity
starts at $\tau_{m-1}=t_\ast-\theta r^p$
and yields positivity at $t_\ast$.
However, it requires the solution
to be defined up to
$\tau_{m-1}+\Lambda r^p
=t_\ast+(\Lambda-\theta)r^p>t_\ast$.
The last inequality in \eqref{eq:choices}
ensures that this time is still less than $2T$.

We now prove by induction that
\begin{equation}\label{eq:induction}
 u(x,\tau_j)\ge\eta^j\ell
 \qquad\text{for every }x\in B_r(z_j),
 \quad j=0,\ldots,m.
\end{equation}
For $j=0$, we have $z_o=y$ and $\tau_o=s$,
so the assertion follows from
\eqref{eq:initial} and $r<\rho$.

Suppose that \eqref{eq:induction} holds for
some $j\in\{0,\ldots,m-1\}$.
Applying \eqref{eq:expansion} with center $z_j$,
initial time $\tau_j$, and lower bound
$\eta^j\ell$, we obtain
\[
 u(x,\tau_{j+1})\ge\eta^{j+1}\ell
 \qquad\text{for every }x\in B_{2r}(z_j),
\]
where we used
$\tau_{j+1}=\tau_j+\theta r^p$.
Moreover, \eqref{eq:choices} gives
$|z_{j+1}-z_j|\le r$, and hence
\( B_r(z_{j+1})\subset B_{2r}(z_j)\subset\Omega \).
Thus, \eqref{eq:induction} holds for $j+1$,
which completes the induction.

Since $z_m=x_\ast$ and $\tau_m=t_\ast$, we conclude
that
\[
 u(x_\ast,t_\ast)\ge\eta^m\ell>0,
\]
contrary to $u(x_\ast,t_\ast)=0$.
The chain contains only finitely many steps,
so the lower bound $\eta^m\ell$ remains
strictly positive; no passage to the limit
as $m\to\infty$ is involved. It follows that $u\equiv0$ in
$\Omega\times(0,t_\ast)$.
By continuity, $u(\cdot,t_\ast)=0$ as well.
\end{proof}


\begin{figure}[htbp]
\centering
\begin{tikzpicture}[
  x={(1.25cm,0cm)},
  y={(0.55cm,0.32cm)},
  z={(0cm,1.25cm)},
  >=Latex,
  every node/.style={font=\small}
]

\path[fill=gray!5,draw=gray!60,thick]
  (-1.5,-1.5,0)
  -- (5.2,-1.5,0)
  -- (5.2,2,0)
  -- (-1.5,2,0)
  -- cycle;

\node at (4.9,1.5,0) {$\Omega$};

\foreach \j/\xx/\yy/\tt in {
  0/0/0/0,
  1/0.8/0.3/0.65,
  2/1.6/0.6/1.30,
  3/2.4/0.3/1.95,
  4/3.2/0/2.60,
  5/4/0.3/3.25
}{
  \coordinate (base\j) at (\xx,\yy,0);
  \coordinate (point\j) at (\xx,\yy,\tt);
}

\draw[blue!70!black,thick]
  (base0) -- (base1) -- (base2)
  -- (base3) -- (base4) -- (base5);

\node[blue!70!black,below] at (2.0,-0.4,0)
  {$\gamma$};

\foreach \j in {1,...,5}{
  \draw[gray!50,densely dotted]
    (base\j) -- (point\j);
}

\draw[->,thick]
  (-1.5,-1.5,-0.15) -- (-1.5,-1.5,3.85)
  node[above] {$\tau$};

\foreach \height/\labeltext in {
  0/{s=\tau_o},
  0.65/{\tau_1},
  1.30/{\tau_2},
  1.95/{\tau_3},
  2.60/{\tau_4},
  3.25/{t_\ast=\tau_5}
}{
  \draw
    (-1.58,-1.5,\height)
    -- (-1.42,-1.5,\height);
  \node[anchor=east] at (-1.65,-1.5,\height)
    {$\labeltext$};
}

\draw[->]
  (-1.5,-1.5,0) -- (5.65,-1.5,0)
  node[right] {$x_1$};

\draw[->]
  (-1.5,-1.5,0) -- (-1.5,2.65,0)
  node[above] {$x_2$};

\foreach \j/\xx/\yy/\tt in {
  0/0/0/0,
  1/0.8/0.3/0.65,
  2/1.6/0.6/1.30,
  3/2.4/0.3/1.95,
  4/3.2/0/2.60,
  5/4/0.3/3.25
}{
  \path[
    fill=orange!15,
    fill opacity=0.65,
    draw=orange!80!black,
    thick
  ]
    plot[
      domain=0:360,
      samples=65,
      variable=\a
    ]
    ({\xx+0.9*cos(\a)},
     {\yy+0.9*sin(\a)},
     \tt)
    -- cycle;
}

\foreach \j/\next in {0/1,1/2,2/3,3/4,4/5}{
  \draw[->,blue!70!black,thick,shorten >=3pt]
    (point\j) -- (point\next);
}

\foreach \j in {0,...,5}{
  \fill[blue!70!black]
    (point\j) circle[radius=2pt];
}

\filldraw[
  fill=green!60!black,
  draw=white,
  line width=0.6pt
]
  (point4) circle[radius=3.5pt];

\filldraw[
  fill=red,
  draw=white,
  line width=0.6pt
]
  (point5) circle[radius=3.5pt];

\node[below left=3pt] at (point0)
  {$(y,s)$};

\node[above=8pt] at (point5)
  {$(x_\ast,t_\ast)$};

\node[
  anchor=west,
  align=left
] at (5.0,0.3,3.25)
  {$u(\cdot,\tau_5)\ge\eta^5\ell$\\
   on $B_r(z_5)$};

\node[
  anchor=west,
  text=blue!70!black
] at (-1.5,-1.5,-0.65)
  {Blue line in the base plane: polygonal path $\gamma$};

\node[
  anchor=west,
  text=orange!80!black
] at (-1.5,-1.5,-1.0)
  {Orange disks: $B_r(z_j)\times\{\tau_j\}$};

\node[
  anchor=west,
  text=green!50!black
] at (-1.5,-1.5,-1.35)
  {Green point: $(z_4,\tau_4)$};

\node[
  anchor=west,
  text=red!75!black
] at (-1.5,-1.5,-1.7)
  {Red point: $(z_5,\tau_5)=(x_\ast,t_\ast)$};

\end{tikzpicture}
\caption{Schematic space-time representation of
the positivity chain for $m=5$.
The horizontal disks represent spatial balls
at successive times, with
$\tau_{j+1}-\tau_j=\theta r^p$.
The dotted lines project their centers onto
the spatial path.
For clarity, the $d$-neighborhood and the
enlarged comparison cylinders are omitted.
}
\label{fig:positivity-chain-spacetime}
\end{figure}
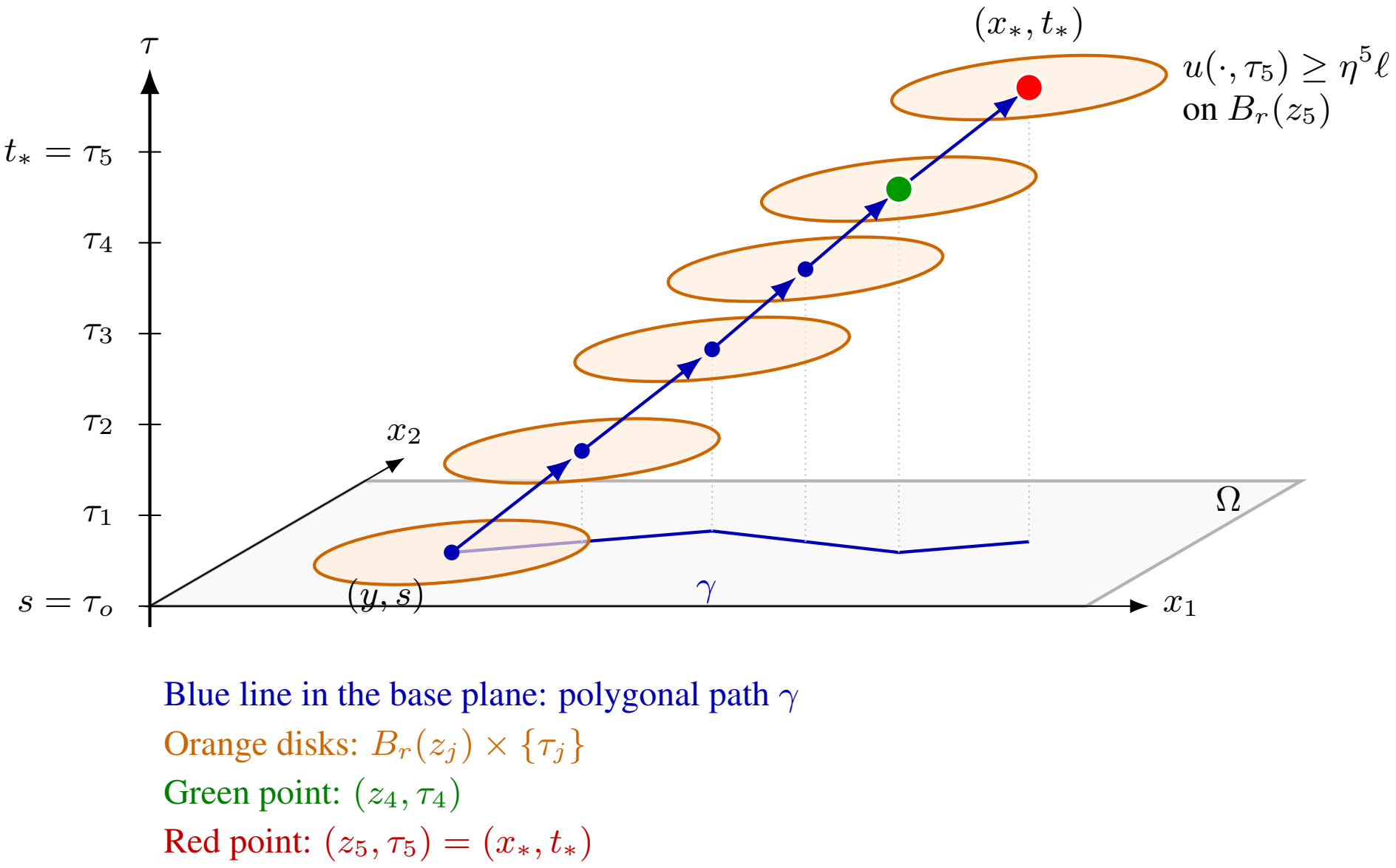

\section{Gradient bounds for smooth solutions}
This section treats smooth positive solutions $u$ to \eqref{PDE} and establishes an upper bound for the gradient of $u^{\frac1k}$, for suitable $k$. 
For the reader's convenience, we outline the four steps of the argument:
\begin{enumerate}
\item[(a)] the power transformation $u=v^k$ turns \eqref{PDE} into 
\[
v_t=a(v)\Delta_pv+b(v)|\nabla v|^p,;
\]
where $a$ and $b$ are defined in~\eqref{nn:eq:coefficients};
\item[(b)] differentiating in $x$ gives a Bochner-type identity 
\[
\calL w=-2a\,w^{\frac{p-2}2}\tr(\sfP(\nabla^2v)^2)+2a'w^{\frac p2}\,\sfP\fprod\nabla^2v+2b'w^{\frac{p+2}2}
\]
for $w=|\nabla v|^2$, where the operator $\calL$ is defined below in \eqref{Def:cal-L}. The first-order term, with $b'<0$, is the only source of coercivity;
\item[(c)] Cauchy-Schwarz and Young's inequalities remove the mixed Hessian term at the price of the term $(N+p-2)(a')^2/a$, which is compensated by $2b'$ as soon as $k\ge k_o$;
\item[(d)] a cutoff function $\eta$ and the maximum principle at a maximum point of $\eta^2w$ localize the estimate.
\end{enumerate}
A more detailed discussion is given at the end of the section.
We consider equation~\eqref{PDE}
in $\Omega_T$, and we set
\begin{equation}\label{nn:eq:k}
 A_*=(N+p-2)(p-2)^2,\qquad
 k_o=2+\frac{A_*}{(p-1)^2},\qquad
 k\ge k_o.
\end{equation}
The following calculation is carried out locally in a region
where $u>0$ and $\D v\ne0$, with $v=u^\frac1k$. We assume
sufficient smoothness for all derivatives below to exist and
for the relevant mixed derivatives to commute. For weak
solutions, this differentiation requires a separate
justification. This is the reason why in Section~\ref{nn:sec:regularization} we consider solutions to a regularized equation. The specific lower bound on $k$ is needed later, but not for the identities derived below.

\subsection{The transformation $u=v^k$}

Let $u$ be a smooth nonnegative solution to \eqref{PDE}.
On the open set where $u>0$, we define
\[
 v=u^\frac1k.
\]
The transformation is smooth on this set. In deriving
the pointwise differentiated identity below, we work
at points where $\nabla v\ne0$. This is a local
restriction for the calculation, not an additional
assumption on the solution. The resulting identity is therefore established on
\[
 \bigl\{(x,t)\in\Omega_T:
          u(x,t)>0,\ |\nabla v(x,t)|>0\bigr\}.
\]
Its use at points where $u=0$ or $\nabla v=0$
requires a separate justification.
Using $v=u^{\frac1k}$ and applying the chain rule, we obtain
\begin{align*}
 v_t
 &=\partial_t\bigl((u^{p-1})^{\frac{1}{k(p-1)}}\bigr)\\
 &=\frac1{k(p-1)}
   u^{\frac1k-(p-1)}\partial_t(u^{p-1})\\
 &=\frac1{k(p-1)}
   v^{1-k(p-1)}
   \diver\bigl(|\D(v^k)|^{p-2}\D(v^k)\bigr).
\end{align*}
Since
\[
 \D(v^k)=kv^{k-1}\D v,
\]
we have
\[
 |\D(v^k)|^{p-2}\D(v^k)
 =k^{p-1}v^{(k-1)(p-1)}|\D v|^{p-2}\D v.
\]
Consequently, the product rule gives
\begin{align*}
 v_t
 &=\frac{k^{p-2}}{p-1}v^{1-k(p-1)}
   \diver\Bigl[
     v^{(k-1)(p-1)}|\nabla v|^{p-2}\nabla v
   \Bigr]\\
 &=\frac{k^{p-2}}{p-1}v^{1-k(p-1)}
   \Bigl[
     v^{(k-1)(p-1)}
     \diver\bigl(|\nabla v|^{p-2}\nabla v\bigr)\\
 &\hspace{36mm}
     +\nabla\bigl(v^{(k-1)(p-1)}\bigr)
      \cdot\bigl(|\nabla v|^{p-2}\nabla v\bigr)
   \Bigr]\\
 &=\frac{k^{p-2}}{p-1}v^{1-k(p-1)}
   \Bigl[
     v^{(k-1)(p-1)}\Delta_pv\\
 &\hspace{36mm}
     +(k-1)(p-1)v^{(k-1)(p-1)-1}
       |\nabla v|^{p-2}
       \underbrace{\nabla v\cdot\nabla v}_{=|\nabla v|^2}
   \Bigr]\\
 &=\frac{k^{p-2}}{p-1}v^{1-k(p-1)}
   \Bigl[
     v^{(k-1)(p-1)}\Delta_pv\\
 &\hspace{36mm}
     +(k-1)(p-1)
       v^{(k-1)(p-1)-1}|\nabla v|^p
   \Bigr].
\end{align*}
The exponents simplify
according to
\begin{equation*}
 1-k(p-1)+(k-1)(p-1)=2-p
\end{equation*}
and 
\begin{equation*}
 1-k(p-1)+(k-1)(p-1)-1=1-p.   
\end{equation*}
Combining the powers of $v$, we obtain
\begin{equation}\label{nn:eq:v}
 v_t=a(v)\Delta_pv+b(v)|\nabla v|^p,
\end{equation}
where
\begin{equation}\label{nn:eq:coefficients}
 a(v)=cv^{2-p},\qquad
 b(v)=c(p-1)(k-1)v^{1-p},\qquad
 c=\frac{k^{p-2}}{p-1}.
\end{equation}
In particular, $a,b>0$, and
\[
 a'(v)=c(2-p)v^{1-p},
 \qquad
 b'(v)=-c(p-1)^2(k-1)v^{-p}.
\]
We point out that the only effect of the power transformation is the new first-order term $b(v)|\nabla v|^p$ (for $k=1$ one has $b\equiv0$). The sign $b>0$, which requires $k>1$, is used in the drift estimate below, and $b'<0$ will provide the coercivity.

\subsection{Expanding the $p$-Laplace operator}

At points where $\nabla v\ne0$, the product rule gives
\begin{align*}
 \Delta_pv
 &=\diver\bigl(|\nabla v|^{p-2}\nabla v\bigr)
 =|\nabla v|^{p-2}\Delta v
   +\nabla\bigl(|\nabla v|^{p-2}\bigr)\cdot\nabla v.
\end{align*}
To compute the second term, we use the chain rule:
\begin{align*}
 \nabla\bigl(|\nabla v|^{p-2}\bigr)
 &=\nabla\bigl((|\nabla v|^2)^{(p-2)/2}\bigr)\\
 &=\frac{p-2}{2}|\nabla v|^{p-4}
   \nabla\bigl(|\nabla v|^2\bigr)\\
 &=(p-2)|\nabla v|^{p-4}\nabla^2v\,\nabla v.
\end{align*}
Here $\nabla^2v$ denotes the Hessian matrix of $v$, and
we have used
\[
 \nabla\bigl(|\nabla v|^2\bigr)
 =2\nabla^2v\,\nabla v.
\]
Substituting this expression yields
\begin{equation}\label{nn:eq:plaplace}
 \Delta_pv
 =|\nabla v|^{p-2}\Delta v
  +(p-2)|\nabla v|^{p-4}
    \bigl(\nabla^2v\,\nabla v\bigr)\cdot\nabla v.
\end{equation}
Consequently, equation \eqref{nn:eq:v} takes the form
\begin{align*}
 v_t
 &=a(v)|\nabla v|^{p-2}\Delta v+(p-2)a(v)|\nabla v|^{p-4}
   \bigl(\nabla^2v\,\nabla v\bigr)\cdot\nabla v+b(v)|\nabla v|^p.
\end{align*}
At points where $\nabla v\ne0$, we can factor out
$a(v)|\nabla v|^{p-2}$ and write
\begin{align*}
 v_t
 &=a(v)|\nabla v|^{p-2}
   \biggl[
     \Delta v
     +(p-2)\nabla^2v\fprod
     \left(
       \frac{\nabla v}{|\nabla v|}
       \otimes
       \frac{\nabla v}{|\nabla v|}
     \right)
   \biggr]
   +b(v)|\nabla v|^p.
\end{align*}
Here the bold dot denotes the Frobenius inner product
of matrices, so that
\[
 \nabla^2v\fprod
 \left(
   \frac{\nabla v}{|\nabla v|}
   \otimes
   \frac{\nabla v}{|\nabla v|}
 \right)
 =
 \frac{(\nabla^2v\,\nabla v)\cdot\nabla v}
      {|\nabla v|^2}.
\]

\subsection{Differentiating the equation}
We take the spatial gradient of~\eqref{nn:eq:v} and then take the scalar product with $2\nabla v$. The product and chain rules give
\begin{align*}
 \partial_t|\nabla v|^2
 &=2\nabla v\cdot\nabla \partial_tv
 =2\nabla v\cdot\nabla
   \bigl(a(v)\Delta_pv+b(v)|\nabla v|^p\bigr)\\
 &=2a'(v)|\nabla v|^2\Delta_pv
   +2a(v)\nabla v\cdot\nabla\Delta_pv\\
 &\phantom{=\,}
   +2b'(v)|\nabla v|^{p+2}
   +2b(v)\nabla v\cdot\nabla\bigl(|\nabla v|^p\bigr).
\end{align*}
Since
\[
 \nabla\bigl(|\nabla v|^p\bigr)
 =\tfrac12 p|\nabla v|^{p-2}
   \nabla\bigl(|\nabla v|^2\bigr),
\]
we obtain
\begin{align}\label{nn:eq:gradient-time}
 \partial_t|\nabla v|^2
&= 2a'(v)|\nabla v|^2\Delta_pv
     +2a(v)\nabla v\cdot\nabla\Delta_pv\notag\\
 &\phantom{=\,}+2b'(v)|\nabla v|^{p+2}+pb(v)|\nabla v|^{p-2}
   \nabla v\cdot\nabla\bigl(|\nabla v|^2\bigr).
\end{align}
We next rewrite the term containing
$\nabla\Delta_pv$. For each $m\in\{1,\ldots,N\}$,
differentiating the flux gives
\begin{align*}
 \partial_{x_m}\Delta_pv
 &=\diver\left(
     \partial_{x_m}
     \bigl(|\nabla v|^{p-2}\nabla v\bigr)
   \right)\\
 &=\diver\Bigl(
     |\nabla v|^{p-2}\nabla(\partial_{x_m}v)
     +(p-2)|\nabla v|^{p-4}
       \bigl(\nabla v\cdot\nabla(\partial_{x_m}v)\bigr)
       \nabla v
   \Bigr).
\end{align*}
Multiplying by $\partial_{x_m}v$ and summing over $m$,
we obtain
\begin{align*}
 &\nabla v\cdot\nabla\Delta_pv\\
 &\quad =\sum_{m=1}^N
   (\partial_{x_m}v)\,
   \diver\Bigl(
     |\nabla v|^{p-2}\nabla(\partial_{x_m}v)
     +(p-2)|\nabla v|^{p-4}
       \bigl(\nabla v\cdot\nabla(\partial_{x_m}v)\bigr)
       \nabla v
   \Bigr).
\end{align*}
We now apply the product rule
\[
 (\partial_{x_m}v)\,\diver X
 =\diver\bigl((\partial_{x_m}v)X\bigr)
  -\nabla(\partial_{x_m}v)\cdot X
\]
to each summand, with $X$ denoting either of the two
vector fields inside the divergence in the preceding
identity. This gives
\begin{align*}
 \nabla v\cdot\nabla\Delta_pv
 &=
 \diver\Biggl(
   |\nabla v|^{p-2}
   \sum_{m=1}^N
     (\partial_{x_m}v)\nabla(\partial_{x_m}v)\\
 &\phantom{\diver\Biggl(\, }
   +(p-2)|\nabla v|^{p-4}
   \left[
     \sum_{m=1}^N
       (\partial_{x_m}v)
       \bigl(\nabla v\cdot\nabla(\partial_{x_m}v)\bigr)
   \right]\nabla v
 \Biggr)\\
 &\phantom{=\,}
 -|\nabla v|^{p-2}
   \sum_{m=1}^N|\nabla(\partial_{x_m}v)|^2
   -(p-2)|\nabla v|^{p-4}
   \sum_{m=1}^N
   \bigl(\nabla v\cdot\nabla(\partial_{x_m}v)\bigr)^2.
\end{align*}
The two sums inside the divergence on the right-hand
side satisfy
\[
 \sum_{m=1}^N
   (\partial_{x_m}v)\nabla(\partial_{x_m}v)
 =\tfrac12\nabla\bigl(|\nabla v|^2\bigr)
\]
and
\[
 \sum_{m=1}^N
   (\partial_{x_m}v)
   \bigl(\nabla v\cdot\nabla(\partial_{x_m}v)\bigr)
 =\tfrac12\nabla v\cdot\nabla\bigl(|\nabla v|^2\bigr).
\]
By the definition of the Frobenius norm,
\[
 \sum_{m=1}^N|\nabla(\partial_{x_m}v)|^2
 =\sum_{m=1}^N\sum_{j=1}^N
   |\partial_{x_j}\partial_{x_m}v|^2
 =|\nabla^2v|^2.
\]
Moreover, since the mixed second derivatives commute,
\[
 \nabla v\cdot\nabla(\partial_{x_m}v)
 =\sum_{j=1}^N
   (\partial_{x_j}v)
   (\partial_{x_j}\partial_{x_m}v)
 =\sum_{j=1}^N
   (\partial_{x_m}\partial_{x_j}v)
   (\partial_{x_j}v)
 =\bigl(\nabla^2v\,\nabla v\bigr)_m.
\]
Squaring and summing over $m$ therefore gives
\[
 \sum_{m=1}^N
   \bigl(\nabla v\cdot\nabla(\partial_{x_m}v)\bigr)^2
 =\sum_{m=1}^N
   \bigl(\nabla^2v\,\nabla v\bigr)_m^2
 =|\nabla^2v\,\nabla v|^2.
\]
Therefore,
\begin{align}\label{nn:eq:gradient-plaplace}
 2
 \nabla v\cdot\nabla\Delta_pv
 &=
 \diver\Bigl(
   |\nabla v|^{p-2}\nabla\bigl(|\nabla v|^2\bigr)
   +(p-2)|\nabla v|^{p-4}
     \bigl(\nabla v\cdot\nabla(|\nabla v|^2)\bigr)
     \nabla v
 \Bigr)\notag\\
 &\phantom{=\,}
 -2|\nabla v|^{p-2}|\nabla^2v|^2
 -2(p-2)|\nabla v|^{p-4}
   |\nabla^2v\,\nabla v|^2.
\end{align}
We first substitute \eqref{nn:eq:gradient-plaplace}
into \eqref{nn:eq:gradient-time}, leaving the first
term on the right-hand side, i.e.
$2a'(v)|\nabla v|^2\Delta_pv$, unchanged. This gives
\begin{align*}
 \partial_t|\nabla v|^2
 &=
 2a'(v)|\nabla v|^2\Delta_pv\\
 &\phantom{=\,}
    +a(v)\diver\Bigl(
   |\nabla v|^{p-2}\nabla\bigl(|\nabla v|^2\bigr)
   +(p-2)|\nabla v|^{p-4}
   \bigl(\nabla v\cdot\nabla(|\nabla v|^2)\bigr)
   \nabla v
 \Bigr)\\
 &\phantom{=\,}-2a(v)|\nabla v|^{p-2}|\nabla^2v|^2-2(p-2)a(v)|\nabla v|^{p-4}
   |\nabla^2v\,\nabla v|^2\\
 &\phantom{=\,}
 +2b'(v)|\nabla v|^{p+2}+pb(v)|\nabla v|^{p-2}
   \nabla v\cdot\nabla\bigl(|\nabla v|^2\bigr).
\end{align*}
The terms appearing in the second and third lines
on the right-hand side together equal
$2a(v)\nabla v\cdot\nabla\Delta_pv$.
Next, we expand the first term on the right-hand side
using \eqref{nn:eq:plaplace}:
\begin{align*}
 2a'(v)|\nabla v|^2\Delta_pv
 &=
 2a'(v)|\nabla v|^2
 \Bigl[
   |\nabla v|^{p-2}\Delta v\\
 &\hspace{30mm}
   +(p-2)|\nabla v|^{p-4}
   \bigl(\nabla^2v\,\nabla v\bigr)\cdot\nabla v
 \Bigr]\\
 &=
 2a'(v)|\nabla v|^p\Delta v\\
 &\quad
 +2(p-2)a'(v)|\nabla v|^{p-2}
   \bigl(\nabla^2v\,\nabla v\bigr)\cdot\nabla v.
\end{align*}
Substituting this expression into the preceding identity
yields
\begin{align*}
 \partial_t|\nabla v|^2
&=
 a(v)\diver\Bigl(
   |\nabla v|^{p-2}\nabla\bigl(|\nabla v|^2\bigr)
   +(p-2)|\nabla v|^{p-4}
   \bigl(\nabla v\cdot\nabla(|\nabla v|^2)\bigr)
   \nabla v
 \Bigr)\\
 &\phantom{=\,}+pb(v)|\nabla v|^{p-2}
   \nabla v\cdot\nabla\bigl(|\nabla v|^2\bigr)\\
 &\phantom{=\,}-2a(v)
 \Big( |\nabla v|^{p-2}|\nabla^2v|^2
 +(p-2)|\nabla v|^{p-4}
   |\nabla^2v\,\nabla v|^2\Big)\\
 &\phantom{=\,}+2a'(v)
 \Big(|\nabla v|^p\Delta v+(p-2)|\nabla v|^{p-2}
   \bigl(\nabla^2v\,\nabla v\bigr)\cdot\nabla v\Big)\\
 &\phantom{=\,}+2b'(v)|\nabla v|^{p+2}.
\end{align*}
Finally, we move the divergence term and the term
containing $b(v)$ to the left-hand side. We obtain
the identity
\begin{align}\label{nn:eq:differentiated}
 \partial_t|\nabla v|^2
 &-a(v)\diver\Bigl(
   |\nabla v|^{p-2}\nabla\bigl(|\nabla v|^2\bigr)
   +(p-2)|\nabla v|^{p-4}
   \bigl(\nabla v\cdot\nabla(|\nabla v|^2)\bigr)
   \nabla v
 \Bigr)\notag\\
 &-\underbrace{pb(v)|\nabla v|^{p-2}
   \nabla v\cdot\nabla\bigl(|\nabla v|^2\bigr)}_{\text{1st-order term}}\notag\\
 =&
 -2a(v)\Bigl(
   |\nabla v|^{p-2}|\nabla^2v|^2
   +(p-2)|\nabla v|^{p-4}
     |\nabla^2v\,\nabla v|^2
 \Bigr)\notag\\
 &+2a'(v)\Bigl(
   |\nabla v|^p\Delta v
   +(p-2)|\nabla v|^{p-2}
     \bigl(\nabla^2v\,\nabla v\bigr)\cdot\nabla v
 \Bigr)\notag\\
 &+2b'(v)|\nabla v|^{p+2}.
\end{align}
The coefficient $a(v)$ remains outside the divergence
throughout the calculation.

\subsection{Identifying the drift}

We work at points where $\nabla v\ne0$. The explicit first-order term on the left side has been marked in \eqref{nn:eq:differentiated}. To identify
the first-order terms hidden in the divergence structure on the left-hand side of
\eqref{nn:eq:differentiated}, we expand the two
contributions to the divergence separately.
For the first contribution, the product and chain
rules give
\begin{align*}
 \diver&\bigl(
   |\nabla v|^{p-2}\nabla(|\nabla v|^2)
 \bigr)\\
 &=
 |\nabla v|^{p-2}\Delta(|\nabla v|^2)
 +\nabla(|\nabla v|^{p-2})
  \cdot\nabla(|\nabla v|^2)\\
 &=
 |\nabla v|^{p-2}\Delta(|\nabla v|^2)
 +(p-2)|\nabla v|^{p-4}
   (\nabla^2v\,\nabla v)\cdot\nabla(|\nabla v|^2).
\end{align*}
For the second contribution, we first leave aside
the constant factor $(p-2)$ and compute
\begin{align*}
 \diver\Bigl(
   |\nabla v|^{p-4}
   \bigl(\nabla v\cdot\nabla(|\nabla v|^2)\bigr)
   \nabla v
 \Bigr)
 &=
 \bigl(\nabla(|\nabla v|^{p-4})\cdot\nabla v\bigr)
 \bigl(\nabla v\cdot\nabla(|\nabla v|^2)\bigr)\\
 &\quad
 +|\nabla v|^{p-4}
 \nabla\bigl(\nabla v\cdot\nabla(|\nabla v|^2)\bigr)
 \cdot\nabla v\\
 &\quad
 +|\nabla v|^{p-4}
 \bigl(\nabla v\cdot\nabla(|\nabla v|^2)\bigr)
 \Delta v.
\end{align*}
Using
\[
 \nabla(|\nabla v|^{p-4})
 =(p-4)|\nabla v|^{p-6}\nabla^2v\,\nabla v
\]
and
\begin{align*}
 \nabla\bigl(\nabla v\cdot\nabla(|\nabla v|^2)\bigr)
 &
 =\nabla^2v\,\nabla(|\nabla v|^2)
  +\nabla^2(|\nabla v|^2)\nabla v,
\end{align*}
we can expand the derivatives in the preceding
identity. Since $\nabla^2v$ is symmetric, we also have
\[
 \bigl(\nabla^2v\,\nabla(|\nabla v|^2)\bigr)\cdot\nabla v
 =
 (\nabla^2v\,\nabla v)\cdot\nabla(|\nabla v|^2).
\]
Consequently, apart from the factor $(p-2)$ the second contribution is
\begin{align*}
 \diver&\Bigl(
   |\nabla v|^{p-4}
   \bigl(\nabla v\cdot\nabla(|\nabla v|^2)\bigr)
   \nabla v
 \Bigr)\\
 &=
 (p-4)|\nabla v|^{p-6}
 \bigl((\nabla^2v\,\nabla v)\cdot\nabla v\bigr)
 \bigl(\nabla v\cdot\nabla(|\nabla v|^2)\bigr)\\
 &\quad
 +|\nabla v|^{p-4}
 (\nabla^2v\,\nabla v)\cdot\nabla(|\nabla v|^2)\\
 &\quad
 +|\nabla v|^{p-4}
 \bigl(\nabla^2(|\nabla v|^2)\nabla v\bigr)
 \cdot\nabla v\\
 &\quad
 +|\nabla v|^{p-4}\Delta v\,
 \nabla v\cdot\nabla(|\nabla v|^2).
\end{align*}
Multiplying this identity by $(p-2)$ and adding
the expansion of the first contribution, we find
\begin{align*}
 \diver&\Bigl(
   |\nabla v|^{p-2}\nabla(|\nabla v|^2)
   +(p-2)|\nabla v|^{p-4}
   \bigl(\nabla v\cdot\nabla(|\nabla v|^2)\bigr)
   \nabla v
 \Bigr)\\
 &=
 |\nabla v|^{p-2}\Delta(|\nabla v|^2)
 +(p-2)|\nabla v|^{p-4}
 \bigl(\nabla^2(|\nabla v|^2)\nabla v\bigr)
 \cdot\nabla v\\
 &\quad
 +2(p-2)|\nabla v|^{p-4}
 (\nabla^2v\,\nabla v)\cdot\nabla(|\nabla v|^2)\\
 &\quad
 +(p-2)|\nabla v|^{p-4}\Delta v\,
 \nabla v\cdot\nabla(|\nabla v|^2)\\
 &\quad
 +(p-2)(p-4)|\nabla v|^{p-6}
 \bigl((\nabla^2v\,\nabla v)\cdot\nabla v\bigr)
 \nabla v\cdot\nabla(|\nabla v|^2)\\
 &=
 |\nabla v|^{p-2}\Delta(|\nabla v|^2)
 +(p-2)|\nabla v|^{p-4}
 \bigl(\nabla^2(|\nabla v|^2)\nabla v\bigr)
 \cdot\nabla v\\
 &\quad
 +\underbrace{(p-2)|\nabla v|^{p-4}
 \bigl(2\nabla^2v\,\nabla v+\Delta v\,\nabla v\bigr)
 \cdot\nabla(|\nabla v|^2)}_{\text{1st-order term}}\\
 &\quad
 +\underbrace{(p-2)(p-4)|\nabla v|^{p-6}
 \bigl((\nabla^2v\,\nabla v)\cdot\nabla v\bigr)
 \nabla v\cdot\nabla(|\nabla v|^2)}_{\text{1st-order term}}.
\end{align*}
 We have identified the first-order terms hidden in the divergence.
The factor $2$ multiplying
$(\nabla^2v\,\nabla v)\cdot\nabla(|\nabla v|^2)$
after the first equality arises from combining
the identical contributions from the two expansions,
after multiplying the second expansion by $(p-2)$.

Now we can identify the first order terms on the left side of \eqref{nn:eq:differentiated}. Indeed, multiplying the last display by $-a(v)$,
as prescribed by \eqref{nn:eq:differentiated}, and
combining them with the explicit first-order term
\[
 -pb(v)|\nabla v|^{p-2}
 \nabla v\cdot\nabla(|\nabla v|^2),
\]
we can write the first-order terms on the left-hand
side as
\[
 -\sfB\cdot\nabla(|\nabla v|^2),
\]
where
\begin{align}\label{nn:eq:B}
 \sfB
 &=(p-2)a(v)|\nabla v|^{p-4}
   \bigl(2\nabla^2v\,\nabla v+\Delta v\,\nabla v\bigr)
   \notag\\
 &\phantom{=\,}
   +(p-2)(p-4)a(v)|\nabla v|^{p-6}
   \bigl((\nabla^2v\,\nabla v)\cdot\nabla v\bigr)
   \nabla v
   \notag\\
 &\phantom{=\,}
   +pb(v)|\nabla v|^{p-2}\nabla v.
\end{align}
The terms containing $a'(v)$ remain on the
right-hand side.
With this definition, equation
\eqref{nn:eq:differentiated} becomes
\begin{align}\label{nn:eq:with-drift}
 \partial_t|\nabla v|^2
 &
 -a(v)|\nabla v|^{p-2}
 \Biggl[
   \Delta(|\nabla v|^2)
   +(p-2)
   \frac{(\nabla^2(|\nabla v|^2)\nabla v)\cdot \nabla v}{|\nabla v|^2}
 \Biggr]
 -\sfB\cdot\nabla(|\nabla v|^2)\notag\\
 &=
 -2a(v)|\nabla v|^{p-2}
 \Biggl[
   |\nabla^2v|^2
   +(p-2)\frac{|\nabla^2v\,\nabla v|^2}{|\nabla v|^2}
 \Biggr]\notag\\
 &\quad
 +2a'(v)|\nabla v|^p
 \Biggl[
   \Delta v
   +(p-2)
   \frac{(\nabla^2v\,\nabla v)\cdot\nabla v}
        {|\nabla v|^2}
 \Biggr]\notag\\
 &\quad
 +2b'(v)|\nabla v|^{p+2}.
\end{align}

\subsection{Estimating the drift}

We estimate the three contributions to $\sfB$
in \eqref{nn:eq:B}. Observe the elementary inequalities
\[
 |\Delta v|\le\sqrt N\,|\nabla^2v|,
 \qquad
 |\nabla^2v\,\nabla v|
 \le|\nabla^2v|\,|\nabla v|;
\]
the latter implies
\[
 \bigl|(\nabla^2v\,\nabla v)\cdot\nabla v\bigr|
 \le|\nabla^2v|\,|\nabla v|^2.
\]
Here $|\nabla^2v|$ denotes the Frobenius norm.

For the first contribution, we obtain
\begin{align*}
 &\left|
 (p-2)a(v)|\nabla v|^{p-4}
 \bigl(2\nabla^2v\,\nabla v+\Delta v\,\nabla v\bigr)
 \right|\\
 &\qquad
 \le |p-2|a(v)|\nabla v|^{p-4}
 \bigl(
   2|\nabla^2v\,\nabla v|
   +|\Delta v|\,|\nabla v|
 \bigr)\\
 &\qquad
 \le |p-2|(2+\sqrt N)a(v)
 |\nabla v|^{p-3}|\nabla^2v|.
\end{align*}
For the second contribution,
\begin{align*}
 &\left|
 (p-2)(p-4)a(v)|\nabla v|^{p-6}
 \bigl((\nabla^2v\,\nabla v)\cdot\nabla v\bigr)
 \nabla v
 \right|\\
 &\qquad
 \le |p-2||p-4|a(v)|\nabla v|^{p-6}
 |\nabla^2v|\,|\nabla v|^3\\
 &\qquad
 =|p-2||p-4|a(v)
 |\nabla v|^{p-3}|\nabla^2v|.
\end{align*}
Finally, since $b(v)>0$,
\[
 \left|pb(v)|\nabla v|^{p-2}\nabla v\right|
 =pb(v)|\nabla v|^{p-1}.
\]
Combining these estimates gives
\begin{align*}
 |\sfB|
 \le{}&
 |p-2|\bigl(2+\sqrt N+|p-4|\bigr)
 a(v)|\nabla v|^{p-3}|\nabla^2v|\\
 &+pb(v)|\nabla v|^{p-1}.
\end{align*}
Thus,
\begin{equation}\label{nn:eq:drift}
 |\sfB|
 \le C(N,p)\Bigl(
   a(v)|\nabla v|^{p-3}|\nabla^2v|
   +b(v)|\nabla v|^{p-1}
 \Bigr),
\end{equation}
where $C(N,p)$ depends only on $N$ and $p$.

\subsection{Simplifying the notation}

We set
\[
 w:=|\nabla v|^2,
 \qquad
 \sfP:=I+(p-2)
 \frac{\nabla v}{|\nabla v|}
 \otimes
 \frac{\nabla v}{|\nabla v|},
\]
where $I$ denotes the identity matrix.
With the convention $C\fprod D=\tr(C^TD)$, we have
\begin{align*}
 \sfP\fprod\nabla^2w
 &=
 \Delta(|\nabla v|^2)
 +
 (p-2)
   \frac{(\nabla^2(|\nabla v|^2)\nabla v)\cdot \nabla v}{|\nabla v|^2},
\end{align*}
and
\[
 \sfP\fprod\nabla^2v
 =
 \Delta v
 +(p-2)
 \frac{(\nabla^2v\,\nabla v)\cdot\nabla v}
      {|\nabla v|^2}.
\]
Moreover,
\[
 \tr\bigl(\sfP(\nabla^2v)^2\bigr)
 =
 |\nabla^2v|^2
 +(p-2)\frac{|\nabla^2v\,\nabla v|^2}{|\nabla v|^2}.
\]
To verify the last identity, we use the definition
of $\sfP$ and the linearity of the trace:
\begin{align*}
 \tr\bigl(\sfP(\nabla^2v)^2\bigr)
 &=
 \tr\bigl((\nabla^2v)^2\bigr)
 +\frac{p-2}{|\nabla v|^2}
 \tr\bigl((\nabla v\otimes\nabla v)(\nabla^2v)^2\bigr).
\end{align*}
Since $\nabla^2v$ is symmetric,
\[
 \tr\bigl((\nabla^2v)^2\bigr)
 =\tr\bigl((\nabla^2v)^T\nabla^2v\bigr)
 =|\nabla^2v|^2.
\]
Furthermore, the identity
$\tr((z\otimes z)A)=z\cdot Az$ gives
\begin{align*}
 \tr\bigl((\nabla v\otimes\nabla v)(\nabla^2v)^2\bigr)
 &=\nabla v\cdot(\nabla^2v)^2\nabla v\\
 &=(\nabla^2v\,\nabla v)\cdot(\nabla^2v\,\nabla v)\\
 &=|\nabla^2v\,\nabla v|^2,
\end{align*}
where we again used the symmetry of $\nabla^2v$.
Combining these identities proves the claimed formula.

We define the operator $\calL$ by
\begin{equation}\label{Def:cal-L}
 \calL\varphi
 :=\partial_t\varphi
   -a(v)w^{\frac{p-2}2}\sfP\fprod\nabla^2\varphi
   -\sfB\cdot\nabla\varphi
\end{equation}
for smooth scalar functions $\varphi$.

Two remarks are in order here. First, $\sfB$ depends on $\nabla^2v$, so $\calL$ is not an operator with a priori controlled coefficients, and no maximum principle for $\calL$ on the whole cylinder is available. We use $\calL$ only \emph{pointwise}, at a maximum point of $\eta^2w$ where $|\nabla v|>0$. This is why half of the negative Hessian term must be retained in the Young step below: it absorbs the $|\nabla^2v|$ coming from $|\sfB|$. Second, $w^{\frac{p-3}2}$ is singular at $\nabla v=0$ for $p<3$; in the localization step it always appears multiplied by $w$, giving $w^{\frac{p-1}2}$ with a positive exponent, so no singularity survives.

With the above notation, equation \eqref{nn:eq:with-drift} becomes
\begin{equation}\label{nn:eq:Lw}
 \calL w
 =-2a(v)w^{\frac{p-2}2}
   \tr\bigl(\sfP(\nabla^2v)^2\bigr)
  +2a'(v)w^{\frac p2}(\sfP\fprod\nabla^2v)
  +2b'(v)w^{\frac{p+2}2}.
\end{equation}
In terms of $w$, the drift estimate
\eqref{nn:eq:drift} reads
\[
 |\sfB|
 \le C(N,p)\Bigl(
   a(v)w^{\frac{p-3}2}|\nabla^2v|
   +b(v)w^{\frac{p-1}{2}}
 \Bigr).
\]

For later use, we record the following properties of $\sfP$.
For every $\xi\in\R^N$,
\[
 \xi\cdot\sfP\xi
 =|\xi|^2
  +(p-2)\left(
    \xi\cdot\frac{\nabla v}{|\nabla v|}
  \right)^2.
\]
Consequently,
\begin{equation}\label{nn:eq:P-bounds}
 \min\{1,p-1\}|\xi|^2
 \le \xi\cdot\sfP\xi
 \le \max\{1,p-1\}|\xi|^2.
\end{equation}
In particular, $\sfP$ is positive definite. Its
eigenvalue in the direction of $\nabla v$ is $p-1$,
and its $N-1$ eigenvalues in the orthogonal
directions are $1$. Hence
\[
 \tr\sfP=N+p-2.
\]
Since $\nabla^2v$ is symmetric, cyclicity of the trace gives
\[
\tr\bigl(\sfP(\nabla^2v)^2\bigr)=\tr\bigl(\nabla^2v\,\sfP\,\nabla^2v\bigr)
 =\sum_{j=1}^N(\nabla^2v\,e_j)\cdot\sfP(\nabla^2v\,e_j),
\]
and applying \eqref{nn:eq:P-bounds} to the vectors $\nabla^2v\,e_j$ yields 
\[
 \min\{1,p-1\}|\nabla^2v|^2
 \le \tr\bigl(\sfP(\nabla^2v)^2\bigr)
 \le \max\{1,p-1\}|\nabla^2v|^2,
\]
because $\sum_j|\nabla^2v\,e_j|^2=|\nabla^2v|^2$.

\subsection{The negative coefficient and the Hessian term}

Applying the Cauchy--Schwarz inequality for the
Frobenius inner product, we obtain
\[
 (\sfP\fprod\nabla^2v)^2
 \le(\tr\sfP)\tr\bigl(\sfP(\nabla^2v)^2\bigr).
\]
Indeed, writing $\sfP^{\frac12}$ for the symmetric
positive definite square root of $\sfP$, we have
\[
 \sfP\fprod\nabla^2v
 =\sfP^{\frac12}\fprod(\sfP^{\frac12}\nabla^2v).
\]
The Frobenius Cauchy--Schwarz inequality therefore
yields
\begin{align*}
 (\sfP\fprod\nabla^2v)^2
 &\le
 |\sfP^{\frac12}|^2\,|\sfP^{\frac12}\nabla^2v|^2=
 (\tr\sfP)\tr\bigl(\sfP(\nabla^2v)^2\bigr).
\end{align*}
Using $\tr\sfP=N+p-2$, we can estimate the mixed
term in \eqref{nn:eq:Lw} as follows:
\begin{align*}
 &2a'(v)w^{\frac{p}2}(\sfP\fprod\nabla^2v)\\
 &\quad\le
 2|a'(v)|w^{\frac{p}2}
 \sqrt{N+p-2}\,
 \sqrt{\tr\bigl(\sfP(\nabla^2v)^2\bigr)}\\
 &\quad=
 2\left[
   \sqrt{a(v)}\,w^{\frac{p-2}4}
   \sqrt{\tr\bigl(\sfP(\nabla^2v)^2\bigr)}
 \right]
 \left[
   \sqrt{N+p-2}\,
   \frac{|a'(v)|}{\sqrt{a(v)}}w^{\frac{p+2}4}
 \right]\\
 &\quad\le
 a(v)w^{\frac{p-2}2}
 \tr\bigl(\sfP(\nabla^2v)^2\bigr)
 +(N+p-2)\frac{(a'(v))^2}{a(v)}w^{\frac{p+2}2},
\end{align*}
where the last step follows from Young's inequality.
Substituting this estimate into \eqref{nn:eq:Lw}
gives
\begin{align}\label{nn:eq:young}
 \calL w
 &\le
 -\underbrace{a(v)w^{\frac{p-2}2}
 \tr\bigl(\sfP(\nabla^2v)^2\bigr)}_{\text{Hessian term}}+
 \left[
   2b'(v)+(N+p-2)\frac{(a'(v))^2}{a(v)}
 \right]w^{\frac{p+2}2}.
\end{align}
We now compute the coefficient in parentheses.
Since
\[
 a(v)=cv^{2-p},
 \qquad
 a'(v)=c(2-p)v^{1-p},
\]
we have
\[
 \frac{(a'(v))^2}{a(v)}
 =
 \frac{c^2(p-2)^2v^{2-2p}}{cv^{2-p}}
 =
 c(p-2)^2v^{-p}.
\]
Furthermore,
\[
 b'(v)=-c(p-1)^2(k-1)v^{-p}.
\]
Recalling $A_*=(N+p-2)(p-2)^2$, we obtain
\begin{align}\label{nn:eq:negative}
 2b'(v)+(N+p-2)\frac{(a'(v))^2}{a(v)}
 &=
 cv^{-p}\bigl[A_*-2(p-1)^2(k-1)\bigr]\notag\\
 &\le
 -\underbrace{c\bigl[A_*+2(p-1)^2\bigr]}_{=\kappa}v^{-p}.
\end{align}
Here we used
\[
 k\ge 2+\frac{A_*}{(p-1)^2},
\]
which implies
\[
 A_*-2(p-1)^2(k-1)
 \le -A_*-2(p-1)^2.
\]
Recalling that $c=k^{p-2}/(p-1)$, we set
\[
 \kappa
 :=\frac{k^{p-2}}{p-1}
   \bigl[A_*+2(p-1)^2\bigr]>0.
\]
Thus, $\kappa$ depends only on $N,p,k$.
We conclude from \eqref{nn:eq:young} that
\begin{equation}\label{nn:eq:coercive}
 \calL w
 \le
 -a(v)w^{\frac{p-2}2}
 \tr\bigl(\sfP(\nabla^2v)^2\bigr)
 -\kappa v^{-p}w^{\frac{p+2}2}.
\end{equation}

\begin{remark}\upshape
Young's inequality separates the negative
Hessian contribution. The remaining term
\[
 -a(v)w^{\frac{p-2}2}
 \tr\bigl(\sfP(\nabla^2v)^2\bigr)
\]
will control the terms involving $|\nabla^2v|$
that arise from the drift when a cutoff function
is introduced.  
\end{remark}

\begin{remark}\upshape
The proof of \eqref{nn:eq:negative} only needs the coefficient $A_*-2(p-1)^2(k-1)$ to be strictly negative, i.e.,\ $k>1+\frac{A_*}{2(p-1)^2}$; then \eqref{nn:eq:coercive} would hold with $\kappa=\frac{k^{p-2}}{p-1}\bigl[2(p-1)^2(k-1)-A_*\bigr]$. The condition $k\ge k_o=2+\frac{A_*}{(p-1)^2}$ is a convenient sufficient choice.
The estimate for $u$ itself is finally obtained with the fixed value $k=k_o(N,p)$, so the dependence of $C$ on $k$ is harmless.
\end{remark}

\subsection{A local gradient estimate via Bernstein's method}\label{nn:sec:bernstein}

\begin{lemma}\label{nn:lem:gradient-estimate}
Let $p>1$ and let $k\ge k_o$, where $k_o$ is
defined in \eqref{nn:eq:k}.
Suppose that $u$ is a smooth nonnegative solution to~\eqref{PDE} in $\Omega_T$. 
Let $r,h>0$, and let
\[
 Q:=B_{2r}(x_o)\times(s-2h,s]
\]
satisfy
\[
 \overline Q\subset
 \{(x,t)\in\Omega_T:u(x,t)>0\}.
\]
Assume that $u\le\sfM$ in $Q$ and set $v=u^{\frac1k}$.
Then
\begin{equation*}
 \|\nabla v\|_{L^\infty(B_r(x_o)\times(s-h,s))}
 \le C\sfM^{\frac1k}\left[\frac1r+\frac1{h^{\frac1p}}\right],
\end{equation*}
where $C>0$ depends only on $N,p,k$.
In particular, $C$ is independent of any positive
lower bound for $u$ in $Q$.
\end{lemma}

\begin{proof}
Throughout the proof, $C$ denotes a positive constant
depending only on $N,p,k$, whose value may change
from line to line.
By assumption, $u>0$ on $\overline Q$, so
$v=u^{\frac1k}$ is smooth in a neighborhood of
$\overline Q$. We do not require $|\nabla v|>0$
throughout $Q$: this condition holds automatically
at any positive maximum of $W=\eta^2|\nabla v|^2$
considered below.
Choose a smooth cutoff function $0\le\eta\le1$
which vanishes on the parabolic boundary of $Q$,
equals one on $B_r(x_o)\times(s-h,s]$, and satisfies
\begin{equation}\label{nn:eq:cutoff}
 |\nabla\eta|\le\frac Cr,
 \qquad
 |\nabla^2\eta|\le\frac C{r^2},
 \qquad
 |\partial_t\eta|\le\frac Ch.
\end{equation} In \eqref{nn:eq:cutoff}, the constants in the gradient, Hessian, and time-derivative bounds can be chosen as $4$, $32\sqrt N$, and $2$, respectively. 

Recall that $w=|\nabla v|^2$ and consider
\[
 W:=\eta^2w=\eta^2|\nabla v|^2
 \qquad\mbox{on \( \overline{B_{2r}}\times[s-2h,s]\).}
 \]
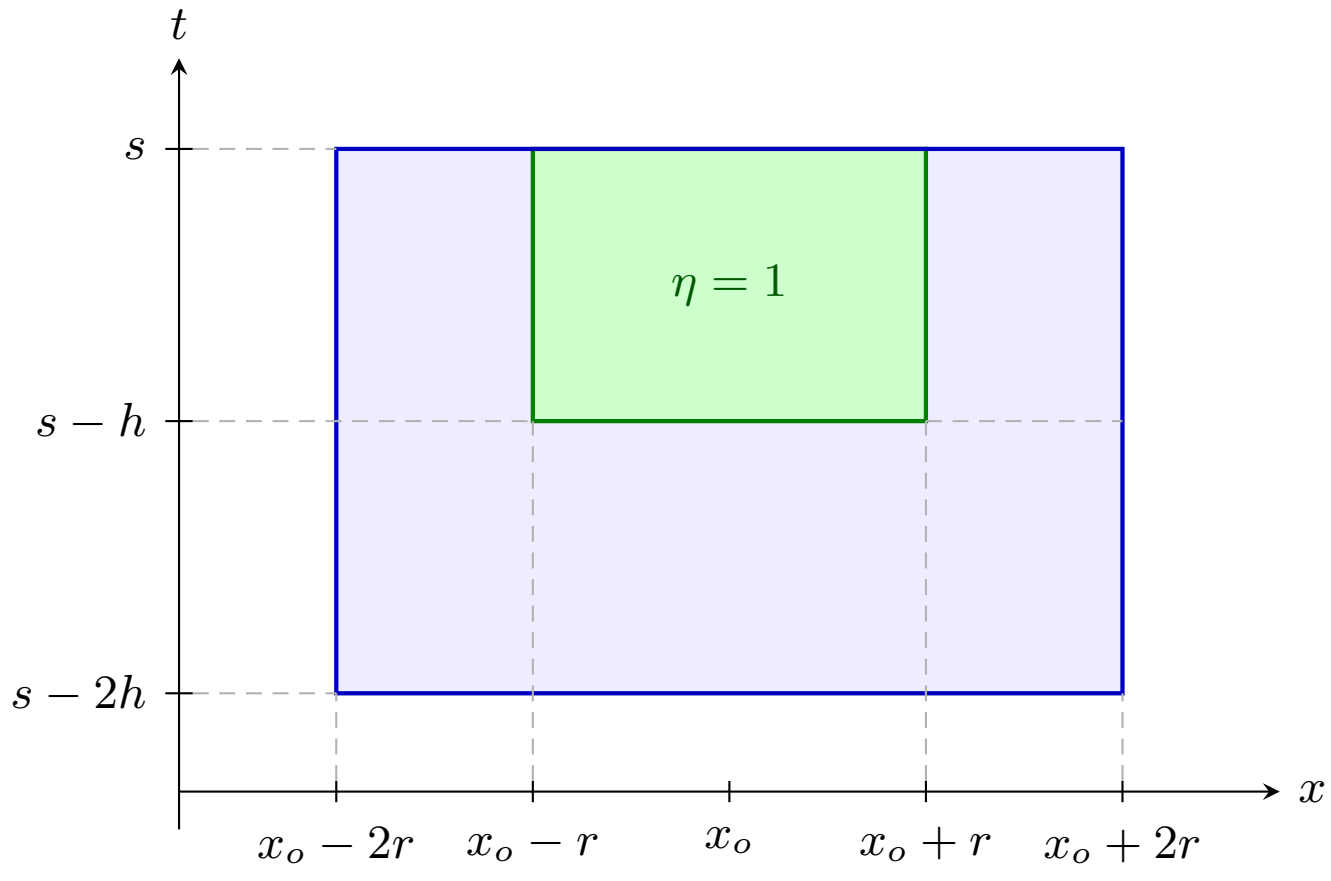
\begin{figure}[htbp]
\centering
\begin{tikzpicture}[
  x=1.3cm,
  y=1cm,
  >=stealth,
  every node/.style={font=\small}
]

\fill[blue!7]
  (-2,0) rectangle (2,3.6);

\filldraw[
  fill=green!20,
  draw=green!50!black,
  thick
]
  (-1,1.8) rectangle (1,3.6);

\draw[blue!75!black,thick]
  (-2,0) rectangle (2,3.6);

\draw[->]
  (-2.8,-0.65) -- (2.8,-0.65)
  node[right] {$x$};

\draw[->]
  (-2.8,-0.9) -- (-2.8,4.2)
  node[above] {$t$};

\draw (-2.87,0) -- (-2.73,0);
\node[left] at (-2.87,0) {$s-2h$};

\draw (-2.87,1.8) -- (-2.73,1.8);
\node[left] at (-2.87,1.8) {$s-h$};

\draw (-2.87,3.6) -- (-2.73,3.6);
\node[left] at (-2.87,3.6) {$s$};

\draw[gray!60,densely dashed]
  (-2.73,0) -- (-2,0);

\draw[gray!60,densely dashed]
  (-2.73,1.8) -- (-1,1.8);

\draw[gray!60,densely dashed]
  (1,1.8) -- (2,1.8);

\draw[gray!60,densely dashed]
  (-2.73,3.6) -- (-2,3.6);

\foreach \x/\label in {
  -2/{x_o-2r},
  -1/{x_o-r},
   0/{x_o},
   1/{x_o+r},
   2/{x_o+2r}
}{
  \draw (\x,-0.72) -- (\x,-0.58);
  \node[below] at (\x,-0.76) {$\label$};
}

\draw[gray!60,densely dashed]
  (-2,-0.58) -- (-2,0);

\draw[gray!60,densely dashed]
  (-1,-0.58) -- (-1,1.8);

\draw[gray!60,densely dashed]
  (1,-0.58) -- (1,1.8);

\draw[gray!60,densely dashed]
  (2,-0.58) -- (2,0);

\node[text=green!35!black] at (0,2.7) {$\eta=1$};

\end{tikzpicture}
\caption{Schematic spatial cross-section of
$Q=B_{2r}(x_o)\times(s-2h,s]$.
The blue rectangle represents
$\overline{B_{2r}(x_o)}\times[s-2h,s]$.
The green region represents
$\overline{B_r(x_o)}\times[s-h,s]$, where $\eta=1$.}
\label{fig:localization-cylinders}
\end{figure}

If the maximum of $W$ is zero, then $W\equiv0$.
Since $\eta=1$ on $B_r(x_o)\times(s-h,s]$, this
implies $\nabla v=0$ there, and the gradient estimate
holds trivially. Otherwise, let $(x_\ast,t_\ast)$ be a point
where $W$ attains a positive maximum. Since $\eta$
vanishes on the parabolic boundary of $Q$, we have
\[
 x_\ast\in B_{2r}(x_o),
 \qquad
 t_\ast\in(s-2h,s],
 \qquad
 \eta(x_\ast,t_\ast)>0,
 \qquad
 w(x_\ast,t_\ast)>0.
\]
The spatial maximum is attained at an interior
point, so $\nabla W=0$ and $\nabla^2W\le0$ at $(x_\ast,t_\ast)$.
If $t_\ast<s$, then $\partial_tW=0$. If $t_\ast=s$,
the maximality of $W$ gives
\[
 \partial_tW(x_\ast,s)
 =\lim_{\tau\downarrow0}
 \frac{W(x_\ast,s)-W(x_\ast,s-\tau)}{\tau}
 \ge0.
\]
Thus, $\partial_tW\ge0$ in either case.
Since $\sfP$ is positive definite and $\nabla^2W$
is negative semidefinite at $(x_\ast,t_\ast)$, we have
\[
 \sfP\fprod\nabla^2W
 =\tr\bigl(\sfP^{1/2}\nabla^2W\,\sfP^{1/2}\bigr)
 \le0.
\]
Since $w(x_\ast,t_\ast)>0$ and $\nabla v$ is continuous,
we have $|\nabla v|>0$ in a neighborhood of
$(x_\ast,t_\ast)$. The identities and inequalities
derived above therefore apply in this neighborhood.
Together with $\nabla W=0$, $\partial_tW\ge0$,
and $a(v)w^{\frac{p-2}2}>0$, the inequality
$\sfP\fprod\nabla^2W\le0$ yields
\begin{equation}\label{nn:eq:maximum}
 \calL W
 =\partial_tW
  -a(v)w^{\frac{p-2}2}\sfP\fprod\nabla^2W
  -\sfB\cdot\nabla W
 \ge0
 \qquad\text{at }(x_\ast,t_\ast).
\end{equation}
Moreover, at $(x_\ast,t_\ast)$ we have
\[
 0=\nabla W
   =\nabla(\eta^2w)
   =2\eta w\nabla\eta+\eta^2\nabla w.
\]
Since $\eta(x_\ast,t_\ast)>0$, we obtain
\begin{equation}\label{nn:eq:max-gradient}
 \nabla w=-2w\frac{\nabla\eta}{\eta}
 \qquad\text{at }(x_\ast,t_\ast).
\end{equation}

Notice that where $\nabla v=0$ the principal coefficient $a\,w^{\frac{p-2}2}\sfP$ vanishes ($p>2$) or blows up ($p<2$), and the differentiated identity is not available; the positivity of $W(x_*,t_*)$ is what keeps the argument away from this set.

To expand $\calL(\eta^2w)$, we first note that,
for smooth scalar functions $f$ and $g$,
\[
 \nabla^2(fg)
 =f\nabla^2g+g\nabla^2f
  +\nabla f\otimes\nabla g
  +\nabla g\otimes\nabla f.
\]
Since $\sfP$ is symmetric, 
we obtain
\[
 \sfP\fprod\nabla^2(fg)
 =f\,\sfP\fprod\nabla^2g+g\,\sfP\fprod\nabla^2f
  +2(\sfP\nabla f)\cdot\nabla g.
\]
The time derivative and the drift term satisfy the
usual product rule. Consequently,
\[
 \calL(fg)
 =f\calL g+g\calL f
  -2a(v)w^{\frac{p-2}2}
   (\sfP\nabla f)\cdot\nabla g.
\]
Applying this identity with $f=\eta^2$ and $g=w$,
we obtain
\[
 \calL(\eta^2w)
 =\eta^2\calL w+w\calL(\eta^2)
  -2a(v)w^{\frac{p-2}2}
   (\sfP\nabla(\eta^2))\cdot\nabla w.
\]
Applying the same product rule to $\eta^2=\eta\eta$
yields
\[
 \calL(\eta^2)
 =2\eta\calL\eta
  -2a(v)w^{\frac{p-2}2}
   (\sfP\nabla\eta)\cdot\nabla\eta.
\]
Substituting this expression and using
$\nabla(\eta^2)=2\eta\nabla\eta$, we find
\begin{align*}
 \calL(\eta^2w)
 &=
 \eta^2\calL w+2\eta w\calL\eta
 -2a(v)w^{\frac{p}2}
   (\sfP\nabla\eta)\cdot\nabla\eta
 -4a(v)\eta w^{\frac{p-2}2}
   (\sfP\nabla\eta)\cdot\nabla w.
\end{align*}
We now evaluate this identity at $(x_\ast,t_\ast)$
and omit the evaluation point for readability.
By \eqref{nn:eq:max-gradient}, we have
\[
 -4a(v)\eta w^{\frac{p-2}2}
   (\sfP\nabla\eta)\cdot\nabla w
 =
 8a(v)w^{\frac{p}2}
   (\sfP\nabla\eta)\cdot\nabla\eta.
\]
Combining this with the preceding identity and using
\eqref{nn:eq:maximum}, we obtain
\begin{align}\label{nn:eq:LW-0}\notag
 0&\le\calL W\\ 
 &=
 \eta^2\calL w
 +2\eta w\Bigl(
   \partial_t\eta
   -a(v)w^{\frac{p-2}{2}}\sfP\fprod\nabla^2\eta
   -\sfB\cdot\nabla\eta
 \Bigr)
 +6a(v)w^{\frac p2}
   (\sfP\nabla\eta)\cdot\nabla\eta.
\end{align}
We apply \eqref{nn:eq:coercive} to the term
$\eta^2\calL w$:
\begin{align*}
 \eta^2\calL w
 &\le
 -a(v)\eta^2w^{\frac{p-2}2}
   \tr\bigl(\sfP(\nabla^2v)^2\bigr)
 -\kappa\eta^2v^{-p}w^{\frac{p+2}2}.
\end{align*}
We estimate the remaining terms separately.
Since $\eta$ and $w$ are nonnegative, the term
containing the time derivative satisfies
\[
 2\eta w\,\partial_t\eta
 \le 2\eta w|\partial_t\eta|.
\]
For the term containing $\nabla^2\eta$, we use
the Cauchy--Schwarz inequality for the Frobenius
inner product: \(|\sfP\fprod\nabla^2\eta|
 \le|\sfP|\,|\nabla^2\eta| \).
Since $\sfP$ is symmetric, its squared Frobenius
norm equals the sum of the squares of its eigenvalues.
Thus,
\[
 |\sfP|^2
 =\tr(\sfP^T\sfP)
 =\tr(\sfP^2)
 =(N-1)+(p-1)^2,
\]
and therefore
\[
 -2a(v)\eta w^{\frac{p}2}\sfP\fprod\nabla^2\eta
 \le
 2\sqrt{N-1+(p-1)^2}\,
 a(v)\eta w^{\frac{p}2}|\nabla^2\eta|.
\]
The constant in this estimate depends only on
$N$ and $p$. For the drift term, the usual Cauchy--Schwarz
inequality gives
\[
 -2\eta w\,\sfB\cdot\nabla\eta
 \le 2\eta w|\sfB|\,|\nabla\eta|.
\]
We leave $|\sfB|$ unchanged at this stage and
estimate it later using \eqref{nn:eq:drift}.
Finally, the upper bound in \eqref{nn:eq:P-bounds}
yields
\[
 (\sfP\nabla\eta)\cdot\nabla\eta
 \le\max\{1,p-1\}|\nabla\eta|^2.
\]
Consequently,
\[
 6a(v)w^{\frac{p}2}(\sfP\nabla\eta)\cdot\nabla\eta
 \le
 6\max\{1,p-1\}\,
 a(v)w^{\frac{p}2}|\nabla\eta|^2.
\]
Here, the constant depends only on $p$.
Substituting these estimates into the preceding
inequality \eqref{nn:eq:LW-0} and moving the two negative terms to
the left-hand side, we obtain
\begin{align}\label{nn:eq:localization-before-drift}
 a(v)\eta^2&w^{\frac{p-2}{2}}
   \tr\bigl(\sfP(\nabla^2v)^2\bigr)
 +\kappa\eta^2v^{-p}w^{\frac{p+2}{2}}\notag\\
 &\le
 C\eta|\partial_t\eta|w
 +Ca(v)w^{\frac p2}
   \bigl(|\nabla\eta|^2+\eta|\nabla^2\eta|\bigr)
 +2\eta w|\sfB|\,|\nabla\eta|.
\end{align}
Here $C$ depends only on $N$ and $p$.
By the drift estimate \eqref{nn:eq:drift},
there exists a constant $C=C(N,p)>0$ such that
\begin{align*}
 2\eta w|\sfB|\,|\nabla\eta|
 &\le
 Ca(v)\eta w^{\frac{p-1}{2}}
 |\nabla^2v|\,|\nabla\eta|
 +Cb(v)\eta w^{\frac{p+1}{2}}|\nabla\eta|.
\end{align*}
After substituting the drift estimate into
\eqref{nn:eq:localization-before-drift}, we absorb
the contribution involving $|\nabla^2v|$ into
the Hessian term on the left-hand side.
Indeed, 
Young's inequality yields
\begin{align*}
 &Ca(v)\eta w^{\frac{p-1}{2}}
   |\nabla^2v|\,|\nabla\eta|\\
    &\qquad=
     \left(
   \sqrt{a(v)}\,\eta w^{\frac{p-2}{4}}|\nabla^2v|
 \right)
 \left(
   C\sqrt{a(v)}\,w^{\frac p4}|\nabla\eta|
 \right)\\
 &\qquad\le
 \varepsilon a(v)\eta^2w^{\frac{p-2}{2}}|\nabla^2v|^2
 +\frac{C^2}{4\varepsilon}
   a(v)w^{\frac p2}|\nabla\eta|^2\\
   &\qquad\le
   \frac{\varepsilon a(v)\eta^2w^{\frac{p-2}2}}{\min\{1,p-1\}} \tr\bigl(\sfP(\nabla^2v)^2\bigr)
   +
   \frac{C^2}{4\varepsilon}
   a(v)w^{\frac p2}|\nabla\eta|^2,
\end{align*}
where we used
\[
 \tr\bigl(\sfP(\nabla^2v)^2\bigr)
 \ge\min\{1,p-1\}|\nabla^2v|^2.
\]
Choosing
\[
 \varepsilon=\tfrac12\min\{1,p-1\},
\]
we obtain
\begin{align*}
 &Ca(v)\eta w^{\frac{p-1}{2}}
   |\nabla^2v|\,|\nabla\eta|\\
 &\qquad\le
 \tfrac12 a(v)\eta^2w^{\frac{p-2}{2}}
 \tr\bigl(\sfP(\nabla^2v)^2\bigr)
 +\frac{C^2}{2\min\{1,p-1\}}
 a(v)w^{\frac p2}|\nabla\eta|^2.
\end{align*}
The first term is absorbed into the left-hand side
of \eqref{nn:eq:localization-before-drift}.
The constant multiplying the second term on the
right-hand side depends only on $N$ and $p$.
Substituting the estimates for the terms on the
right-hand side of
\eqref{nn:eq:localization-before-drift},
absorbing the Hessian contribution into the
left-hand side, and then dropping the remaining
nonnegative Hessian term, we obtain
\begin{align}\label{nn:eq:cutoff-ineq}
 \kappa\eta^2v^{-p}w^{\frac{p+2}{2}}
 &\le C\Bigl[
   \eta|\partial_t\eta|w
   +a(v)w^{\frac p2}
     \bigl(|\nabla\eta|^2+\eta|\nabla^2\eta|\bigr)
   +b(v)\eta w^{\frac{p+1}{2}}|\nabla\eta|
 \Bigr].
\end{align}
Here $C$ depends only on $N$ and $p$.
To estimate the maximum of $W$, we set
\[
 \sfY:=\eta\sqrt w=\sqrt W.
\]
Multiplying \eqref{nn:eq:cutoff-ineq} by
$v^p\eta^p$ and recalling $\sfY=\eta\sqrt w$, we obtain
\begin{align*}
 \kappa
 \underbrace{\eta^{p+2}w^{\frac{p+2}{2}}}_{=\sfY^{p+2}}
 &\le C\Bigl[
   v^p|\partial_t\eta|\,
   \underbrace{\eta^{p+1}w}_{=\eta^{p-1}\sfY^2}\\
 &\qquad
   +a(v)v^p
   \underbrace{\eta^pw^{\frac p2}}_{=\sfY^p}
   \bigl(|\nabla\eta|^2+\eta|\nabla^2\eta|\bigr)\\
 &\qquad
   +b(v)v^p
   \underbrace{\eta^{p+1}w^{\frac{p+1}{2}}}_{=\sfY^{p+1}}
   |\nabla\eta|
 \Bigr].
\end{align*}
Thus,
\begin{align*}
 \kappa\sfY^{p+2}
 &\le C\Bigl[
   v^p\eta^{p-1}|\partial_t\eta|\sfY^2
   +a(v)v^p\sfY^p
     \bigl(|\nabla\eta|^2+\eta|\nabla^2\eta|\bigr)
   +b(v)v^p\sfY^{p+1}|\nabla\eta|
 \Bigr].
\end{align*}
Since $p>1$ and $0\le\eta\le1$, we have
$\eta^{p-1}\le1$. Recalling that
$c=k^{p-2}/(p-1)$, we substitute
$a(v)v^p=cv^2$ and
$b(v)v^p=c(p-1)(k-1)v$.
Using $v\le\sfM^{\frac1k}$ and the derivative bounds in
\eqref{nn:eq:cutoff}, and dividing by $\kappa$,
we obtain
\begin{equation}\label{nn:eq:Y}
 \sfY^{p+2}
 \le C\left[
   \frac{\sfM^{\frac pk}}{h}\sfY^2
   +\frac{\sfM^{\frac 2k}}{r^2}\sfY^p
   +\frac{\sfM^{\frac1k}}r\sfY^{p+1}
 \right],
\end{equation}
where $C=C(N,p,k)>0$ is a new constant.
Young's inequality gives, for every $\varepsilon>0$,
\begin{align*}
 \frac{\sfM^{\frac pk}}{h}\sfY^2
 &\le
 \varepsilon\sfY^{p+2}
 +C(p)\varepsilon^{-\frac2p}
 \frac{\sfM^{\frac{p+2}{k}}}{h^{\frac{p+2}{p}}},\\
 \frac{\sfM^{\frac 2k}}{r^2}\sfY^p
 &\le
 \varepsilon\sfY^{p+2}
 +C(p)\varepsilon^{-\frac p2}
 \frac{\sfM^{\frac{p+2}{k}}}{r^{p+2}},
 \end{align*}
 and
 \begin{align*}
 \frac{\sfM^{\frac1k}}r\sfY^{p+1}
 &\le
 \varepsilon\sfY^{p+2}
 +C(p)\varepsilon^{-(p+1)}
 \frac{\sfM^{\frac{p+2}{k}}}{r^{p+2}}.
\end{align*}
Choosing $\varepsilon=\frac1{6C}$, where $C$ is
the constant in \eqref{nn:eq:Y}, we can absorb
the three contributions
$C\varepsilon\sfY^{p+2}$ into the left-hand side.
We obtain
\[
 \sfY^{p+2}
 \le C \sfM^{\frac{p+2}{k}}
 \left[
   \frac1{h^{\frac{p+2}{p}}}
   +\frac1{r^{p+2}}
 \right],
\]
where $C=C(N,p,k)>0$. Taking the $(p+2)$-th root, we conclude that
\[
 \sfY(x_\ast,t_\ast)
 \le C \sfM^{\frac1k}
\left[\frac1r+\frac1{h^{\frac1p}}
 \right].
\]
Since $(x_\ast,t_\ast)$ is a maximum point of
$W=\sfY^2$ on
$\overline{B_{2r}(x_o)}\times[s-2h,s]$
and $\sfY\ge0$, we have
\[
 \sfY(x,t)\le\sfY(x_\ast,t_\ast)
 \qquad
 \text{for all }(x,t)\in
 \overline{B_{2r}(x_o)}\times[s-2h,s].
\]
In particular, this inequality holds on
$B_r(x_o)\times(s-h,s)$, where $\eta=1$.
Thus, for every $(x,t)$ in this smaller cylinder,
\[
 |\nabla v(x,t)|
 =\eta(x,t)\sqrt{w(x,t)}
 =\sfY(x,t)
 \le\sfY(x_\ast,t_\ast)
 \le C \sfM^{\frac1k}\left[\frac1r+\frac1{h^{\frac1p}}\right].
\]
Taking the supremum over the smaller cylinder yields
\begin{equation*}
 \|\nabla v\|_{L^\infty(B_r(x_o)\times(s-h,s))}
 \le C \sfM^{\frac1k}\left[\frac1r+\frac1{h^{\frac1p}}\right],
\end{equation*}
where the constant $C=C(N,p,k)$ is independent of $s$. 
No positive lower bound for $v$ enters the constant
in this estimate. 
\end{proof}

\begin{remark}[The idea of the gradient estimate]\upshape
The argument begins with the transformation $u=v^k$,
which gives
\[
 v_t=a(v)\Delta_pv+b(v)|\nabla v|^p.
\]
The exponent $k$ is chosen so that, after differentiating
the equation and applying Young's inequality, the
coefficient of the highest power of $|\nabla v|$
is strictly negative.
More precisely, for $w=|\nabla v|^2$ we obtain
\[
 \calL w
 \le
 -a(v)w^{\frac{p-2}{2}}
   \tr\bigl(\sfP(\nabla^2v)^2\bigr)
 -\kappa v^{-p}w^{\frac{p+2}{2}}.
\]
The two negative terms serve different purposes.
The first controls the Hessian, whereas the second
provides the term needed to bound the gradient.
We then localize the argument by introducing
$W=\eta^2w$. At a positive maximum of $W$, we have
\[
 \calL W\ge0,
 \qquad
 \nabla w=-2w\frac{\nabla\eta}{\eta}.
\]
These relations allow us to control the additional
terms arising from the cutoff. In particular, the
negative Hessian term absorbs the Hessian-dependent
contribution of the drift.
After this absorption, the bound $v\le\sfM^{\frac1k}$ and
the cutoff estimates reduce the differential
inequality to the algebraic inequality
\[
 \sfY^{p+2}
 \le C\left[
   \frac{\sfM^{\frac pk}}{h}\sfY^2
   +\frac{\sfM^{\frac 2k}}{r^2}\sfY^p
   +\frac{\sfM^{\frac1k}}{r}\sfY^{p+1}
 \right]
\]
at the maximum point, where $\sfY=\eta\sqrt w$.
All powers of $\sfY$ on the right-hand side are
strictly smaller than $p+2$. Therefore, Young's inequality
gives a bound for $\sfY$, and the fact
that $\eta=1$ on the smaller cylinder yields
\[
 \|\nabla v\|_{L^\infty(B_r(x_o)\times(s-h,s))}
 \le C\sfM^{\frac1k}\left[\frac1r+\frac1{h^{\frac1p}}\right].
\]
This is a localized Bernstein argument (see \cite{Bernstein-1904,Bernstein-1912}), in the sense that a differential
inequality for the squared gradient is combined with
the maximum principle to obtain a gradient bound.
The homogeneity of the equation accounts for the
linear dependence on $\sfM^{\frac1k}$. The constant $C$ depends
only on $N,p,k$, and not on a positive lower bound
for $v$. 
\end{remark}

\section{Gradient estimate for regularized equation}
\label{nn:sec:regularization}

For $\varepsilon>0$, we consider the regularized equation
\begin{equation}\label{nn:eq:reg}
 \partial_t(u_\varepsilon^{p-1})
 =
 \diver\Bigl(
   (\varepsilon^2+|\nabla u_\varepsilon|^2)^{\frac{p-2}{2}}
   \nabla u_\varepsilon
 \Bigr).
\end{equation}
We work with smooth solutions on a cylinder $Q=B_{2r}(x_o)\times(s-2h,s]$ with $\overline Q\subset\Omega_T$, satisfying $0<\sfm\le u_\varepsilon\le\sfM$ in $Q$, where $\sfm,\sfM$ do not depend on $\varepsilon$. All constants below are independent of $\varepsilon$, $\sfm$, $\sfM$, $r$, $h$. The treatment of points where the limiting solution vanishes requires a separate step.
Throughout this section, $k$ satisfies
\eqref{nn:eq:k}.

\subsection{The transformed equation}
The argument follows the same strategy as in the
previous section for smooth positive solutions.
We set
\[
 v_\varepsilon:=u_\varepsilon^{\frac1k},
 \qquad
 w_\varepsilon:=|\nabla v_\varepsilon|^2,
\]
and recall
\[
 c:=\frac{k^{p-2}}{p-1},
 \qquad
 a(v):=cv^{2-p}.
\]
The regularization introduces an additional dependence on $v_\varepsilon$ in the transformed coefficients, which requires a separate calculation. 
Since $u_\varepsilon=v_\varepsilon^k$, we have
\[
 \nabla u_\varepsilon
 =kv_\varepsilon^{k-1}\nabla v_\varepsilon
\]
and
\[
 \nabla^2u_\varepsilon
 =
 kv_\varepsilon^{k-1}\nabla^2v_\varepsilon
 +k(k-1)v_\varepsilon^{k-2}
   \nabla v_\varepsilon\otimes\nabla v_\varepsilon.
\]
Moreover,
\[
 \partial_t(u_\varepsilon^{p-1})
 =
 k(p-1)v_\varepsilon^{k(p-1)-1}
 \partial_t v_\varepsilon.
\]
To expand the regularized divergence, we use
\begin{align*}
 &\diver\Bigl(
   (\varepsilon^2+|\nabla u_\varepsilon|^2)^{\frac{p-2}{2}}
   \nabla u_\varepsilon
 \Bigr)\\
 &\qquad=
 (\varepsilon^2+|\nabla u_\varepsilon|^2)^{\frac{p-2}{2}}
 \left(
   I+(p-2)
   \frac{\nabla u_\varepsilon\otimes\nabla u_\varepsilon}
        {\varepsilon^2+|\nabla u_\varepsilon|^2}
 \right)\fprod\nabla^2u_\varepsilon.
\end{align*}
We introduce
\[
 \sfw_\varepsilon
 :=
 w_\varepsilon
 +\frac{\varepsilon^2}
       {k^2v_\varepsilon^{2(k-1)}},
 \qquad
 \boldsymbol{\theta}_\varepsilon
 :=\frac{w_\varepsilon}{\sfw_\varepsilon}.
\]
With this notation,
\[
 \varepsilon^2+|\nabla u_\varepsilon|^2
 =
 k^2v_\varepsilon^{2(k-1)}\sfw_\varepsilon
\]
and
\[
 \frac{\nabla u_\varepsilon\otimes\nabla u_\varepsilon}
      {\varepsilon^2+|\nabla u_\varepsilon|^2}
 =
 \frac{\nabla v_\varepsilon\otimes\nabla v_\varepsilon}
      {\sfw_\varepsilon}.
\]
Substituting these expressions into
\eqref{nn:eq:reg} and dividing by
$k(p-1)v_\varepsilon^{k(p-1)-1}$ yields
\begin{align*}
 \partial_t v_\varepsilon
 &=
 a(v_\varepsilon)\sfw_\varepsilon^{\frac{p-2}{2}}
 \left(
   I+(p-2)
   \frac{\nabla v_\varepsilon\otimes\nabla v_\varepsilon}
        {\sfw_\varepsilon}
 \right)\fprod\nabla^2v_\varepsilon\\
 &\quad
 +\frac{a(v_\varepsilon)(k-1)}{v_\varepsilon}
   \sfw_\varepsilon^{\frac{p-2}{2}}
 \left(
   I+(p-2)
   \frac{\nabla v_\varepsilon\otimes\nabla v_\varepsilon}
        {\sfw_\varepsilon}
 \right)\fprod
 (\nabla v_\varepsilon\otimes\nabla v_\varepsilon).
\end{align*}
The contraction in the last term satisfies
\begin{align*}
 \left(
   I+(p-2)
   \frac{\nabla v_\varepsilon\otimes\nabla v_\varepsilon}
        {\sfw_\varepsilon}
 \right)\fprod
 (\nabla v_\varepsilon\otimes\nabla v_\varepsilon)
 &=w_\varepsilon
  +(p-2)\frac{w_\varepsilon^2}{\sfw_\varepsilon}\\
 &
 =w_\varepsilon
  \bigl(1+(p-2)\boldsymbol{\theta}_\varepsilon\bigr).
\end{align*}
Consequently,
\begin{align}\label{nn:eq:reg-v}
 \partial_t v_\varepsilon
 &=
 a(v_\varepsilon)\sfw_\varepsilon^{\frac{p-2}{2}}
 \left(
   I+(p-2)
   \frac{\nabla v_\varepsilon\otimes\nabla v_\varepsilon}
        {\sfw_\varepsilon}
 \right)\fprod\nabla^2v_\varepsilon\notag\\
 &\quad
 +\frac{a(v_\varepsilon)(k-1)}{v_\varepsilon}
   \sfw_\varepsilon^{\frac{p-2}{2}}
   w_\varepsilon
   \bigl(1+(p-2)\boldsymbol{\theta}_\varepsilon\bigr).
\end{align}
Notice that the regularization contributes the term
\[
 \frac{\varepsilon^2}
      {k^2v_\varepsilon^{2(k-1)}}
\]
to $\sfw_\varepsilon$. Since this term depends on
$v_\varepsilon$, differentiating the coefficients
in \eqref{nn:eq:reg-v} produces additional terms
that we compute below.

\subsection{The coefficients and their derivatives}

We now work at points where $w_\varepsilon>0$.
Then $0<\boldsymbol{\theta}_\varepsilon<1$ and
\( \sfw_\varepsilon
 =w_\varepsilon/{\boldsymbol{\theta}_\varepsilon}\).
To keep the dependence on the gradient direction
explicit, for $\theta\in(0,1]$ and
$\xi\in\R^N\setminus\{0\}$ we define
\begin{align*}
 \sfC_\theta(\xi)
 &:=
 \theta^{-\frac{p-2}{2}}
 \left(
   I+(p-2)\theta
   \frac{\xi\otimes \xi}{|\xi|^2}
 \right)
\end{align*}
and 
\begin{align*}
 \sfj(\theta)
 &:=
 \theta^{-\frac{p-2}{2}}
 \bigl(1+(p-2)\theta\bigr).
\end{align*}
When differentiating $\sfC_\theta(\xi)$ with respect
to $\theta$, we keep $\xi$ fixed. For independent variables $v>0$ and
$\xi\in\R^N$, we define the principal coefficient
matrix by
\begin{equation*}
    \sfA_\varepsilon(v,\xi)
 :=
 a(v)
 \left(
   |\xi|^2+\frac{\varepsilon^2}{k^2v^{2(k-1)}}
 \right)^{\frac{p-2}{2}}
 \left(
   I+(p-2)
   \frac{\xi\otimes\xi}
        {|\xi|^2+\frac{\varepsilon^2}{k^2v^{2(k-1)}}}
 \right),
\end{equation*}
and the scalar term by
\begin{equation*}
  \sfF_\varepsilon(v,\xi)
 :=
 \frac{a(v)(k-1)}v
 \left(
   |\xi|^2+\frac{\varepsilon^2}{k^2v^{2(k-1)}}
 \right)^{\frac{p-2}{2}}
 |\xi|^2
 \left(
   1+(p-2)
   \frac{|\xi|^2}
        {|\xi|^2+\frac{\varepsilon^2}{k^2v^{2(k-1)}}}
 \right).   
\end{equation*}
Here $v$ and $\xi$ are independent arguments of
the coefficient functions. They do not replace
the notation $v_\varepsilon$ for the transformed
solution. At points where $\nabla v_\varepsilon\ne0$,
evaluating at
$(v,\xi)=(v_\varepsilon,\nabla v_\varepsilon)$
gives the principal coefficient matrix in the form
\begin{equation*}
    \sfA_\varepsilon(v_\varepsilon,\nabla v_\varepsilon)
 =
 a(v_\varepsilon)
 w_\varepsilon^{\frac{p-2}{2}}
\sfC_{\boldsymbol{\theta}_\varepsilon}(\nabla v_\varepsilon),
\end{equation*}
while the scalar term is given by
\begin{equation*}
    \sfF_\varepsilon
    (v_\varepsilon,\nabla v_\varepsilon)
 =
 \frac{a(v_\varepsilon)(k-1)}{v_\varepsilon}
 w_\varepsilon^{\frac p2}\sfj(\boldsymbol{\theta}_\varepsilon).
\end{equation*}
 Therefore, equation \eqref{nn:eq:reg-v} can be re-written as
\begin{equation}\label{nn:eq:reg-v-short}
 \partial_t v_\varepsilon
 =
 \sfA_\varepsilon(v_\varepsilon,\nabla v_\varepsilon)
 \fprod\nabla^2v_\varepsilon
 +\sfF_\varepsilon(v_\varepsilon,\nabla v_\varepsilon).
\end{equation}

To carry out the Bernstein argument, we derive
an equation for $w_\varepsilon=|\nabla v_\varepsilon|^2$
by differentiating \eqref{nn:eq:reg-v-short} with
respect to the spatial variables. The chain rule
then produces derivatives of $\sfA_\varepsilon$
and $\sfF_\varepsilon$ with respect to both their
scalar and vector arguments. We first compute
the partial derivatives with respect to the scalar
argument, keeping the vector argument fixed. For this purpose, we introduce the function
\[
 \boldsymbol{\Theta}_\varepsilon(v,\xi)
 :=
 \frac{|\xi|^2}
      {|\xi|^2+\frac{\varepsilon^2}{k^2v^{2(k-1)}}},
 \qquad v>0,\quad \xi\in\R^N.
\]
The previously defined quantity
$\boldsymbol{\theta}_\varepsilon$ is then given by
\[
 \boldsymbol{\theta}_\varepsilon
  =\boldsymbol{\Theta}_\varepsilon
   (v_\varepsilon,\nabla v_\varepsilon).
\]
Since
\[
 \partial_v
 \left(
   |\xi|^2+\frac{\varepsilon^2}{k^2v^{2(k-1)}}
 \right)
 =
 -\frac{2(k-1)}v
 \frac{\varepsilon^2}{k^2v^{2(k-1)}},
\]
we obtain
\begin{align}\label{nn:eq:Theta-v}
\partial_v\boldsymbol{\Theta}_\varepsilon(v,\xi)
 &=
 \frac{2(k-1)}v
 \frac{|\xi|^2}
      {|\xi|^2+\frac{\varepsilon^2}{k^2v^{2(k-1)}}}
 \frac{\frac{\varepsilon^2}{k^2v^{2(k-1)}}}
      {|\xi|^2+\frac{\varepsilon^2}{k^2v^{2(k-1)}}}
 \notag\\
 &=
 \frac{2(k-1)}v
 \boldsymbol{\Theta}_\varepsilon(v,\xi)
 \bigl(1-\boldsymbol{\Theta}_\varepsilon(v,\xi)\bigr).
\end{align}
We also compute the derivative of $\sfC_\theta(\xi)$
with respect to $\theta$, keeping $\xi\ne0$ fixed.
From its definition, the product rule gives
\begin{align}\label{nn:eq:C-theta}
 \partial_\theta\sfC_\theta(\xi)
 &=
 -\tfrac12 (p-2)\theta^{-\frac p2}
 \left(
   I+(p-2)\theta
   \frac{\xi\otimes\xi}{|\xi|^2}
 \right)
 +(p-2)\theta^{-\frac{p-2}{2}}
\frac{\xi\otimes\xi}{|\xi|^2}
 \notag\\
 &=
 -\tfrac12 (p-2)\theta^{-\frac p2}I
 +\tfrac12 (p-2)(4-p)
 \theta^{-\frac{p-2}{2}}
  \frac{\xi\otimes\xi}{|\xi|^2}.
\end{align}
Similarly, differentiating $\sfj(\theta)$, we obtain
\begin{align}\label{nn:eq:j-theta}
 \sfj'(\theta)
 &=
 -\tfrac12 (p-2)\theta^{-\frac p2}
 +\tfrac12(p-2)(4-p)
  \theta^{-\frac{p-2}{2}}.
\end{align}
With these derivatives at hand, we first compute
$\partial_v\sfA_\varepsilon(v,\xi)$, keeping the
vector argument $\xi$ fixed.
For $v>0$ and $\xi\ne0$, we have
\[
 \sfA_\varepsilon(v,\xi)
 =
 a(v)|\xi|^{p-2}
 \sfC_{\boldsymbol{\Theta}_\varepsilon(v,\xi)}(\xi).
\]
The product and chain rules give
\begin{align*}
 \partial_v\sfA_\varepsilon(v,\xi)
 &=
 a'(v)|\xi|^{p-2}
 \sfC_{\boldsymbol{\Theta}_\varepsilon(v,\xi)}(\xi)
 +a(v)|\xi|^{p-2} 
   \partial_\theta\sfC_\theta(\xi)
 \big|_{\theta=\boldsymbol{\Theta}_\varepsilon(v,\xi)}
 \partial_v\boldsymbol{\Theta}_\varepsilon(v,\xi).
\end{align*}
Using $a'(v)=(2-p)a(v)/v$ and
\eqref{nn:eq:Theta-v}, and factoring out
$\frac{a(v)}v|\xi|^{p-2}$, we obtain
\begin{align*}
 \partial_v\sfA_\varepsilon(v,\xi)
 &=
 \frac{a(v)}v|\xi|^{p-2}
 \Biggl[
   (2-p)\sfC_{\boldsymbol{\Theta}_\varepsilon(v,\xi)}(\xi)\\
 &\hspace{23mm}
   +2(k-1)\boldsymbol{\Theta}_\varepsilon(v,\xi)
   \bigl(1-\boldsymbol{\Theta}_\varepsilon(v,\xi)\bigr)
     \partial_\theta\sfC_\theta(\xi)
   \big|_{\theta=\boldsymbol{\Theta}_\varepsilon(v,\xi)}
 \Biggr].
\end{align*}
The preceding calculation treats $v$ and $\xi$ as independent variables. To use the resulting coefficient derivative in the equation for $v_\varepsilon$, we now evaluate it at $(v,\xi)=(v_\varepsilon,\nabla v_\varepsilon)$.
Here $(\partial_v\sfA_\varepsilon)(v_\varepsilon,\nabla v_\varepsilon)$ denotes the partial derivative with respect to the scalar argument, computed with $\xi$ fixed and then evaluated along the solution. Since $\boldsymbol\Theta_\varepsilon(v_\varepsilon,\nabla v_\varepsilon)=\boldsymbol\theta_\varepsilon$ and $|\nabla v_\varepsilon|^{p-2}=w_\varepsilon^{\frac{p-2}2}$, we obtain
\begin{align*}
 &(\partial_v\sfA_\varepsilon)
   (v_\varepsilon,\nabla v_\varepsilon)\\
 &\quad=
 \frac{a(v_\varepsilon)}{v_\varepsilon}
 w_\varepsilon^{\frac{p-2}{2}}
 \Biggl[
   (2-p)
   \sfC_{\boldsymbol{\theta}_\varepsilon}
   (\nabla v_\varepsilon)
   +2(k-1)\boldsymbol{\theta}_\varepsilon
   (1-\boldsymbol{\theta}_\varepsilon)
     \partial_\theta\sfC_\theta(\nabla v_\varepsilon)
   \big|_{\theta=\boldsymbol{\theta}_\varepsilon}
 \Biggr].
\end{align*}
We next differentiate the scalar term
\[
 \sfF_\varepsilon(v,\xi)
 =
 (k-1)\frac{a(v)}v|\xi|^p
 \sfj\bigl(\boldsymbol{\Theta}_\varepsilon(v,\xi)\bigr)
\]
with respect to $v$, keeping $\xi\ne0$ fixed.
First, since $a'(v)=(2-p)a(v)/v$, we have
\begin{align*}
 \frac{\mathrm d}{\mathrm dv}
 \left(\frac{a(v)}v\right)
 &=
 (1-p)\frac{a(v)}{v^2}.
\end{align*}
The product and chain rules give
\begin{align*}
 \partial_v&\sfF_\varepsilon(v,\xi)\\
 &=
 (k-1)(1-p)\frac{a(v)}{v^2}|\xi|^p
 \sfj\bigl(\boldsymbol{\Theta}_\varepsilon(v,\xi)\bigr)
 +(k-1)\frac{a(v)}v|\xi|^p
 \sfj'\bigl(\boldsymbol{\Theta}_\varepsilon(v,\xi)\bigr)
 \partial_v\boldsymbol{\Theta}_\varepsilon(v,\xi)\\
 &=
 \frac{a(v)(k-1)}{v^2}|\xi|^p
 \Bigl[
   (1-p)\sfj\bigl(\boldsymbol{\Theta}_\varepsilon(v,\xi)\bigr)
   +v\,\sfj'\bigl(\boldsymbol{\Theta}_\varepsilon(v,\xi)\bigr)
   \partial_v\boldsymbol{\Theta}_\varepsilon(v,\xi)
 \Bigr].
\end{align*}
Substituting \eqref{nn:eq:Theta-v}, we obtain
\begin{align*}
 \partial_v\sfF_\varepsilon(v,\xi)
 &=
 \frac{a(v)(k-1)}{v^2}|\xi|^p
 \Bigl[
   (1-p)
   \sfj\bigl(\boldsymbol{\Theta}_\varepsilon(v,\xi)\bigr)\\
 &\hspace{30mm}
   +2(k-1)\boldsymbol{\Theta}_\varepsilon(v,\xi)
   \bigl(1-\boldsymbol{\Theta}_\varepsilon(v,\xi)\bigr)
   \sfj'\bigl(\boldsymbol{\Theta}_\varepsilon(v,\xi)\bigr)
 \Bigr].
\end{align*}
Evaluating along the solution as before, we obtain
\begin{align*}
 (\partial_v\sfF_\varepsilon)
   (v_\varepsilon,\nabla v_\varepsilon)
 &=
 \frac{a(v_\varepsilon)(k-1)}{v_\varepsilon^2}
 w_\varepsilon^{\frac p2}
 \Bigl[
   (1-p)\sfj(\boldsymbol{\theta}_\varepsilon)
   +2(k-1)\boldsymbol{\theta}_\varepsilon
   (1-\boldsymbol{\theta}_\varepsilon)
   \sfj'(\boldsymbol{\theta}_\varepsilon)
 \Bigr].
\end{align*} 
To write the preceding formulas for
$(\partial_v\sfA_\varepsilon)
(v_\varepsilon,\nabla v_\varepsilon)$ and
$(\partial_v\sfF_\varepsilon)
(v_\varepsilon,\nabla v_\varepsilon)$ more compactly,
we introduce, for $\theta\in(0,1]$ and $\xi\ne0$,
\begin{equation}\label{nn:eq:reg-D}
 \sfD_\theta(\xi)
 :=
 (2-p)\sfC_\theta(\xi)
 +2(k-1)\theta(1-\theta)
   \partial_\theta\sfC_\theta(\xi),
\end{equation}
and
\begin{equation*}
 \sfJ(\theta)
 :=
 (1-p)\sfj(\theta)
 +2(k-1)\theta(1-\theta)\sfj'(\theta).
\end{equation*}
Thus,
\begin{equation}\label{nn:eq:reg-A-v}
 (\partial_v\sfA_\varepsilon)
 (v_\varepsilon,\nabla v_\varepsilon)
 =
 \frac{a(v_\varepsilon)}{v_\varepsilon}
 w_\varepsilon^{\frac{p-2}{2}}
 \sfD_{\boldsymbol{\theta}_\varepsilon}
 (\nabla v_\varepsilon),
\end{equation}
and
\begin{equation}\label{nn:eq:reg-F-v}
 (\partial_v\sfF_\varepsilon)
 (v_\varepsilon,\nabla v_\varepsilon)
 =
 \frac{a(v_\varepsilon)(k-1)}{v_\varepsilon^2}
 w_\varepsilon^{\frac p2}
 \sfJ(\boldsymbol{\theta}_\varepsilon).
\end{equation}

\subsection{The differentiated equation}
For $\ell=1,\ldots,N$, we define the components
of the drift by
\begin{equation}\label{nn:eq:reg-B}
 (\sfB_\varepsilon)_\ell
 :=
 (\partial_{\xi_\ell}\sfA_\varepsilon)
 (v_\varepsilon,\nabla v_\varepsilon)
 \fprod\nabla^2v_\varepsilon
 +(\partial_{\xi_\ell}\sfF_\varepsilon)
 (v_\varepsilon,\nabla v_\varepsilon),
\end{equation}
and introduce the operator
\begin{equation}\label{nn:eq:reg-operator}
 \calL_\varepsilon\varphi
 :=
 \partial_t\varphi
 -\sfA_\varepsilon(v_\varepsilon,\nabla v_\varepsilon)
  \fprod\nabla^2\varphi
 -\sfB_\varepsilon\cdot\nabla\varphi
\end{equation}
for smooth scalar functions $\varphi$.
We differentiate \eqref{nn:eq:reg-v-short} with
respect to $x_\ell$, multiply by
$2\partial_{x_\ell}v_\varepsilon$, and sum over
$\ell=1,\ldots,N$. In the following calculation,
all coefficient functions and their partial
derivatives are evaluated at
$(v_\varepsilon,\nabla v_\varepsilon)$.
The chain rule, the definition of $\sfB_\varepsilon$,
and the commutation of mixed spatial derivatives give
\begin{align*}
 \partial_t(\partial_{x_\ell}v_\varepsilon)
 &=
 \sfA_\varepsilon\fprod
 \nabla^2(\partial_{x_\ell}v_\varepsilon)
 +(\partial_v\sfA_\varepsilon)\fprod
   \nabla^2v_\varepsilon\,
   \partial_{x_\ell}v_\varepsilon
    +(\partial_v\sfF_\varepsilon)
   \partial_{x_\ell}v_\varepsilon\\
 &\quad
   +\sum_{j=1}^N
 \bigl((\partial_{\xi_j}\sfA_\varepsilon)\fprod
       \nabla^2v_\varepsilon\bigr)
 \partial_{x_\ell}\partial_{x_j}v_\varepsilon
    +\sum_{j=1}^N
   (\partial_{\xi_j}\sfF_\varepsilon)
   \partial_{x_\ell}\partial_{x_j}v_\varepsilon\\
   &=
 \sfA_\varepsilon\fprod
 \nabla^2(\partial_{x_\ell}v_\varepsilon)
 +\Big[ (\partial_v\sfA_\varepsilon)\fprod
   \nabla^2v_\varepsilon
    +(\partial_v\sfF_\varepsilon)\Big]
   \partial_{x_\ell}v_\varepsilon\\
 &\quad
   +\sum_{j=1}^N \Big[\underbrace{
 \bigl((\partial_{\xi_j}\sfA_\varepsilon)\fprod
       \nabla^2v_\varepsilon\bigr)
    +
   (\partial_{\xi_j}\sfF_\varepsilon)}_{= (\sfB_\varepsilon)_j}\Big]
   \underbrace{\partial_{x_\ell}\partial_{x_j}v_\varepsilon}_{=\partial_{x_j}\partial_{x_\ell}v_\varepsilon}\\
   &=
    \sfA_\varepsilon\fprod
 \nabla^2(\partial_{x_\ell}v_\varepsilon) 
 +
  \sfB_\varepsilon\cdot
   \nabla(\partial_{x_\ell}v_\varepsilon)
   +
   \Bigl[
   (\partial_v\sfA_\varepsilon)\fprod\nabla^2v_\varepsilon
   +\partial_v\sfF_\varepsilon
 \Bigr]\partial_{x_\ell}v_\varepsilon.
\end{align*}
Multiplying by $2\partial_{x_\ell}v_\varepsilon$
and summing over $\ell$, we obtain
\begin{align*}
 \partial_t w_\varepsilon
 &=
 2\sum_{\ell=1}^N
 (\partial_{x_\ell}v_\varepsilon)
 \sfA_\varepsilon\fprod
 \nabla^2(\partial_{x_\ell}v_\varepsilon)
 +\sfB_\varepsilon\cdot\nabla w_\varepsilon
 +2w_\varepsilon
 \Bigl(
   (\partial_v\sfA_\varepsilon)\fprod\nabla^2v_\varepsilon
   +\partial_v\sfF_\varepsilon
 \Bigr).
\end{align*}
To rewrite the first term, we use
\[
 \nabla^2w_\varepsilon
 =
 2\sum_{\ell=1}^N
 \nabla(\partial_{x_\ell}v_\varepsilon)
 \otimes\nabla(\partial_{x_\ell}v_\varepsilon)
 +
 2\sum_{\ell=1}^N
 (\partial_{x_\ell}v_\varepsilon)
 \nabla^2(\partial_{x_\ell}v_\varepsilon).
\]
Let $e_1,\ldots,e_N$ denote the standard basis of
$\R^N$. Expanding the gradients in this basis and
rearranging the sums, we obtain
\begin{align*}
 \sum_{\ell=1}^N
 \nabla(\partial_{x_\ell}v_\varepsilon)
 \otimes\nabla(\partial_{x_\ell}v_\varepsilon)
 &=
 \sum_{\ell=1}^N
 \left(
   \sum_{i=1}^N
   \partial_{x_i}\partial_{x_\ell}v_\varepsilon\,e_i
 \right)
 \otimes
 \left(
   \sum_{j=1}^N
   \partial_{x_j}\partial_{x_\ell}v_\varepsilon\,e_j
 \right)\\
 &=
 \sum_{i,j=1}^N
 \left(
   \sum_{\ell=1}^N
   (\partial_{x_i}\partial_{x_\ell}v_\varepsilon)
   (\partial_{x_j}\partial_{x_\ell}v_\varepsilon)
 \right)e_i\otimes e_j\\
 &=
 \sum_{i,j=1}^N
 \left(
   \sum_{\ell=1}^N
   (\partial_{x_i}\partial_{x_\ell}v_\varepsilon)
   (\partial_{x_\ell}\partial_{x_j}v_\varepsilon)
 \right)e_i\otimes e_j\\
 &=
 (\nabla^2v_\varepsilon)^2.
\end{align*}
In the penultimate step, we used the symmetry
of $\nabla^2v_\varepsilon$.
Here $(\nabla^2v_\varepsilon)^2$ denotes the matrix
product of the Hessian with itself.
Rearranging the identity for $\nabla^2w_\varepsilon$,
we obtain
\[
 2\sum_{\ell=1}^N
 (\partial_{x_\ell}v_\varepsilon)
 \nabla^2(\partial_{x_\ell}v_\varepsilon)
 =
 \nabla^2w_\varepsilon
 -2(\nabla^2v_\varepsilon)^2.
\]
Taking the Frobenius inner product with
$\sfA_\varepsilon$ and using its linearity gives
\begin{align*}
 2\sum_{\ell=1}^N
 (\partial_{x_\ell}v_\varepsilon)
 \sfA_\varepsilon\fprod
 \nabla^2(\partial_{x_\ell}v_\varepsilon)
 &=
 \sfA_\varepsilon\fprod\nabla^2w_\varepsilon
 -2\sfA_\varepsilon\fprod(\nabla^2v_\varepsilon)^2\\
 &=
 \sfA_\varepsilon\fprod\nabla^2w_\varepsilon
 -2\tr\bigl(
   \sfA_\varepsilon(\nabla^2v_\varepsilon)^2
 \bigr).
\end{align*}
In the last step, we used the definition of the
Frobenius inner product and the symmetry of
$\sfA_\varepsilon$.
Substituting the preceding identity into the
equation for $\partial_t w_\varepsilon$ and
subtracting
$\sfA_\varepsilon\fprod\nabla^2w_\varepsilon$
and $\sfB_\varepsilon\cdot\nabla w_\varepsilon$
from both sides gives
\begin{align*}
 &\underbrace{
   \partial_t w_\varepsilon
   -\sfA_\varepsilon\fprod\nabla^2w_\varepsilon
   -\sfB_\varepsilon\cdot\nabla w_\varepsilon
 }_{=\calL_\varepsilon w_\varepsilon}\\
 &\qquad=
 -2\tr\bigl(
   \sfA_\varepsilon(\nabla^2v_\varepsilon)^2
 \bigr)
 +2w_\varepsilon
   (\partial_v\sfA_\varepsilon)\fprod\nabla^2v_\varepsilon
 +2w_\varepsilon\,\partial_v\sfF_\varepsilon.
\end{align*}
Here $\sfA_\varepsilon$,
$\partial_v\sfA_\varepsilon$, and
$\partial_v\sfF_\varepsilon$ are evaluated at
$(v_\varepsilon,\nabla v_\varepsilon)$, while
$\sfB_\varepsilon$ is the drift defined in
\eqref{nn:eq:reg-B}.
All quantities are evaluated at the same
space-time point.
The partial derivatives with respect to $v$
are taken with the gradient argument fixed.
We now substitute
\[
 \sfA_\varepsilon(v_\varepsilon,\nabla v_\varepsilon)
 =
 a(v_\varepsilon)w_\varepsilon^{\frac{p-2}{2}}
 \sfC_{\boldsymbol{\theta}_\varepsilon}
 (\nabla v_\varepsilon),
\]
together with the expressions for
$\partial_v\sfA_\varepsilon$ and
$\partial_v\sfF_\varepsilon$ from
\eqref{nn:eq:reg-A-v} and \eqref{nn:eq:reg-F-v}.
Combining the powers of $w_\varepsilon$, we obtain
\begin{align}\label{nn:eq:reg-differentiated}
 \calL_\varepsilon w_\varepsilon
 &=
 -2a(v_\varepsilon)w_\varepsilon^{\frac{p-2}{2}}
 \tr\Bigl(
   \sfC_{\boldsymbol{\theta}_\varepsilon}
   (\nabla v_\varepsilon)
   (\nabla^2v_\varepsilon)^2
 \Bigr)
 +2\frac{a(v_\varepsilon)}{v_\varepsilon}
   w_\varepsilon^{\frac p2}
   \Bigl(
     \sfD_{\boldsymbol{\theta}_\varepsilon}
     (\nabla v_\varepsilon)
     \fprod\nabla^2v_\varepsilon
   \Bigr)\notag\\
 &\quad
 +2\frac{a(v_\varepsilon)(k-1)}{v_\varepsilon^2}
   w_\varepsilon^{\frac{p+2}{2}}
   \sfJ(\boldsymbol{\theta}_\varepsilon).
\end{align}

For $\varepsilon>0$, $\sfA_\varepsilon$ and $\sfF_\varepsilon$ are smooth in $(v,\xi)$ even at $\xi=0$, so the differentiated equation holds at every point (with $u_\varepsilon\ge\sfm>0$), not only where $\nabla v_\varepsilon\ne0$; only the representation through $\sfC_\theta,\sfD_\theta,\sfJ$ requires $w_\varepsilon>0$. This is the difference from the previous section, where the coefficients degenerate or blow up at $\nabla v=0$.  

On the other hand, when $\varepsilon=0$, we have $\boldsymbol\theta=1$; in this case $\sfC_1=\sfP$, $\sfD_1=(2-p)\sfP$, and $\sfJ(1)=-(p-1)^2$. Since $a'=(2-p)a/v$ and $b'=-c(p-1)^2(k-1)v^{-p}$, the middle term of \eqref{nn:eq:reg-differentiated} becomes $2a'w^{\frac p2}\,\sfP\fprod\nabla^2v$ and the last one $2b'w^{\frac{p+2}2}$. Therefore, \eqref{nn:eq:reg-differentiated} reduces to \eqref{nn:eq:Lw}.
\subsection{Estimating the term linear in the Hessian}We next estimate the second term on the right-hand
side of \eqref{nn:eq:reg-differentiated}, which is
linear in $\nabla^2v_\varepsilon$ and has no
prescribed sign. Using the Cauchy--Schwarz and
Young inequalities, we absorb its Hessian-dependent
contribution into the negative first term, retaining
half of that term for the localization argument.
Recall that, for $\theta\in(0,1]$ and $\xi\ne0$,
\[
 \sfC_\theta(\xi)
 =
 \theta^{-\frac{p-2}{2}}
 \left(
   I+(p-2)\theta
   \frac{\xi\otimes\xi}{|\xi|^2}
 \right).
\]
If $\zeta\perp\xi$, then \( \sfC_\theta(\xi)\zeta
 =\theta^{-\frac{p-2}{2}}\zeta \),
whereas
\[
 \sfC_\theta(\xi)\frac{\xi}{|\xi|}
 =
 \theta^{-\frac{p-2}{2}}
 \bigl(1+(p-2)\theta\bigr)\frac{\xi}{|\xi|}.
\]
Thus, its eigenvalues are
$\theta^{-\frac{p-2}{2}}$ in directions orthogonal
to $\xi$ and
$\theta^{-\frac{p-2}{2}}(1+(p-2)\theta)$
in the direction of $\xi$.
Since
\[
 1+(p-2)\theta\ge\min\{1,p-1\}>0,
\]
the symmetric matrix $\sfC_\theta(\xi)$ is positive
definite and hence invertible. Its positive definite
square root $\sfC_\theta(\xi)^{\frac12}$ and the
inverse $\sfC_\theta(\xi)^{-\frac12}$ are therefore
well defined.
For every symmetric matrix $H\in\R^{N\times N}$,
inserting these factors and applying the
Cauchy--Schwarz inequality for the Frobenius
inner product gives
\begin{align}\label{nn:eq:reg-matrix-cauchy}
 \bigl(\sfD_\theta(\xi)\fprod H\bigr)^2
 &=
 \Bigl[
   \bigl(\sfC_\theta(\xi)^{-\frac12}\sfD_\theta(\xi)\bigr)
   \fprod
   \bigl(\sfC_\theta(\xi)^{\frac12}H\bigr)
 \Bigr]^2\notag\\
 &\le
 \left|\sfC_\theta(\xi)^{-\frac12}\sfD_\theta(\xi)\right|^2
 \left|\sfC_\theta(\xi)^{\frac12}H\right|^2\notag\\
 &=
 \tr\bigl(\sfC_\theta(\xi)^{-1}\sfD_\theta(\xi)^2\bigr)
 \tr\bigl(\sfC_\theta(\xi)H^2\bigr).
\end{align}
To justify the last equality, we use the definition
of the Frobenius norm, the symmetry of the matrices,
and the cyclic property of the trace.%
\footnote{For square matrices of the same size,
$\tr(ABC)=\tr(BCA)=\tr(CAB)$.
More generally, the trace of a product is invariant
under cyclic permutations of its factors, but not
under arbitrary permutations.}
For the first factor, these properties give
\begin{align*}
 &\left|
   \sfC_\theta(\xi)^{-\frac12}\sfD_\theta(\xi)
 \right|^2\\
 &\qquad=
 \tr\Bigl(
   \bigl(\sfC_\theta(\xi)^{-\frac12}
         \sfD_\theta(\xi)\bigr)^T
   \sfC_\theta(\xi)^{-\frac12}\sfD_\theta(\xi)
 \Bigr)\\
 &\qquad=
 \tr\Bigl(
   \sfD_\theta(\xi)\sfC_\theta(\xi)^{-1}
   \sfD_\theta(\xi)
 \Bigr)\\
 &\qquad=
 \tr\bigl(
   \sfC_\theta(\xi)^{-1}\sfD_\theta(\xi)^2
 \bigr).
\end{align*}
Similarly, for the second factor,
\begin{align*}
 \left|\sfC_\theta(\xi)^{\frac12}H\right|^2
 &=
 \tr\Bigl(
   \bigl(\sfC_\theta(\xi)^{\frac12}H\bigr)^T
   \sfC_\theta(\xi)^{\frac12}H
 \Bigr)\\
 &=
 \tr\bigl(H\sfC_\theta(\xi)H\bigr)\\
 &=
 \tr\bigl(\sfC_\theta(\xi)H^2\bigr).
\end{align*}
Applying \eqref{nn:eq:reg-matrix-cauchy} with
$\theta=\boldsymbol{\theta}_\varepsilon$,
$\xi=\nabla v_\varepsilon$, and
$H=\nabla^2v_\varepsilon$, followed by Young's
inequality, yields
\begin{align*}
 &2\frac{a(v_\varepsilon)}{v_\varepsilon}
   w_\varepsilon^{\frac p2}
   \Bigl(
     \sfD_{\boldsymbol{\theta}_\varepsilon}
     (\nabla v_\varepsilon)
     \fprod\nabla^2v_\varepsilon
   \Bigr)\\
 &\quad\le
 2\frac{a(v_\varepsilon)}{v_\varepsilon}
 w_\varepsilon^{\frac p2}
 \sqrt{
   \tr\Bigl(
     \sfC_{\boldsymbol{\theta}_\varepsilon}
     (\nabla v_\varepsilon)^{-1}
     \sfD_{\boldsymbol{\theta}_\varepsilon}
     (\nabla v_\varepsilon)^2
   \Bigr)
 }
 \sqrt{
   \tr\Bigl(
     \sfC_{\boldsymbol{\theta}_\varepsilon}
     (\nabla v_\varepsilon)
     (\nabla^2v_\varepsilon)^2
   \Bigr)
 }\\
 &\quad\le
 a(v_\varepsilon)w_\varepsilon^{\frac{p-2}{2}}
 \tr\Bigl(
   \sfC_{\boldsymbol{\theta}_\varepsilon}
   (\nabla v_\varepsilon)
   (\nabla^2v_\varepsilon)^2
 \Bigr)
 +\frac{a(v_\varepsilon)}{v_\varepsilon^2}
 w_\varepsilon^{\frac{p+2}{2}}
 \tr\Bigl(
   \sfC_{\boldsymbol{\theta}_\varepsilon}
   (\nabla v_\varepsilon)^{-1}
   \sfD_{\boldsymbol{\theta}_\varepsilon}
   (\nabla v_\varepsilon)^2
 \Bigr).
\end{align*}
Consequently,
\begin{align}\label{nn:eq:reg-coercivity}
 \calL_\varepsilon w_\varepsilon
 &\le
 -a(v_\varepsilon)w_\varepsilon^{\frac{p-2}{2}}
 \tr\Bigl(
   \sfC_{\boldsymbol{\theta}_\varepsilon}
   (\nabla v_\varepsilon)
   (\nabla^2v_\varepsilon)^2
 \Bigr)\notag\\
 &\quad
 +\frac{a(v_\varepsilon)}{v_\varepsilon^2}
 w_\varepsilon^{\frac{p+2}{2}}
 \Bigl[
   \tr\Bigl(
     \sfC_{\boldsymbol{\theta}_\varepsilon}
     (\nabla v_\varepsilon)^{-1}
     \sfD_{\boldsymbol{\theta}_\varepsilon}
     (\nabla v_\varepsilon)^2
   \Bigr)
   +2(k-1)\sfJ(\boldsymbol{\theta}_\varepsilon)
 \Bigr].
\end{align}

\subsection{A coercive inequality for large gradients}

For each $\xi\ne0$, the matrix $\sfC_\theta(\xi)$
is a linear combination of $I$ and
\( \xi\otimes \xi/|\xi|^2\).
By \eqref{nn:eq:C-theta}, its derivative
$\partial_\theta\sfC_\theta(\xi)$ has the same
structure. Hence, by \eqref{nn:eq:reg-D},
so does $\sfD_\theta(\xi)$.
Their eigenvalues depend only on $\theta,p,k$,
and not on the direction of $\xi$.%
\footnote{A matrix of the form
$\lambda I+\mu\,\frac{\xi}{|\xi|}\otimes\frac{\xi}{|\xi|}$
has eigenvalue $\lambda+\mu$ in the direction of
$\xi$ and eigenvalue $\lambda$ in each orthogonal
direction. For $\sfC_\theta(\xi)$ and
$\sfD_\theta(\xi)$, the scalar coefficients
$\lambda$ and $\mu$ depend only on $\theta,p,k$.} Since $\sfC_\theta(\xi)$ and $\sfD_\theta(\xi)$
have a common orthonormal eigenbasis, the trace
of $\sfC_\theta(\xi)^{-1}\sfD_\theta(\xi)^2$
is the sum of the squared eigenvalues of
$\sfD_\theta(\xi)$ divided by the corresponding
eigenvalues of $\sfC_\theta(\xi)$.
Consequently,
\begin{equation}\label{eq:trace+J}
\tr\bigl(\sfC_\theta(\xi)^{-1}\sfD_\theta(\xi)^2\bigr)
 +2(k-1)\sfJ(\theta)
\end{equation}
depends only on $\theta,N,p,k$. The dependence on $N$ arises from the multiplicity
$N-1$ of the eigenvalues corresponding to directions
orthogonal to $\xi$.

For each fixed $\xi\ne0$, the functions
$\sfC_\theta(\xi)$, $\sfD_\theta(\xi)$,
$\sfj(\theta)$, and $\sfJ(\theta)$ are continuous
at $\theta=1$. At this endpoint, we have
\[
 \sfC_1(\xi)
 =
 I+(p-2)\frac{\xi}{|\xi|}\otimes\frac{\xi}{|\xi|},
 \qquad
 \sfD_1(\xi)=(2-p)\sfC_1(\xi),
\]
and
\[
 \sfj(1)=p-1,
 \qquad
 \sfJ(1)=-(p-1)^2.
\]
Since $\sfD_1(\xi)=(2-p)\sfC_1(\xi)$, we obtain
\begin{align*}
 \tr\bigl(\sfC_1(\xi)^{-1}\sfD_1(\xi)^2\bigr)
 &=
 (2-p)^2
 \tr\bigl(\sfC_1(\xi)^{-1}\sfC_1(\xi)^2\bigr)\\
 &=
 (p-2)^2\tr\sfC_1(\xi)\\
 &=
 (p-2)^2\bigl((N-1)+(p-1)\bigr)\\
 &=
 (N+p-2)(p-2)^2
 =A_*.
\end{align*}
Here we used that the trace of $\sfC_1(\xi)$
is the sum of its eigenvalues, namely $p-1$
in the direction of $\xi$ and $1$ in each
of the $N-1$ orthogonal directions.
Using $\sfJ(1)=-(p-1)^2$ and the preceding
trace identity, we obtain
\begin{align}\label{nn:eq:reg-margin}
 &\tr\bigl(\sfC_1(\xi)^{-1}\sfD_1(\xi)^2\bigr)
   +2(k-1)\sfJ(1)\notag\\
 &\qquad
 =A_*-2(k-1)(p-1)^2\notag\\
 &\qquad
 \le A_*
 -2\left(1+\frac{A_*}{(p-1)^2}\right)(p-1)^2
 \notag\\
 &\qquad
 =-A_*-2(p-1)^2<0.
\end{align}
Here we used the choice of $k$ in
\eqref{nn:eq:k}. The last inequality is strict
because $A_*\ge0$ and $p>1$. In fact, the stated choice of $k$ is not sharp, but this is irrelevant for the argument, since only the existence of some admissible $k=k(N,p)$ matters. The expression on the left-hand side of~\eqref{eq:trace+J} is a continuous
function of $\theta$ depending only on $N,p,k$;
in particular, it is independent of $\xi$.
By \eqref{nn:eq:reg-margin}, its value at $\theta=1$
is at most $-[A_*+2(p-1)^2]$.
Continuity therefore allows us to choose
$\theta_o=\theta_o(N,p,k)\in(0,1)$ such that
\begin{align*}
\tr\bigl(\sfC_\theta(\xi)^{-1}\sfD_\theta(\xi)^2\bigr)
   +2(k-1)\sfJ(\theta)
 &
 \le-\tfrac12\bigl[A_*+2(p-1)^2\bigr]
 \quad
 \forall\,\theta\in[\theta_o,1], \,\xi\ne0.
\end{align*}
The choice of $\theta_o$ is independent of
$\varepsilon$ and the lower bound $\sfm$.
We next obtain a uniform positive lower bound
for the eigenvalues of $\sfC_\theta(\xi)$ on
$[\theta_o,1]$. Recall that these eigenvalues are \( \theta^{-\frac{p-2}{2}}\) in directions orthogonal to \(\xi\), and \( \theta^{-\frac{p-2}{2}}\bigl(1+(p-2)\theta\bigr) \)
in the direction of \(\xi\).
For $\theta\in[\theta_o,1]$, we have
\[
 \theta^{-\frac{p-2}{2}}
 \ge\min\left\{1,\theta_o^{\frac{2-p}{2}}\right\},
 \qquad
 1+(p-2)\theta\ge\min\{1,p-1\}.
\]
For the eigenvalues in directions orthogonal to
$\xi$, it follows that
\[
 \theta^{-\frac{p-2}{2}}
 \ge
 \min\left\{1,\theta_o^{\frac{2-p}{2}}\right\}
 \ge
 \min\{1,p-1\}
 \min\left\{1,\theta_o^{\frac{2-p}{2}}\right\},
\]
where the last inequality uses
$0<\min\{1,p-1\}\le1$.
For the eigenvalue in the direction of $\xi$,
multiplying the two positive lower bounds gives
\[
 \theta^{-\frac{p-2}{2}}
 \bigl(1+(p-2)\theta\bigr)
 \ge
 \min\{1,p-1\}
 \min\left\{1,\theta_o^{\frac{2-p}{2}}\right\}.
\]
Thus, both eigenvalues are bounded below by
\[
 c_1:=
 \min\{1,p-1\}
 \min\left\{1,\theta_o^{\frac{2-p}{2}}\right\}>0.
\]
Since $\theta_o$ depends only on $N,p,k$, the same
is true of $c_1$. In particular, $c_1$ is independent
of $\varepsilon$, $\sfm$, and $\xi$.
Consequently, for every $\zeta\in\R^N$,
\[
 \zeta\cdot\sfC_\theta(\xi)\zeta
 \ge c_1|\zeta|^2
 \qquad
 \text{for }\theta\in[\theta_o,1],\ \xi\ne0.
\]
Suppose that
$\boldsymbol{\theta}_\varepsilon\ge\theta_o$.
The $j$-th column of $\nabla^2v_\varepsilon$ is
$\nabla(\partial_{x_j}v_\varepsilon)$. Applying
the preceding inequality to this vector gives
\[
 \nabla(\partial_{x_j}v_\varepsilon)\cdot
 \sfC_{\boldsymbol{\theta}_\varepsilon}
 (\nabla v_\varepsilon)
 \nabla(\partial_{x_j}v_\varepsilon)
 \ge
 c_1|\nabla(\partial_{x_j}v_\varepsilon)|^2.
\]
Using the symmetry of $\nabla^2v_\varepsilon$
and the cyclic property of the trace, we obtain
\begin{align*}
 \tr\Bigl(
   \sfC_{\boldsymbol{\theta}_\varepsilon}
   (\nabla v_\varepsilon)
   (\nabla^2v_\varepsilon)^2
 \Bigr)&=
 \tr\Bigl(
   (\nabla^2v_\varepsilon)^T
   \sfC_{\boldsymbol{\theta}_\varepsilon}
   (\nabla v_\varepsilon)
   \nabla^2v_\varepsilon
 \Bigr)\\
 &=
 \sum_{j=1}^N
 \nabla(\partial_{x_j}v_\varepsilon)\cdot
 \sfC_{\boldsymbol{\theta}_\varepsilon}
 (\nabla v_\varepsilon)
 \nabla(\partial_{x_j}v_\varepsilon)\\
 &\ge
 c_1\sum_{j=1}^N
 |\nabla(\partial_{x_j}v_\varepsilon)|^2\\
 &=
 c_1|\nabla^2v_\varepsilon|^2.
\end{align*}
We now substitute the lower bound for the Hessian
term and the negative upper bound for the coefficient
in square brackets into \eqref{nn:eq:reg-coercivity}.
For $\boldsymbol{\theta}_\varepsilon\ge\theta_o$,
this gives
\begin{align*}
 \calL_\varepsilon w_\varepsilon
 &\le
 -c_1a(v_\varepsilon)
 w_\varepsilon^{\frac{p-2}{2}}
 |\nabla^2v_\varepsilon|^2
 -\tfrac12\bigl[A_*+2(p-1)^2\bigr]
 \frac{a(v_\varepsilon)}{v_\varepsilon^2}
 w_\varepsilon^{\frac{p+2}{2}}.
\end{align*}
Here the lower bound for the trace gives an upper
bound for its negative multiple, since
$a(v_\varepsilon)w_\varepsilon^{\frac{p-2}{2}}>0$. Recalling $a(v)=cv^{2-p}$, we have \(a(v_\varepsilon)/v_\varepsilon^2 =cv_\varepsilon^{-p} \).
Consequently,
\begin{align*}
 \calL_\varepsilon w_\varepsilon
 &\le
 -c_1a(v_\varepsilon)
 w_\varepsilon^{\frac{p-2}{2}}
 |\nabla^2v_\varepsilon|^2
 -\tfrac12 c\bigl[A_*+2(p-1)^2\bigr]
 v_\varepsilon^{-p}
 w_\varepsilon^{\frac{p+2}{2}}.
\end{align*}
Setting
\[
 \kappa_o
 :=\tfrac12 c\bigl[A_*+2(p-1)^2\bigr]
 =\frac{k^{p-2}}{2(p-1)}
   \bigl[A_*+2(p-1)^2\bigr]>0,
\]
we obtain
\begin{equation}\label{nn:eq:reg-high}
 \calL_\varepsilon w_\varepsilon
 \le
 -c_1a(v_\varepsilon)
 w_\varepsilon^{\frac{p-2}{2}}|\nabla^2v_\varepsilon|^2
 -\kappa_o v_\varepsilon^{-p}
 w_\varepsilon^{\frac{p+2}{2}}
 \qquad
 \text{if }\boldsymbol{\theta}_\varepsilon\ge\theta_o.
\end{equation}
The constants $c_1$ and $\kappa_o$ depend only on
$N,p,k$ and are independent of $\varepsilon$
and the lower bound $\sfm$.
Recalling the definitions of $w_\varepsilon$ and
$\boldsymbol{\theta}_\varepsilon$, the condition
$\boldsymbol{\theta}_\varepsilon\ge\theta_o$
under which \eqref{nn:eq:reg-high} holds is
equivalent to
\[
 w_\varepsilon
 =|\nabla v_\varepsilon |^2
 \ge
 \frac{\theta_o}{1-\theta_o}
 \frac{\varepsilon^2}
      {k^2v_\varepsilon^{2(k-1)}}.
\]
\begin{remark}[The gradient threshold and the lower bound]\upshape
Since
\[
 |\nabla u_\varepsilon|^2
 =k^2v_\varepsilon^{2(k-1)}w_\varepsilon,
\]
the condition
$\boldsymbol{\theta}_\varepsilon\ge\theta_o$
can equivalently be expressed as
\[
 |\nabla u_\varepsilon|
 \ge
 \sqrt{\frac{\theta_o}{1-\theta_o}}\,\varepsilon.
\]
Thus, \eqref{nn:eq:reg-high} applies where the
magnitude of $\nabla u_\varepsilon$ is at least
a fixed multiple of $\varepsilon$. This multiple
depends only on $N,p,k$ and is independent of
the lower bound $\sfm$.

The lower bound enters when estimating the
transformed gradient in the complementary region.
Indeed, if $\boldsymbol{\theta}_\varepsilon<\theta_o$,
then
\[
 |\nabla v_\varepsilon|
 <
 \sqrt{\frac{\theta_o}{1-\theta_o}}\,
 \frac{\varepsilon}{kv_\varepsilon^{k-1}}
 \le
 \sqrt{\frac{\theta_o}{1-\theta_o}}\,
 \frac{\varepsilon}{k\sfm^{\frac{k-1}{k}}},
\]
where we used
$v_\varepsilon=u_\varepsilon^{\frac1k}\ge\sfm^{\frac1k}$.
For fixed $\sfm>0$, this bound tends to zero
as $\varepsilon\downarrow0$.
\end{remark}

\subsection{Bounds for the principal matrix and the drift}

We continue to assume
$\boldsymbol{\theta}_\varepsilon\ge\theta_o$.
For $\theta\in[\theta_o,1]$, we have
\[
 \theta^{-\frac{p-2}{2}}
 \le
 \max\left\{1,\theta_o^{-\frac{p-2}{2}}\right\},
 \qquad
 1+(p-2)\theta\le\max\{1,p-1\}.
\]
Thus, both eigenvalues of $\sfC_\theta(\xi)$
are bounded above by
\[
 \max\{1,p-1\}
 \max\left\{1,\theta_o^{-\frac{p-2}{2}}\right\},
\]
uniformly in $\theta\in[\theta_o,1]$ and $\xi\ne0$.
Since $\theta_o$ depends only on $N,p,k$,
this upper bound has the same dependence.
Since the squared Frobenius norm of a symmetric
matrix equals the sum of its squared eigenvalues,
the preceding upper bound gives
\[
 |\sfC_\theta(\xi)|
 \le
 \sqrt N\,\max\{1,p-1\}
 \max\left\{1,\theta_o^{-\frac{p-2}{2}}\right\}
 \qquad
 \text{for }\theta\in[\theta_o,1],\ \xi\ne0.
\]
Using the representation of $\sfA_\varepsilon$
and the homogeneity of the Frobenius norm,
we therefore obtain
\begin{align*}
 |\sfA_\varepsilon(v_\varepsilon,\nabla v_\varepsilon)|
 &=
 a(v_\varepsilon)w_\varepsilon^{\frac{p-2}{2}}
 \left|
   \sfC_{\boldsymbol{\theta}_\varepsilon}
   (\nabla v_\varepsilon)
 \right|\le
 C a(v_\varepsilon)w_\varepsilon^{\frac{p-2}{2}},
\end{align*}
provided that $\boldsymbol{\theta}_\varepsilon
\ge\theta_o$. Here we may take
\[
 C=
 \sqrt N\,\max\{1,p-1\}
 \max\left\{1,\theta_o^{-\frac{p-2}{2}}\right\},
\]
which depends only on $N,p,k$.

To estimate $\sfB_\varepsilon$, we compute
the partial derivatives of $\sfA_\varepsilon(v,\xi)$
and $\sfF_\varepsilon(v,\xi)$ with respect to the
gradient argument $\xi$, keeping $v>0$ fixed.
These are precisely the derivatives appearing
in the definition of $\sfB_\varepsilon$.
For $\xi\ne0$, differentiation gives
\[
 \nabla_\xi\boldsymbol{\Theta}_\varepsilon(v,\xi)
 =
 2\boldsymbol{\Theta}_\varepsilon(v,\xi)
 \bigl(1-\boldsymbol{\Theta}_\varepsilon(v,\xi)\bigr)
 \frac{\xi}{|\xi|^2}.
\]
Since $0<\boldsymbol{\Theta}_\varepsilon(v,\xi)<1$,
we have $\boldsymbol{\Theta}_\varepsilon(v,\xi)\bigl(1-\boldsymbol{\Theta}_\varepsilon(v,\xi)\bigr)\le\frac14$, and we obtain
\begin{equation*}
 \left|
   \partial_{\xi_\ell}
   \boldsymbol{\Theta}_\varepsilon(v,\xi)
 \right|
 \le
 \left|
   \nabla_\xi\boldsymbol{\Theta}_\varepsilon(v,\xi)
 \right|
 \le\frac1{2|\xi|}.
\end{equation*}

We next consider the explicit dependence
of $\sfC_\theta(\xi)$ on the direction of $\xi$.
For $\ell=1,\ldots,N$,
\[
 \partial_{\xi_\ell}\left(\frac{\xi}{|\xi|}\right)
 =
 \frac{e_\ell}{|\xi|}
 -\frac{\xi_\ell\xi}{|\xi|^3},
\]
where $e_\ell$ is the $\ell$-th standard basis vector.
In particular,
\[
 \left|
 \partial_{\xi_\ell}\left(\frac{\xi}{|\xi|}\right)
 \right|
 \le\frac1{|\xi|}.
\]
Keeping $\theta$ fixed and applying the product
rule to the tensor product in $\sfC_\theta(\xi)$
therefore yields
\begin{equation*}
 |\partial_{\xi_\ell}\sfC_\theta(\xi)|
 \le
 \frac{2|p-2|}{|\xi|}
 \theta^{-\frac{p-2}{2}}\theta
 \le\frac{C_1}{|\xi|}
 \qquad\text{for }\theta\in[\theta_o,1].
\end{equation*}
To account also for the dependence through
$\boldsymbol{\Theta}_\varepsilon(v,\xi)$, we use
\eqref{nn:eq:C-theta} and \eqref{nn:eq:j-theta}.
These formulas imply
\[
 |\partial_\theta\sfC_\theta(\xi)|\le C_2,
 \qquad
 |\sfj'(\theta)|\le C_3
 \qquad
 \text{for }\theta\in[\theta_o,1],\ \xi\ne0.
\]
The constants $C_1,C_2,C_3$ depend only on $N,p,k$.
By their definitions, $\sfC_\theta(\xi)$ and
$\sfj(\theta)$ are bounded on the same interval,
with the same parameter dependence.

Suppose now that
$\boldsymbol{\Theta}_\varepsilon(v,\xi)\ge\theta_o$.
The chain rule gives
\begin{align*}
 &\partial_{\xi_\ell}
 \Bigl(
   \sfC_{\boldsymbol{\Theta}_\varepsilon(v,\xi)}(\xi)
 \Bigr)\\
 &\quad=
 \left.
   \partial_{\xi_\ell}\sfC_\theta(\xi)
 \right|_{\theta=\boldsymbol{\Theta}_\varepsilon(v,\xi)}
 +
 \left.
   \partial_\theta\sfC_\theta(\xi)
 \right|_{\theta=\boldsymbol{\Theta}_\varepsilon(v,\xi)}
 \partial_{\xi_\ell}
   \boldsymbol{\Theta}_\varepsilon(v,\xi).
\end{align*}
Here the first term accounts for the explicit
dependence on $\xi$, with $\theta$ held fixed,
and the second for the dependence through
$\boldsymbol{\Theta}_\varepsilon(v,\xi)$.
Combining the preceding bounds, we obtain
\begin{align*}
 \left|
 \partial_{\xi_\ell}
 \Bigl(
   \sfC_{\boldsymbol{\Theta}_\varepsilon(v,\xi)}(\xi)
 \Bigr)
 \right|
 &\le
\frac{C_1+2C_2}{|\xi|}.
\end{align*}
The function $\sfj$ depends on $\xi$ only through
$\boldsymbol{\Theta}_\varepsilon(v,\xi)$.
Hence the chain rule gives
\[
 \nabla_\xi
 \Bigl(
   \sfj(\boldsymbol{\Theta}_\varepsilon(v,\xi))
 \Bigr)
 =
 \sfj'(\boldsymbol{\Theta}_\varepsilon(v,\xi))
 \nabla_\xi\boldsymbol{\Theta}_\varepsilon(v,\xi).
\]
Since $|\sfj'(\theta)|\le C_3$ on $[\theta_o,1]$,
with $C_3$ depending only on $N,p,k$, we obtain
\begin{align*}
 \left|
 \nabla_\xi
 \Bigl(
   \sfj(\boldsymbol{\Theta}_\varepsilon(v,\xi))
 \Bigr)
 \right|
 &=
 \left|
   \sfj'(\boldsymbol{\Theta}_\varepsilon(v,\xi))
 \right|
 \left|
   \nabla_\xi\boldsymbol{\Theta}_\varepsilon(v,\xi)
 \right|\le\frac{2 C_3}{|\xi|}.
\end{align*}
Taking $C=\max\{C_1+2C_2,2C_3\}$ proves both
estimates. This constant depends only on $N,p,k$
and is independent of $\varepsilon$ and $v$.
The first estimate includes both the dependence
through $\boldsymbol{\Theta}_\varepsilon(v,\xi)$
and the explicit dependence on the direction of $\xi$.

We now differentiate the representations of
$\sfA_\varepsilon$ and $\sfF_\varepsilon$.
Using
\[
 \nabla_\xi|\xi|^{p-2}
 =(p-2)|\xi|^{p-4}\xi,
 \qquad
 \nabla_\xi|\xi|^p=p|\xi|^{p-2}\xi,
\]
the product rule and the preceding bounds give
\[
 |\partial_{\xi_\ell}\sfA_\varepsilon(v,\xi)|
 \le C a(v)|\xi|^{p-3},
 \qquad
 |\nabla_\xi\sfF_\varepsilon(v,\xi)|
 \le C\frac{a(v)}v|\xi|^{p-1}\qquad\text{for }\boldsymbol\Theta_\varepsilon(v,\xi)\ge\theta_o.
\]
Here and below, $C$ depends only on $N,p,k$.
Evaluating at
$(v,\xi)=(v_\varepsilon,\nabla v_\varepsilon)$
and applying the Frobenius Cauchy--Schwarz inequality
in the definition of $\sfB_\varepsilon$, we obtain
\begin{align*}
 |\sfB_\varepsilon|
 &\le
 \left(
   \sum_{\ell=1}^N
   \big|
     (\partial_{\xi_\ell}\sfA_\varepsilon)
     (v_\varepsilon,\nabla v_\varepsilon)
   \big|^2
 \right)^{\frac12}
 |\nabla^2v_\varepsilon|
 +\big|
   (\nabla_\xi\sfF_\varepsilon)
   (v_\varepsilon,\nabla v_\varepsilon)
 \big|.
\end{align*}
Consequently,
\begin{equation}\label{nn:eq:reg-drift}
 |\sfB_\varepsilon|
 \le C\left(
   a(v_\varepsilon)
   w_\varepsilon^{\frac{p-3}{2}}|\nabla^2v_\varepsilon|
   +\frac{a(v_\varepsilon)}{v_\varepsilon}
   w_\varepsilon^{\frac{p-1}{2}}
 \right)\qquad\text{for }\boldsymbol\theta_\varepsilon\ge\theta_o.
\end{equation}
The constant is independent of $\varepsilon$
and the lower bound $\sfm$.

\subsection{Two cases at the maximum}

Let
\[
 Q:=B_{2r}(x_o)\times(s-2h,s],
 \qquad
 \overline Q\subset\Omega\times I,
\]
and suppose that
\[
 0<\sfm\le u_\varepsilon\le\sfM
 \qquad\text{in }Q,
\]
where $\sfm$ and $\sfM$ are independent of
$\varepsilon$.
Since $v_\varepsilon=u_\varepsilon^{\frac1k}$,
we have
\[
 \sfm^{\frac1k}
 \le v_\varepsilon
 \le \sfM^{\frac1k}.
\]
Choose a smooth cutoff function $\eta$ as in
Section~\ref{nn:sec:bernstein}. More precisely,
we require $0\le\eta\le1$, with $\eta=0$ on
the parabolic boundary of $Q$ and $\eta=1$ on
$B_r(x_o)\times(s-h,s]$, such that
\[
 |\nabla\eta|\le\frac Cr,
 \qquad
 |\nabla^2\eta|\le\frac C{r^2},
 \qquad
 |\partial_t\eta|\le\frac Ch,
\]
where $C$ depends only on $N$.
Consider
\[
 W_\varepsilon:=\eta^2w_\varepsilon
 \qquad\text{on }
 \overline{B_{2r}(x_o)}\times[s-2h,s].
\]
If its maximum is zero, then
$\nabla v_\varepsilon=0$ wherever $\eta=1$,
and the desired estimate follows immediately.
Otherwise, let $(x_\ast,t_\ast)$ be a point
where $W_\varepsilon$ attains a positive maximum.
In particular,
\[
 \eta(x_\ast,t_\ast)>0,
 \qquad
 w_\varepsilon(x_\ast,t_\ast)>0.
\]

Throughout the following calculation, all quantities are evaluated at $(x_\ast,t_\ast)$. We first derive an inequality valid for every value of $\boldsymbol\theta_\varepsilon(x_\ast,t_\ast)$, and only afterwards distinguish between two cases.

Since $w_\varepsilon(x_\ast,t_\ast)>0$,
the gradient $\nabla v_\varepsilon$ is nonzero
in a neighborhood of this point, so the
differentiated equation applies there. At the maximum point, we have
\[
 \nabla W_\varepsilon=0,
 \qquad
 \nabla^2W_\varepsilon\le0,
 \qquad
 \partial_tW_\varepsilon\ge0,
\]
where the time derivative is understood from
the left if $t_\ast=s$.
Since the principal matrix $\sfA_\varepsilon$
is positive definite, these relations and
the definition \eqref{nn:eq:reg-operator} imply
\[
 \calL_\varepsilon W_\varepsilon
 =
 \partial_tW_\varepsilon
 -\sfA_\varepsilon(v_\varepsilon,\nabla v_\varepsilon)
  \fprod\nabla^2W_\varepsilon
 -\sfB_\varepsilon\cdot\nabla W_\varepsilon
 \ge0
 \qquad\mbox{at \( (x_\ast,t_\ast)\).}
\]
Moreover,
\[
 0=\nabla W_\varepsilon
 =2\eta w_\varepsilon\nabla\eta
  +\eta^2\nabla w_\varepsilon.
\]
As $\eta>0$ at the maximum point \( (x_\ast, t_\ast)\), we obtain
\begin{equation}\label{nn:eq:reg-max-gradient}
 \nabla w_\varepsilon
 =-2w_\varepsilon\frac{\nabla\eta}{\eta}.
\end{equation}
In the following formulas, we write
$\sfA_\varepsilon$ for
$\sfA_\varepsilon(v_\varepsilon,\nabla v_\varepsilon)$.
The product rule for $\calL_\varepsilon$ gives
\begin{align*}
 \calL_\varepsilon W_\varepsilon
 &=
 \eta^2\calL_\varepsilon w_\varepsilon
 +2\eta w_\varepsilon\calL_\varepsilon\eta
 -2w_\varepsilon
   (\sfA_\varepsilon\nabla\eta)\cdot\nabla\eta
 -4\eta
   (\sfA_\varepsilon\nabla\eta)\cdot
   \nabla w_\varepsilon.
\end{align*}
Using \eqref{nn:eq:reg-max-gradient}, the last
term becomes
\[
 -4\eta
   (\sfA_\varepsilon\nabla\eta)\cdot
   \nabla w_\varepsilon
 =
 8w_\varepsilon
   (\sfA_\varepsilon\nabla\eta)\cdot\nabla\eta.
\]
Consequently,
\begin{align*}
 0&\le\calL_\varepsilon W_\varepsilon=
 \eta^2\calL_\varepsilon w_\varepsilon
 +2\eta w_\varepsilon
 \Bigl(
   \partial_t\eta
   -\sfA_\varepsilon\fprod\nabla^2\eta
   -\sfB_\varepsilon\cdot\nabla\eta
 \Bigr)
 +6w_\varepsilon
   (\sfA_\varepsilon\nabla\eta)\cdot\nabla\eta.
\end{align*}

Now we distinguish two cases according to whether
$\boldsymbol{\theta}_\varepsilon(x_\ast,t_\ast)
\ge\theta_o$ or
$\boldsymbol{\theta}_\varepsilon(x_\ast,t_\ast)
<\theta_o$.
Suppose first that
$\boldsymbol{\theta}_\varepsilon(x_\ast,t_\ast)
\ge\theta_o$. Then the coercive inequality
\eqref{nn:eq:reg-high}, the bound for the principal
matrix, and the drift estimate
\eqref{nn:eq:reg-drift} allow us to carry out
the localization argument from
Section~\ref{nn:sec:bernstein}. In particular, \eqref{nn:eq:reg-high} yields
\begin{align*}
 \eta^2\calL_\varepsilon w_\varepsilon
 &\le
 -c_1a(v_\varepsilon)\eta^2
   w_\varepsilon^{\frac{p-2}{2}}
   |\nabla^2v_\varepsilon|^2
 -\kappa_o\eta^2v_\varepsilon^{-p}
   w_\varepsilon^{\frac{p+2}{2}}.
\end{align*}
For the terms containing derivatives of $\eta$,
the bound
\[
 |\sfA_\varepsilon|
 \le Ca(v_\varepsilon)
      w_\varepsilon^{\frac{p-2}{2}}
\]
and the Cauchy--Schwarz inequality give
\[
 -2\eta w_\varepsilon
   \sfA_\varepsilon\fprod\nabla^2\eta
 \le
 Ca(v_\varepsilon)\eta
 w_\varepsilon^{\frac p2}|\nabla^2\eta|
\]
and
\[
 6w_\varepsilon
   (\sfA_\varepsilon\nabla\eta)\cdot\nabla\eta
 \le
 Ca(v_\varepsilon)w_\varepsilon^{\frac p2}
 |\nabla\eta|^2.
\]
Estimating the time derivative and the drift
term by their absolute values, and moving
the two negative terms to the left-hand side,
we obtain
\begin{align}\label{nn:eq:reg-localization-before-drift}
 c_1a(v_\varepsilon)&\eta^2
   w_\varepsilon^{\frac{p-2}{2}}
   |\nabla^2v_\varepsilon|^2
 +\kappa_o\eta^2v_\varepsilon^{-p}
   w_\varepsilon^{\frac{p+2}{2}}
 \notag\\
 &\le
 2\eta|\partial_t\eta|w_\varepsilon
 +Ca(v_\varepsilon)w_\varepsilon^{\frac p2}
   \bigl(
     |\nabla\eta|^2+\eta|\nabla^2\eta|
   \bigr)
 +2\eta w_\varepsilon
   |\sfB_\varepsilon|\,|\nabla\eta|.
\end{align}
Here and below, $C$ denotes a positive constant
depending only on $N,p,k$, whose value may
change from line to line.
We next estimate the drift term
$2\eta w_\varepsilon|\sfB_\varepsilon|\,|\nabla\eta|$
using \eqref{nn:eq:reg-drift}. This gives
\begin{align*}
 2\eta w_\varepsilon
 |\sfB_\varepsilon|\,|\nabla\eta|
 &\le
 C_{\mathrm d}a(v_\varepsilon)\eta
 w_\varepsilon^{\frac{p-1}{2}}
 |\nabla^2v_\varepsilon|\,|\nabla\eta|
 +C_{\mathrm d}
 \frac{a(v_\varepsilon)}{v_\varepsilon}
 \eta w_\varepsilon^{\frac{p+1}{2}}
 |\nabla\eta|.
\end{align*}
Here $C_{\mathrm d}$ depends only on $N,p,k$.
We denote this constant explicitly to track its
contribution when applying Young's inequality
to absorb the Hessian-dependent term into
the Hessian term on the left-hand side of
\eqref{nn:eq:reg-localization-before-drift}.
Indeed, Young's inequality gives
\begin{align*}
 C_{\mathrm d}a(v_\varepsilon)\eta
 w_\varepsilon^{\frac{p-1}{2}}
 |\nabla^2v_\varepsilon|\,|\nabla\eta|
 &=
 \left(
   \sqrt{a(v_\varepsilon)}\,\eta
   w_\varepsilon^{\frac{p-2}{4}}
   |\nabla^2v_\varepsilon|
 \right)
 \left(
   C_{\mathrm d}\sqrt{a(v_\varepsilon)}
   w_\varepsilon^{\frac p4}|\nabla\eta|
 \right)\\
 &\le
 \tfrac12 c_1 a(v_\varepsilon)\eta^2
 w_\varepsilon^{\frac{p-2}{2}}
 |\nabla^2v_\varepsilon|^2
 +\frac{C_{\mathrm d}^2}{2c_1}
 a(v_\varepsilon)w_\varepsilon^{\frac p2}
 |\nabla\eta|^2.
\end{align*}
Since $c_1$ and $C_{\mathrm d}$ depend only on
$N,p,k$, so does the constant in front of the
second term on the right-hand side. Substituting these estimates into
\eqref{nn:eq:reg-localization-before-drift},
we absorb the Hessian-dependent term from
the right-hand side into half of the Hessian
term on the left-hand side.
Dropping the remaining nonnegative Hessian
term on the left-hand side, we obtain
\begin{align}
 \kappa_o&\eta^2v_\varepsilon^{-p}
 w_\varepsilon^{\frac{p+2}{2}}\notag\\
 &\le
 C\Bigl[
   \eta|\partial_t\eta|w_\varepsilon
+a(v_\varepsilon)w_\varepsilon^{\frac p2}
   \bigl(
     |\nabla\eta|^2+\eta|\nabla^2\eta|
   \bigr)
   +\frac{a(v_\varepsilon)}{v_\varepsilon}
   \eta w_\varepsilon^{\frac{p+1}{2}}
   |\nabla\eta|
 \Bigr].
 \label{nn:eq:reg-cutoff-ineq}
\end{align}
We now set
\[
 \sfY_\varepsilon
 :=\eta\sqrt{w_\varepsilon}
 =\sqrt{W_\varepsilon}.
\]
Multiplying \eqref{nn:eq:reg-cutoff-ineq} by
$v_\varepsilon^p\eta^p$ yields
\begin{align*}
 \kappa_o&\sfY_\varepsilon^{p+2}\\
 &\le
 C\Bigl[
   v_\varepsilon^p\eta^{p-1}
   |\partial_t\eta|\sfY_\varepsilon^2
   +a(v_\varepsilon)v_\varepsilon^p
   \sfY_\varepsilon^p
   \bigl(
     |\nabla\eta|^2+\eta|\nabla^2\eta|
   \bigr)
   +a(v_\varepsilon)v_\varepsilon^{p-1}
   \sfY_\varepsilon^{p+1}|\nabla\eta|
 \Bigr].
\end{align*}
Since $a(v)=cv^{2-p}$, where
$c=k^{p-2}/(p-1)$, we have
\[
 a(v_\varepsilon)v_\varepsilon^p
 =cv_\varepsilon^2,
 \qquad
 a(v_\varepsilon)v_\varepsilon^{p-1}
 =cv_\varepsilon.
\]
Using $p>1$, $0\le\eta\le1$,
$v_\varepsilon\le\sfM^{\frac1k}$,
and the bounds for the derivatives of $\eta$,
we conclude that
\begin{equation}\label{nn:eq:reg-Y}
 \sfY_\varepsilon^{p+2}
 \le C\left(
   \frac{\sfM^{\frac pk}}h\sfY_\varepsilon^2
   +\frac{\sfM^{\frac2k}}{r^2}\sfY_\varepsilon^p
   +\frac{\sfM^{\frac1k}}r\sfY_\varepsilon^{p+1}
 \right).
\end{equation}
For every $\delta>0$, Young's inequality gives
\[
 \frac{\sfM^{\frac pk}}h\sfY_\varepsilon^2
 \le
 \delta\sfY_\varepsilon^{p+2}
 +C(p)\delta^{-\frac2p}
 \frac{\sfM^{\frac{p+2}{k}}}
      {h^{\frac{p+2}{p}}},
\]
as well as
\[
 \frac{\sfM^{\frac2k}}{r^2}\sfY_\varepsilon^p
 \le
 \delta\sfY_\varepsilon^{p+2}
 +C(p)\delta^{-\frac p2}
 \frac{\sfM^{\frac{p+2}{k}}}{r^{p+2}},
\]
and
\[
 \frac{\sfM^{\frac1k}}r\sfY_\varepsilon^{p+1}
 \le
 \delta\sfY_\varepsilon^{p+2}
 +C(p)\delta^{-(p+1)}
 \frac{\sfM^{\frac{p+2}{k}}}{r^{p+2}}.
\]
Choosing $\delta=1/(6C)$, with $C$ as in
\eqref{nn:eq:reg-Y}, allows us to absorb the
three terms containing $\sfY_\varepsilon^{p+2}$ on the right-hand side into the left-hand side.
Combining the remaining constants, we obtain
\[
 \sfY_\varepsilon^{p+2}
 \le
 C\sfM^{\frac{p+2}{k}}
 \left(
   \frac1{h^{\frac{p+2}{p}}}
   +\frac1{r^{p+2}}
 \right).
\]
Taking the $(p+2)$-th root gives
\[
 \sqrt{W_\varepsilon(x_\ast,t_\ast)}
 \le
 C\sfM^{\frac1k}
 \left(
   \frac1r+\frac1{h^{\frac1p}}
 \right),
\]
where $C$ depends only on $N,p,k$. This completes the first case.

On the other hand, if $\boldsymbol{\theta}_\varepsilon(x_\ast,t_\ast)
<\theta_o$, its definition gives
\[
 w_\varepsilon(x_\ast,t_\ast)
 <
 \frac{\theta_o}{1-\theta_o}
 \frac{\varepsilon^2}
 {k^2v_\varepsilon(x_\ast,t_\ast)^{2(k-1)}}.
\]
Since $\eta\le1$ and
$v_\varepsilon\ge\sfm^{\frac1k}$, it follows that
\begin{align*}
 \sqrt{W_\varepsilon(x_\ast,t_\ast)}
 &\le
 \sqrt{w_\varepsilon(x_\ast,t_\ast)}\le
 \frac1k\sqrt{\frac{\theta_o}{1-\theta_o}}\,
 \varepsilon\sfm^{-\frac{k-1}{k}}.
\end{align*}
Notice that the second case uses no differential inequality, only the definition of $\boldsymbol\theta_\varepsilon$ and $v_\varepsilon\ge\sfm^{\frac1k}$; this is the only place where $\sfm$ enters quantitatively.

Combining the two cases, we find
\[
 \sqrt{W_\varepsilon(x_\ast,t_\ast)}
 \le
 C\sfM^{\frac1k}
 \left(\frac1r+\frac1{h^{\frac1p}}\right)
 +C\varepsilon\sfm^{-\frac{k-1}{k}}.
\]
For every $(x,t)\in B_r(x_o)\times(s-h,s]$,
we have $\eta(x,t)=1$ and hence
\[
 |\nabla v_\varepsilon(x,t)|
 =\sqrt{W_\varepsilon(x,t)}
 \le\sqrt{W_\varepsilon(x_\ast,t_\ast)}\le
  C\sfM^{\frac1k}
 \left(\frac1r+\frac1{h^{\frac1p}}\right)
 +C\varepsilon\sfm^{-\frac{k-1}{k}}.
\]
Recalling
$v_\varepsilon=u_\varepsilon^{\frac1k}$,
we conclude that
\begin{equation}\label{nn:eq:reg-est}
 \left\|
   \nabla\bigl(u_\varepsilon^{\frac1k}\bigr)
 \right\|_{L^\infty(B_r(x_o)\times(s-h,s))}
 \le
 C\sfM^{\frac1k}
 \left(\frac1r+\frac1{h^{\frac1p}}\right)
 +C\varepsilon\sfm^{-\frac{k-1}{k}},
\end{equation}
where $C$ depends only on $N,p,k$ and is
independent of $\varepsilon,\sfm,\sfM,r,h$.

We summarize the estimate just established
in the following proposition.

\begin{proposition}[Gradient estimate for regularized solutions]
\label{nn:prop:reg-gradient}
Let $p>1$ and
\[
 k\ge
 2+\frac{(N+p-2)(p-2)^2}{(p-1)^2}.
\]
Let
\[
 Q:=B_{2r}(x_o)\times(s-2h,s],
 \qquad r,h>0,
 \qquad \overline Q\subset\Omega\times I.
\]
For $\varepsilon>0$, let $u_\varepsilon$ be a
smooth positive solution to~\eqref{nn:eq:reg}
in a neighborhood of $\overline Q$.
Suppose that
\[
 0<\sfm\le u_\varepsilon\le\sfM
 \qquad\text{in }Q.
\]
Then
\[
 \left\|
   \nabla\bigl(u_\varepsilon^{\frac1k}\bigr)
 \right\|_{L^\infty(B_r(x_o)\times(s-h,s))}
 \le
 C\sfM^{\frac1k}
 \left(
   \frac1r+\frac1{h^{\frac1p}}
 \right)
 +C\varepsilon\sfm^{-\frac{k-1}{k}},
\]
where $C$ depends only on $N,p,k$ and is
independent of $\varepsilon,\sfm,\sfM,r,h$.
In particular, if $\sfm$ and $\sfM$ are independent
of $\varepsilon$, the second term tends to zero
as $\varepsilon\downarrow0$.
\end{proposition}

\begin{remark}\upshape
For $\varepsilon=0$ one has $\boldsymbol\theta\equiv1$, the second case discussed above is void, the last term is absent, and Proposition~\ref{nn:prop:reg-gradient} becomes Lemma~\ref{nn:lem:gradient-estimate}. 
\end{remark}

\begin{remark}[The mechanism behind the gradient estimate]\upshape
The guiding idea is that a large gradient produces
a negative contribution in the differentiated
equation which limits its growth.
To reveal this mechanism, we pass from
$u_\varepsilon$ to
$v_\varepsilon=u_\varepsilon^{\frac1k}$ and study
$w_\varepsilon=|\nabla v_\varepsilon|^2$.
The power transformation is essential:
the choice of $k$ ensures that, when the gradient
is sufficiently large relative to the regularization
scale, the differentiated equation satisfies
\[
 \calL_\varepsilon w_\varepsilon
 \le
 -c_1a(v_\varepsilon)
 w_\varepsilon^{\frac{p-2}{2}}
 |\nabla^2v_\varepsilon|^2
 -\kappa_ov_\varepsilon^{-p}
 w_\varepsilon^{\frac{p+2}{2}}.
\]
The last term becomes increasingly negative
as the gradient grows and thus provides the
main control on its size.

To use this control locally without prescribing
the gradient on the boundary, we introduce
$W_\varepsilon=\eta^2w_\varepsilon$.
The cutoff function removes the influence of boundary
values of the gradient, but its derivatives
produce additional terms in the equation.
These terms represent the cost of localization:
their bounds increase as the spatial radius
or the time length of the cylinder decreases.

At a positive maximum of $W_\varepsilon$,
the maximum principle gives
$\calL_\varepsilon W_\varepsilon\ge0$.
On the other hand, the differentiated equation
provides the negative contributions displayed
above. Consequently, the terms generated by
the cutoff function must balance these contributions.
The retained negative Hessian term is used
to absorb the Hessian-dependent part of the
drift. After this absorption, setting
$\sfY_\varepsilon=\eta\sqrt{w_\varepsilon}$
and using $v_\varepsilon\le\sfM^{\frac1k}$
leads to
\[
 \sfY_\varepsilon^{p+2}
 \le C\left(
   \frac{\sfM^{\frac pk}}h\sfY_\varepsilon^2
   +\frac{\sfM^{\frac2k}}{r^2}\sfY_\varepsilon^p
   +\frac{\sfM^{\frac1k}}r\sfY_\varepsilon^{p+1}
 \right).
\]
The decisive point is that all powers on the
right-hand side are strictly smaller than
$p+2$. Therefore, the cutoff terms cannot balance
the left-hand side for arbitrarily large
$\sfY_\varepsilon$. This bounds the maximum
of $W_\varepsilon$ and hence the gradient
on the inner cylinder, where $\eta=1$.

The regularization requires a distinction between
two regimes at the maximum point.
If the gradient is sufficiently large relative
to the regularization scale, the coercive
inequality above applies.
In the complementary regime, the gradient is
already controlled directly by that scale:
\[
 |\nabla u_\varepsilon|
 <
 \sqrt{\frac{\theta_o}{1-\theta_o}}\,\varepsilon.
\]
The lower bound $u_\varepsilon\ge\sfm>0$
then prevents the power transformation from
amplifying this small gradient uncontrollably.
Indeed,
\[
 |\nabla v_\varepsilon|
 =
 \frac1k u_\varepsilon^{-\frac{k-1}{k}}
 |\nabla u_\varepsilon|
 \le C\varepsilon\sfm^{-\frac{k-1}{k}}.
\]
If $\sfm>0$ and $\sfM$ are independent of
$\varepsilon$, this additional term vanishes
as $\varepsilon\downarrow0$, while the remaining
bound is uniform in $\varepsilon$.

Thus, the coercive differential inequality is
needed only in the regime where the gradient
might be large. In the complementary regime,
the condition defining that regime already
provides the required control.
\end{remark}

\section{Approximation of positive solutions}

The goal of this section is the proof of

\begin{lemma}\label{lem:approximation}
Let $u$ be a continuous weak solution to
\eqref{PDE} in $\Omega_T$, and let
$Q_{R,S}\Subset\Omega_T$ be a cylinder such that
\[
 0<\sfm\le u\le\sfM
 \qquad\text{in }Q_{R,S}.
\]
Then there exist a sequence $\varepsilon_j\downarrow0$
and functions \( u_j\in C^\infty(Q_{R,S})
 \cap C^0(\overline{Q_{R,S}}) \)
solving
\begin{equation}\label{PDE-reg}
 \partial_t(u_j^{p-1})
 =
 \Div\bigl(\sfa_{\varepsilon_j}(\nabla u_j)\bigr)
 \qquad\text{in }Q_{R,S},
\end{equation}
where
\[
 \sfa_\varepsilon(\xi)
 :=
 (\varepsilon^2+|\xi|^2)^{\frac{p-2}{2}}\xi,
\]
such that
\[
 \tfrac12\sfm\le u_j\le2\sfM
 \qquad\text{in }\overline{Q_{R,S}},
\]
and, as $j\to\infty$,
\begin{equation*}
 u_j\to u
 \qquad\text{locally uniformly in }Q_{R,S},
\end{equation*}
and
\begin{equation*}
 u_j\to u
 \qquad\text{uniformly on }\partial_{\mathrm{par}}Q_{R,S}.
\end{equation*}
\end{lemma}

\begin{proof}
By translating the spatial and temporal coordinates,
we may assume that
\(  Q_{R,S}=B_R(0)\times(-S,0]\). We divide the proof into several steps.\\

\noindent\textbf{Step 1: Construction of solutions to the regularized problem.}
By mollifying $u$ in a neighborhood of
$\overline{Q_{R,S}}$, we obtain a sequence
$g_j\in C^\infty(\overline{Q_{R,S}})$ such that
\[
 \tfrac12\sfm\le g_j\le2\sfM
 \qquad\text{in }\overline{Q_{R,S}},
\]
and
\[
 g_j\to u
 \qquad\text{uniformly on }\overline{Q_{R,S}}
 \quad\text{as }j\to\infty.
\]
In particular,
\[
 g_j|_{\partial_{\mathrm{par}}Q_{R,S}}
 \to
 u|_{\partial_{\mathrm{par}}Q_{R,S}}
 \qquad\text{uniformly on }
 \partial_{\mathrm{par}}Q_{R,S}.
\]
Our first goal is to construct positive smooth
solutions to \eqref{PDE-reg}, with
$\varepsilon_j\in(0,1)$, satisfying
$u_j=g_j$ on $\partial_{\mathrm{par}}Q_{R,S}$.
For positive smooth functions, the identities
\[
 \partial_t(u_j^{p-1})
 =(p-1)u_j^{p-2}\partial_tu_j
\]
and
\[
 \Div\bigl(\sfa_{\varepsilon_j}(\nabla u_j)\bigr)
 =
 D_\xi\sfa_{\varepsilon_j}(\nabla u_j)
 \fprod \nabla^2u_j
\]
allow us to rewrite \eqref{PDE-reg} as
\begin{equation}\label{PDE2}
 \partial_tu_j
 -
 \frac{u_j^{2-p}}{p-1}
 D_\xi\sfa_{\varepsilon_j}(\nabla u_j)
 \fprod
 \nabla^2u_j
 =0
 \qquad\text{in }Q_{R,S}.
\end{equation}
Here $D_\xi\sfa_{\varepsilon_j}$ denotes
the Jacobian matrix of $\sfa_{\varepsilon_j}$
with respect to $\xi$.
To construct solutions to \eqref{PDE2}, we first
modify its coefficients to obtain a uniformly
parabolic auxiliary problem.
We begin with the scalar factor
$\frac{u^{2-p}}{p-1}$, extending it smoothly
to all real arguments while preserving its
values on the intended range of the solutions. More precisely, we choose a smooth function
$a_o\colon\R\to(0,\infty)$ such that
\[
 a_o(u)=\frac{u^{2-p}}{p-1}
 \qquad\text{for }
 u\in\left[\tfrac12\sfm,2\sfM\right],
\]
and
\[
 0<c_-\le a_o(u)\le c_+<\infty
 \qquad\text{for every }u\in\R,
 \qquad
 \|a_o'\|_{L^\infty(\R)}\le C.
\]
Here $u$ denotes an independent scalar variable,
and $c_-,c_+$ and $C$ depend only on
$p,\sfm,\sfM$.
The function $a_o$ should not be confused with
the coefficient $a(v)$ introduced earlier.
This modification ensures that the scalar
factor is bounded and bounded away from zero.
To obtain uniform parabolicity, we next
truncate the dependence on the gradient. For fixed $j$, choose
\[
 \sfK_\ast>
 \|\nabla g_j\|_{L^\infty(Q_{R,S})},
\]
to be specified sufficiently large later.
Let $\chi_{\sfK_\ast}\in C_0^\infty(\R^N)$ satisfy
\[
 0\le\chi_{\sfK_\ast}\le1,
 \qquad
 \chi_{\sfK_\ast}(\xi)=1
 \quad\text{if }|\xi|\le\sfK_\ast,
 \qquad
 \chi_{\sfK_\ast}(\xi)=0
 \quad\text{if }|\xi|\ge2\sfK_\ast.
\]
For independent variables $u\in\R$ and
$\xi\in\R^N$, we define
\[
 \sfE_{j,\sfK_\ast}(u,\xi)
 :=
 \chi_{\sfK_\ast}(\xi)a_o(u)
 D_\xi\sfa_{\varepsilon_j}(\xi)
 +(1-\chi_{\sfK_\ast}(\xi))I.
\]
Thus,
\[
 \sfE_{j,\sfK_\ast}(u,\xi)
 =
 \frac{u^{2-p}}{p-1}
 D_\xi\sfa_{\varepsilon_j}(\xi)
 \qquad\text{if }
 \tfrac12\sfm\le u\le2\sfM
 \ \text{and }|\xi|\le\sfK_\ast.
\]
For fixed $j$ and $\sfK_\ast$, the matrix
$\sfE_{j,\sfK_\ast}$ is uniformly positive definite
and bounded on $\R\times\R^N$.
The corresponding ellipticity bounds may depend
on $p,\sfm,\sfM,\varepsilon_j,\sfK_\ast$. Classical existence theory
(see, for example, \cite[Theorem~12.16]{Lieberman})
yields a classical solution
\( u_{j,\sfK_\ast}
 \in C^\infty(Q_{R,S})
 \cap C^1(\overline{Q_{R,S}})\)
of the Cauchy--Dirichlet problem
\begin{equation*}
\left\{
 \begin{array}{cl}
 \displaystyle
 \partial_tu_{j,\sfK_\ast}
 =
 \sfE_{j,\sfK_\ast}
 \bigl(u_{j,\sfK_\ast},\nabla u_{j,\sfK_\ast}\bigr)
 \fprod\nabla^2u_{j,\sfK_\ast}
 &\text{in }Q_{R,S},\\[1ex]
 u_{j,\sfK_\ast}=g_j
 &\text{on }\partial_{\mathrm{par}}Q_{R,S}.
 \end{array}\right.
\end{equation*}
Since $\frac12\sfm\le g_j\le2\sfM$ on
$\partial_{\mathrm{par}}Q_{R,S}$, the maximum
principle yields
\[
 \tfrac12\sfm\le u_{j,\sfK_\ast}\le2\sfM
 \qquad\text{in }\overline{Q_{R,S}}.
\]
Consequently, the definition of $a_o$ gives
\[
 a_o(u_{j,\sfK_\ast})
 =
 \frac{u_{j,\sfK_\ast}^{2-p}}{p-1}
 \qquad\text{in }Q_{R,S}.
\]

We claim that, for $\sfK_\ast$ sufficiently large,
\[
 |\nabla u_{j,\sfK_\ast}|<\sfK_\ast
 \qquad\text{throughout }Q_{R,S}.
\]
Then $\chi_{\sfK_\ast}(\nabla u_{j,\sfK_\ast})=1$,
so, by the preceding identity for
$a_o(u_{j,\sfK_\ast})$, the function
$u_{j,\sfK_\ast}$ solves the untruncated
regularized equation \eqref{PDE-reg}.
Suppose, towards a contradiction, that
$|\nabla u_{j,\sfK_\ast}|\ge\sfK_\ast$
at some point of $Q_{R,S}$.
Since $\nabla u_{j,\sfK_\ast}$ is continuous
on $\overline{Q_{R,S}}$, there is a first time
$t_\ast\in(-S,0)$ such that
\[
 \max_{x\in\overline{B_R}}
 |\nabla u_{j,\sfK_\ast}(x,t_\ast)|
 =\sfK_\ast.
\]
Here $t_\ast>-S$ because
$u_{j,\sfK_\ast}(\cdot,-S)=g_j(\cdot,-S)$
and
\[
 \max_{x\in\overline{B_R}}
 |\nabla g_j(x,-S)|<\sfK_\ast.
\]
By the definition of $t_\ast$,
\[
 |\nabla u_{j,\sfK_\ast}|<\sfK_\ast
 \qquad\text{in }
 Q_\ast:=B_R\times(-S,t_\ast).
\]
Thus, $u_{j,\sfK_\ast}$ solves the untruncated
regularized equation \eqref{PDE-reg} in $Q_\ast$. Set $v_j:=u_{j,\sfK_\ast}^{\frac1k}$, with
$k\ge k_o$ as in \eqref{nn:eq:k}.
We claim that
\begin{equation}\label{max-principle}
 \max_{\overline{Q_\ast}}|\nabla v_j|
 \le
 \max\left\{
   \max_{\partial_{\mathrm{par}}Q_\ast}|\nabla v_j|,
   C\varepsilon_j\sfm^{-\frac{k-1}{k}}
 \right\},
\end{equation}
where $C=C(N,p,k)$ is independent of
$j$ and $\sfK_\ast$.  Let
\[
 w_j:=|\nabla v_j|^2,
 \qquad
 \boldsymbol{\theta}_j
 :=
 \frac{w_j}{
   w_j+\frac{\varepsilon_j^2}
   {k^2v_j^{2(k-1)}}
 }.
\]
We apply the maximum argument from
Section~\ref{nn:sec:regularization} directly to
$w_j$, without introducing a cutoff function.
We show that the constant in
\eqref{max-principle} can be chosen as
\[
 C_\ast
 :=
 \frac{2^{\frac{k-1}{k}}}{k}
 \sqrt{\frac{\theta_o}{1-\theta_o}}.
\]
Suppose, towards a contradiction, that
\eqref{max-principle} fails with $C=C_\ast$.
Let $(\hat x,\hat t)$ be a point where $w_j$
attains its maximum on $\overline{Q_\ast}$.
Since this maximum strictly exceeds the maximum
on the parabolic boundary, we have
\(  \hat x\in B_R\) and \( \hat t\in(-S,t_\ast]\). At this point,
we have \(  \nabla w_j=0\), \( \nabla^2w_j\le0 \), and \( \partial_tw_j\ge0 \),
where the time derivative is understood from
the left if $\hat t=t_\ast$.
Consequently, the positive definiteness of
the principal matrix gives
\[
 \calL_{\varepsilon_j}w_j(\hat x,\hat t)\ge0.
\]
On the other hand, the assumed failure of
\eqref{max-principle} implies
\[
 w_j(\hat x,\hat t)
 >
 C_\ast^2\varepsilon_j^2
 \sfm^{-\frac{2(k-1)}{k}}.
\]
Using the definition of $C_\ast$ and the lower
bound $v_j\ge(\frac12 \sfm)^{\frac1k}$, we obtain
\begin{align*}
 w_j(\hat x,\hat t)
 &>
 \frac{\theta_o}{1-\theta_o}
 \frac{\varepsilon_j^2}{
   k^2(\frac12 \sfm)^{\frac{2(k-1)}{k}}
 }\ge
 \frac{\theta_o}{1-\theta_o}
 \frac{\varepsilon_j^2}{
   k^2v_j(\hat x,\hat t)^{2(k-1)}
 }.
\end{align*}
Equivalently,
$\boldsymbol{\theta}_j(\hat x,\hat t)>\theta_o$. The function $u_{j,\sfK_\ast}$ satisfies
the untruncated regularized equation
\eqref{PDE-reg} in $Q_\ast$.
If $\hat t=t_\ast$, the differentiated identity
for $w_j$ also holds at $(\hat x,t_\ast)$
by smoothness and passage to the limit
from earlier times.
Thus, \eqref{nn:eq:reg-high} applies and yields,
at $(\hat x,\hat t)$,
\begin{align*}
 \calL_{\varepsilon_j}w_j
 &\le
 -c_1a(v_j)w_j^{\frac{p-2}{2}}
   |\nabla^2v_j|^2
 -\kappa_o v_j^{-p}w_j^{\frac{p+2}{2}}<0,
\end{align*}
because $w_j(\hat x,\hat t)>0$.
This contradicts
$\calL_{\varepsilon_j}w_j(\hat x,\hat t)\ge0$
and proves \eqref{max-principle}.
The constant $C_\ast$ depends only on $N,p,k$.

To control the boundary term in \eqref{max-principle},
we first estimate the normal derivatives of
$u_{j,\sfK_\ast}$ on the lateral boundary of $Q_\ast$.
The corresponding bounds for
$v_j=u_{j,\sfK_\ast}^{\frac1k}$ then follow
by the chain rule.
Together with the prescribed initial values,
these boundary estimates allow us to apply
\eqref{max-principle} on the whole cylinder. For the boundary estimate, we introduce
\[
 w_\pm(x,t):=g_j(x,t)\pm\sfK f(x),
 \qquad
 f(x):=1-e^{|x|-R},
\]
where $\sfK>1$ will be chosen sufficiently large.
We compare these functions with $u_{j,\sfK_\ast}$
in the annular cylinder
\(\bigl(B_R\setminus\overline{B_{R/2}}\bigr)
 \times(-S,t_\ast) \),
where $f$ is smooth.
Since $f\ge0$ in $B_R$ and $f=0$ on $\partial B_R$,
we have \( w_-\le g_j\le w_+ \)
in $Q_\ast$
and
\( w_\pm=g_j=u_{j,\sfK_\ast}  \)
on \(\partial B_R\times(-S,t_\ast)\). On the inner lateral boundary, we have
$f=1-e^{-R/2}>0$.
We therefore impose the additional condition
\[
 \sfK\bigl(1-e^{-R/2}\bigr)
 \ge 2\sfM-\tfrac12\sfm.
\]
Using the bounds
$\frac12\sfm\le g_j,u_{j,\sfK_\ast}\le2\sfM$,
we obtain
\[
 w_-
 \le 2\sfM-\sfK\bigl(1-e^{-R/2}\bigr)
 \le\tfrac12\sfm
 \le u_{j,\sfK_\ast}
\]
and
\[
 w_+
 \ge\tfrac12\sfm+\sfK\bigl(1-e^{-R/2}\bigr)
 \ge2\sfM
 \ge u_{j,\sfK_\ast}
\]
on $\partial B_{R/2}\times(-S,t_\ast)$.
At the initial time, the same ordering follows
from $u_{j,\sfK_\ast}=g_j$ and $f\ge0$.
Thus,
\[
 w_-\le u_{j,\sfK_\ast}\le w_+
\]
on the entire parabolic boundary of the
annular cylinder.

For $0<|x|<R$, we have \( e^{-R}\le|\nabla f(x)|\le1\).
Thus, choosing
\[
 \sfK>
 (e^R+1)\|\nabla g_j\|_{L^\infty(Q_\ast)},
\]
the triangle inequality yields
\[
 \frac{\sfK}{e^R(e^R+1)}
 \le|\nabla w_\pm|
 \le2\sfK
\]
throughout the annular cylinder.

Differentiating the regularized flux gives
\[
 D_\xi\sfa_{\varepsilon_j}(\xi)
 =
 (\varepsilon_j^2+|\xi|^2)^{\frac{p-2}{2}}
 \left[
 I+(p-2)
 \frac{\xi\otimes\xi}{\varepsilon_j^2+|\xi|^2}
 \right].
\]
For $\xi\ne0$, the matrix in brackets has
eigenvalue $1$ in directions orthogonal to $\xi$
and eigenvalue
\[
 1+(p-2)\frac{|\xi|^2}{\varepsilon_j^2+|\xi|^2}
 \ge\min\{1,p-1\}
\]
in the direction of $\xi$.
Evaluating at $\xi=\nabla w_\pm$, we conclude
that the eigenvalues of
$D_\xi\sfa_{\varepsilon_j}(\nabla w_\pm)$
are bounded below by
\[
 \min\{1,p-1\}
 \bigl(\varepsilon_j^2+|\nabla w_\pm|^2\bigr)^
 {\frac{p-2}{2}}
 \ge c(p,R)\sfK^{p-2}.
\]
For $p\ge2$, the last inequality follows from
the lower bound for $|\nabla w_\pm|$.
For $1<p<2$, we instead use
\( \varepsilon_j^2+|\nabla w_\pm|^2
 \le 1+4\sfK^2
 \le5\sfK^2 \),
since $0<\varepsilon_j<1$ and $\sfK>1$,
and obtain
\[
 \bigl(\varepsilon_j^2+|\nabla w_\pm|^2\bigr)^
 {\frac{p-2}{2}}
 \ge
 5^{\frac{p-2}{2}}\sfK^{p-2}.
\]
Thus, the lower bound  for the eigenvalues of
$D_\xi\sfa_{\varepsilon_j}(\nabla w_\pm)$ holds for every $p>1$,
with $c(p,R)>0$ independent of
$j,\sfK$, and $\sfK_\ast$.

Moreover, for $0<|x|<R$, we have
\[
 \nabla^2f(x)
 =
 -e^{|x|-R}
 \left[
   \frac{x\otimes x}{|x|^2}
   +\frac1{|x|}
   \left(I-\frac{x\otimes x}{|x|^2}\right)
 \right].
\]
Its eigenvalues are $-e^{|x|-R}$ in the radial
direction and $-e^{|x|-R}/|x|$ in tangential
directions.
Since
$\nabla^2w_\pm=\nabla^2g_j\pm\sfK\nabla^2f$,
every eigenvalue of $\nabla^2w_+$ is bounded
above by
\[
 \|\nabla^2g_j\|_{L^\infty(Q_\ast)}
 -\sfK e^{-R}\min\Big\{1,\frac1R\Big\},
\]
whereas every eigenvalue of $\nabla^2w_-$
is bounded below by
\[
 -\|\nabla^2g_j\|_{L^\infty(Q_\ast)}
 +\sfK e^{-R}\min\Big\{1,\frac1R\Big\}.
\]
Choosing $\sfK$ so large that
\[
 \sfK e^{-R}\min\Big\{1,\frac1R\Big\}
 >
 \|\nabla^2g_j\|_{L^\infty(Q_{R,S})},
\]
we conclude that $\nabla^2w_+$ is negative
definite and $\nabla^2w_-$ is positive definite
throughout the annular cylinder.
Combining the lower eigenvalue bound for
$D_\xi\sfa_{\varepsilon_j}(\nabla w_\pm)$
with the preceding Hessian estimates, we obtain
\begin{align*}
 -D_\xi\sfa_{\varepsilon_j}(\nabla w_+)
 \fprod\nabla^2w_+
 &\ge
 c(p,R)\sfK^{p-2}
 \left[
   c(R)\sfK
   -\|\nabla^2g_j\|_{L^\infty(Q_\ast)}
 \right]
 >0
\end{align*}
and
\begin{align*}
 D_\xi\sfa_{\varepsilon_j}(\nabla w_-)
 \fprod\nabla^2w_-
 &\ge
 c(p,R)\sfK^{p-2}
 \left[
   c(R)\sfK
   -\|\nabla^2g_j\|_{L^\infty(Q_\ast)}
 \right]
 >0
\end{align*}
throughout the annular cylinder, where
\[
 c(R):=e^{-R}\min\Big\{1,\frac1R\Big\}.
\]
The strict positivity follows from our choice of $\sfK$.
Since $f$ is independent of time, we have \( \partial_tw_\pm=\partial_tg_j\), and therefore 
\(  |\partial_tw_\pm|
 \le\|\partial_tg_j\|_{L^\infty(Q_{R,S})}\). 
In particular, this bound is independent of $\sfK$. For the upper barrier, the preceding estimate gives
\begin{align*}
 &\partial_tw_+
 -\frac{w_+^{2-p}}{p-1}
 D_\xi\sfa_{\varepsilon_j}(\nabla w_+)
 \fprod\nabla^2w_+\\
 &\qquad\ge
 -\|\partial_tg_j\|_{L^\infty(Q_{R,S})}
 +\frac{w_+^{2-p}}{p-1}
 c(p,R)\sfK^{p-2}
 \left[
   c(R)\sfK
   -\|\nabla^2g_j\|_{L^\infty(Q_{R,S})}
 \right].
\end{align*}
By increasing $\sfK$ if necessary, we ensure that
the expression in brackets is at least
$\frac12c(R)\sfK$.
If $1<p\le2$, then $w_+\ge g_j\ge\frac12\sfm$
implies
\[
 \frac{w_+^{2-p}}{p-1}
 \ge\frac{(\frac12 \sfm)^{2-p}}{p-1}>0.
\]
Since $\sfK^{p-1}\to\infty$ as $\sfK\to\infty$,
the positive term dominates the time derivative
for sufficiently large $\sfK$.
Thus, $w_+$ is a supersolution to \eqref{PDE2}
throughout the annular cylinder.
If $p>2$, we instead use
\[
 \frac{w_+^{2-p}}{p-1}
 \ge\frac{(2\sfM)^{2-p}}{p-1}>0
 \qquad\text{wherever }w_+\le2\sfM.
\]
The same argument shows that, for sufficiently
large $\sfK$, the supersolution inequality holds
at every point of the annular cylinder where
$w_+\le2\sfM$. Consequently, by \cite[Lemma A.1]{BDL-holder} the truncated function
$\min\{w_+,2\sfM\}$ is a weak supersolution
to \eqref{PDE-reg} in this case. 
Similarly, for sufficiently large $\sfK$, the
function $w_-$ satisfies the subsolution
inequality wherever $w_->0$ if $p\ge2$,
and wherever $w_->\frac12\sfm$ if $1<p<2$.
Therefore, in the same way, the truncated function
$\max\{w_-,\frac12\sfm\}$ is a weak subsolution
to \eqref{PDE-reg} in the annular cylinder
for every $p>1$.
The functions
$\max\{w_-,\frac12\sfm\}$,
$u_{j,\sfK_\ast}$, and $\min\{w_+,2\sfM\}$
are bounded and bounded away from zero.
Together with the ordering on the parabolic
boundary of the annular cylinder, the comparison
principle from \cite{BS} therefore yields
\[
 \max\{w_-,\tfrac12\sfm\}
 \le u_{j,\sfK_\ast}
 \le \min\{w_+,2\sfM\}
 \qquad\text{in }
 \bigl(B_R\setminus\overline{B_{R/2}}\bigr)
 \times(-S,t_\ast).
\]
For $|x|\le \frac12 R$, we have
$f(x)\ge1-e^{-R/2}$.
Recalling that $\sfK$ was chosen so that
\[
 \sfK(1-e^{-R/2})
 \ge 2\sfM-\tfrac12\sfm,
\]
we also obtain
\[
 g_j-\sfK f
 \le\tfrac12\sfm
 \le u_{j,\sfK_\ast}
 \le2\sfM
 \le g_j+\sfK f
\]
in this inner region.
Consequently,
\[
 g_j(x,t)-\sfK f(x)
 \le u_{j,\sfK_\ast}(x,t)
 \le g_j(x,t)+\sfK f(x)
 \qquad\text{for every }(x,t)\in Q_\ast.
\]
At every point
$(x,t)\in\partial B_R\times(-S,t_\ast)$,
the boundary condition and $f=0$ give
\[
 w_-(x,t)
 =u_{j,\sfK_\ast}(x,t)
 =w_+(x,t)
 =g_j(x,t).
\]
Fix an arbitrary point
$(x,t)\in\partial B_R\times(-S,t_\ast)$.
Since both barriers agree with
$u_{j,\sfK_\ast}$ at this boundary point,
taking radial difference quotients from
inside $B_R$ in the preceding comparison
inequality gives
\[
 \partial_r w_+(x,t)
 \le
 \partial_r u_{j,\sfK_\ast}(x,t)
 \le
 \partial_r w_-(x,t).
\]
Here $\partial_r$ denotes differentiation
in the outward radial direction.
Using $\partial_r f=-1$ on $\partial B_R$,
we obtain
\[
 \partial_r g_j(x,t)-\sfK
 \le
 \partial_r u_{j,\sfK_\ast}(x,t)
 \le
 \partial_r g_j(x,t)+\sfK.
\]
Consequently,
\[
 \left|
   \partial_r u_{j,\sfK_\ast}(x,t)
 \right|
 \le
 \left|\partial_r g_j(x,t)\right|+\sfK
 \le
 \|\nabla g_j\|_{L^\infty(Q_{R,S})}+\sfK.
\]
Since the boundary point was arbitrary,
this estimate holds on
$\partial B_R\times(-S,t_\ast)$
and, by continuity, also at $t=t_\ast$. Since $u_{j,\sfK_\ast}=g_j$ on the lateral
boundary, their tangential derivatives agree
there. Hence the difference
$\nabla u_{j,\sfK_\ast}-\nabla g_j$
is normal to $\partial B_R$.
The preceding estimate therefore gives
\[
 |\nabla u_{j,\sfK_\ast}-\nabla g_j|
 =
 |\partial_r u_{j,\sfK_\ast}-\partial_r g_j|
 \le\sfK
 \qquad\text{on }
 \partial B_R\times(-S,t_\ast].
\]
On the initial boundary, we have
\[
 \nabla u_{j,\sfK_\ast}(\cdot,-S)
 =\nabla g_j(\cdot,-S).
\]
Consequently,
\[
 \max_{\partial_{\mathrm{par}}Q_\ast}
 |\nabla u_{j,\sfK_\ast}|
 \le
 \|\nabla g_j\|_{L^\infty(Q_{R,S})}+\sfK.
\]
Since $v_j=u_{j,\sfK_\ast}^{\frac1k}$
and $u_{j,\sfK_\ast}\ge\frac12\sfm$,
the chain rule also yields
\[
 \max_{\partial_{\mathrm{par}}Q_\ast}
 |\nabla v_j|
 \le
 \frac1k
 \left(\frac2{\sfm}\right)^{\frac{k-1}{k}}
 \left(
   \|\nabla g_j\|_{L^\infty(Q_{R,S})}+\sfK
 \right).
\]
This bound is independent of $\sfK_\ast$,
since the barrier parameter $\sfK$ was chosen
independently of the truncation threshold. The preceding estimate bounds the boundary term
independently of $\sfK_\ast$.
For fixed $j$, we may therefore choose
the truncation threshold $\sfK_\ast$
sufficiently large, possibly depending on $j$,
such that
\[
 \max\left\{
   \max_{\partial_{\mathrm{par}}Q_\ast}
     |\nabla v_j|,
   C_\ast\varepsilon_j
     \sfm^{-\frac{k-1}{k}}
 \right\}
 <
 \frac{\sfK_\ast}
      {2k(2\sfM)^{\frac{k-1}{k}}}.
\]
With this choice of $\sfK_\ast$, the chain rule,
the bound $u_{j,\sfK_\ast}\le2\sfM$, and
\eqref{max-principle} give
\begin{align*}
 \max_{\overline Q_\ast}
 |\nabla u_{j,\sfK_\ast}|
 &\le
 k(2\sfM)^{\frac{k-1}{k}}
 \max_{\overline Q_\ast}|\nabla v_j|<
 \tfrac12\sfK_\ast.
\end{align*}
This contradicts the choice of $t_\ast$, since
\[
 \max_{\overline{B_R}}
 |\nabla u_{j,\sfK_\ast}(\cdot,t_\ast)|
 =\sfK_\ast.
\]
Consequently,
\[
 |\nabla u_{j,\sfK_\ast}|<\sfK_\ast
 \qquad\text{throughout }Q_{R,S}.
\]
Thus, the gradient cutoff is identically one
along the solution, and $u_{j,\sfK_\ast}$
solves the untruncated regularized equation
\eqref{PDE-reg} on the whole cylinder $Q_{R,S}$. For each $j$, we fix such a truncation threshold
and set
\( u_j:=u_{j,\sfK_\ast}\).
This gives the desired smooth solutions
to \eqref{PDE-reg}. \\

\noindent\textbf{Step 2: Locally uniform convergence.}
Let $Q\Subset Q_{R,S}$ be a cylinder.
Choose another cylinder $\widehat Q$ such that \( Q\Subset\widehat Q\Subset Q_{R,S}\).
Recall that \( \frac12\sfm\le u_j\le2\sfM\) in  \( Q_{R,S}\).
We cover $\overline{\widehat Q}$ by finitely
many cylinders of the form
$B_r(x)\times(s-h,s]$, whose corresponding
larger cylinders
$B_{2r}(x)\times(s-2h,s]$
are compactly contained in $Q_{R,S}$.
Applying Proposition~\ref{nn:prop:reg-gradient}
on these cylinders gives
\[
 \big\|
   \nabla\bigl(u_j^{\frac1k}\bigr)
 \big\|_{L^\infty(\widehat Q)}
 \le
 C\left(
   \sfM^{\frac1k}
   +\frac{\varepsilon_j}
          {\sfm^{\frac{k-1}{k}}}
 \right).
\]
Here $C$ depends only on $N,p,k$ and the
fixed cylinders $\widehat Q$ and $Q_{R,S}$;
the factors arising from the bounds
$\frac12\sfm\le u_j\le2\sfM$ have been
absorbed into $C$. By the chain rule, we have \( \nabla u_j
 =
 k u_j^{\frac{k-1}{k}}
 \nabla\bigl(u_j^{\frac1k}\bigr) \).
Using $u_j\le2\sfM$ and
$0<\varepsilon_j<1$, we therefore obtain
\begin{align*}
 \sfL
 &:=
 \sup_{j\in\N}
 \|\nabla u_j\|_{L^\infty(\widehat Q)}
 \le
 Ck(2\sfM)^{\frac{k-1}{k}}
 \left(
   \sfM^{\frac1k}
   +\frac1{\sfm^{\frac{k-1}{k}}}
 \right)
 <\infty,
\end{align*}
where we used that $\sfm>0$ and $\sfM$
are fixed.

To prove that the sequence $(u_j)$ is
equicontinuous on $Q$, it remains to establish
a modulus of continuity in time that is
independent of $j$. Let $\eta\in C^\infty_0(B_1)$ be a standard
mollifying kernel satisfying $\eta\ge0$ and
$\int_{B_1}\eta(y)\,\dy=1$, and set
$\eta_\delta(y):=\delta^{-N}\eta(y/\delta)$
for $\delta>0$. For
\(0<\delta<\dist(Q,\partial\widehat Q) \),
we define the spatial mollification by
\[
 (u_j^{p-1})_\delta(x,t)
 :=
 \int_{B_\delta(x)}
 u_j(y,t)^{p-1}\eta_\delta(x-y)\,\dy,
 \qquad (x,t)\in Q.
\]
The restriction on $\delta$ ensures that
the integral involves only values of $u_j$
in $\widehat Q$. The spatial Lipschitz bounds for $u_j$ on
$\widehat Q$ and the normalization of
$\eta_\delta$ give
\begin{align}\label{molli-converge}
 \big|
   (u_j^{p-1})_\delta(x,t)-u_j^{p-1}(x,t)
 \big|
 &\le
 \delta
 \|\nabla(u_j^{p-1})\|_{L^\infty(\widehat Q)}
 \le
 C(p,\sfm,\sfM)\sfL\delta
\end{align}
for every $(x,t)\in Q$.
Similarly,
\begin{align}\label{Lipschitz-in-x}
 \big|
   u_j^{p-1}(x_1,t)&-u_j^{p-1}(x_2,t)
 \big|\notag\\
 &\le
 \|\nabla(u_j^{p-1})\|_{L^\infty(Q)}
 |x_1-x_2|
\le
 C(p,\sfm,\sfM)\sfL|x_1-x_2|
\end{align}
whenever $(x_1,t),(x_2,t)\in Q$.
In both estimates, we use the chain rule
\[
 \nabla(u_j^{p-1})
 =(p-1)u_j^{p-2}\nabla u_j.
\]
Together with the bounds
$\frac12\sfm\le u_j\le2\sfM$, this allows us
to choose
\[
 C(p,\sfm,\sfM)
 :=
 (p-1)\max\left\{
   \left(\tfrac12 \sfm\right)^{p-2},
   (2\sfM)^{p-2}
 \right\}.
\]

Fix $x$ and times $t_1<t_2$ such that
$(x,t_1),(x,t_2)\in Q$.
In the weak formulation of \eqref{PDE-reg},
we use test functions of the form
\( \varphi(y,\tau)
 :=
 \eta_\delta(x-y)\zeta(\tau) \),
where $\zeta$ is a smooth temporal cutoff
approximating the characteristic function
of $[t_1,t_2]$.
Passing to the limit in this approximation
gives
\begin{align}\label{glueing}
 &\big|
   (u_j^{p-1})_\delta(x,t_2)
   -(u_j^{p-1})_\delta(x,t_1)
 \big|
 \notag\\
 &\quad=
 \left|
   \int_{t_1}^{t_2}\int_{B_\delta(x)}
   \sfa_{\varepsilon_j}(\nabla u_j(y,\tau))
   \cdot\nabla_y\eta_\delta(x-y)
   \,\dy\,\dtau
 \right|
 \notag\\
 &\quad\le
 C(\sfL+\varepsilon_j)^{p-1}
 \frac{t_2-t_1}{\delta}.
\end{align}
Here we used
\[
 \big|
   \sfa_{\varepsilon_j}(\nabla u_j)
 \big|
 \le
 \bigl(\varepsilon_j+|\nabla u_j|\bigr)^{p-1}
 \le
 (\varepsilon_j+\sfL)^{p-1}
 \qquad\text{in }\widehat Q,
\]
together with
\[
 \int_{B_\delta(x)}
 |\nabla_y\eta_\delta(x-y)|\,\dy
 =
 \frac1\delta
 \int_{B_1}|\nabla\eta(y)|\,\dy.
\]
Thus, we may take
$C=\|\nabla\eta\|_{L^1(B_1)}$,
which is independent of $j,\delta,t_1,t_2$.

Let $(x_1,t_1),(x_2,t_2)\in Q$.
Since $Q$ is a cylinder, we also have
$(x_1,t_2)\in Q$.
Combining \eqref{molli-converge},
\eqref{Lipschitz-in-x}, and \eqref{glueing},
and using $0<\varepsilon_j<1$, we obtain
\begin{align*}
 &\big|
   u_j^{p-1}(x_1,t_1)-u_j^{p-1}(x_2,t_2)
 \big|\\
 &\qquad\le
 \big|
   u_j^{p-1}(x_1,t_1)
   -(u_j^{p-1})_\delta(x_1,t_1)
 \big|
 +\big|
   (u_j^{p-1})_\delta(x_1,t_1)
   -(u_j^{p-1})_\delta(x_1,t_2)
 \big|\\
 &\qquad\phantom{\le\,}
 +\big|
   (u_j^{p-1})_\delta(x_1,t_2)
   -u_j^{p-1}(x_1,t_2)
 \big|
 +\big|
   u_j^{p-1}(x_1,t_2)-u_j^{p-1}(x_2,t_2)
 \big|\\
 &\qquad\le
 C\left[
   \delta
   +\frac{|t_2-t_1|}{\delta}
   +|x_1-x_2|
 \right]
\end{align*}
for every
$0<\delta<\dist(Q,\partial\widehat Q)$.
Here $C$ depends only on $p,\sfm,\sfM,\sfL$
and the fixed mollifying kernel $\eta$,
and is independent of $j$, $\delta$,
and the chosen points. Choosing first $\delta$ and then
$|t_2-t_1|$ and $|x_1-x_2|$ sufficiently small,
we conclude that $(u_j^{p-1})_{j\in\N}$ is
equicontinuous on $Q$.
Since 
$\sigma\mapsto\sigma^{\frac1{p-1}}$
is Lipschitz on the fixed interval
\( \bigl[
   \left(\tfrac12 \sfm\right)^{p-1},
   (2\sfM)^{p-1}
 \bigr]\),
the sequence $(u_j)_{j\in\N}$ is equicontinuous
on $Q$ as well. The sequence is also uniformly bounded.
Applying the Arzel\'a--Ascoli theorem on
compact cylinders and using a diagonal
argument along an exhaustion of $Q_{R,S}$,
we therefore obtain a subsequence, not
relabeled, such that
\[
 u_j\to w
 \qquad\text{locally uniformly in }Q_{R,S}
 \;\text{as }j\to\infty,
\]
where
$w\colon Q_{R,S}\to[\frac12\sfm,2\sfM]$
is continuous.\\

\noindent\textbf{Step 3: Initial and boundary
values of the limit map.}
To show that the limit $w$ attains the
prescribed boundary values, we fix a lateral
boundary point
\[
 (x_{\mathrm b},t_{\mathrm b})
 \in\partial B_R\times[-S,0].
\] 
We choose an exterior ball
$B_\varrho(y)\subset\Omega\setminus\overline{B_R}$
that is externally tangent to $B_R$
at $x_{\mathrm b}$.
For a parameter $\gamma>0$ satisfying
\[
 (\gamma+1)\min\{1,p-1\}>N-1,
\]
we define
\[
 f(x)
 :=
 1-\left(\frac{\varrho}{|x-y|}\right)^\gamma,
 \qquad x\in B_R.
\]
We note that $f(x_{\mathrm b})=0$
and $f>0$ in $B_R$.
Let $\sfK>0$ be a constant
to be chosen later. For each $x\in B_R$, we choose an orthonormal
basis of $\R^N$ consisting of the radial
unit vector $\frac{x-y}{|x-y|}$ and $N-1$
vectors orthogonal to $x-y$.
These remaining vectors are tangent to the
sphere $\partial B_{|x-y|}(y)$ at $x$. The vectors of this basis are eigenvectors
of both matrices, as follows from the
explicit formulas
\[
 \nabla^2f(x)
 =
 \frac{\gamma\varrho^\gamma}{|x-y|^{\gamma+2}}
 \Biggl[
   I-(\gamma+2)
   \frac{x-y}{|x-y|}
   \otimes
   \frac{x-y}{|x-y|}
 \Biggr]
\]
and
\begin{align*}
 &D_\xi\sfa_{\varepsilon_j}
   \bigl(\sfK\nabla f(x)\bigr)\\
 &\quad=
 \bigl(
   \varepsilon_j^2+\sfK^2|\nabla f(x)|^2
 \bigr)^{\frac{p-2}{2}}
 \Biggl[
   I+(p-2)
   \frac{\sfK^2|\nabla f(x)|^2}
        {\varepsilon_j^2+\sfK^2|\nabla f(x)|^2}
   \frac{x-y}{|x-y|}
   \otimes
   \frac{x-y}{|x-y|}
 \Biggr].
\end{align*}
Indeed, the tensor product appearing in
both formulas acts as the identity in the
radial direction and vanishes in every
tangential direction.
Computing the Frobenius inner product in
this common orthonormal eigenbasis, we obtain
\begin{align}\label{supersol}
 &-D_\xi\sfa_{\varepsilon_j}
   \bigl(\sfK\nabla f(x)\bigr)
   \fprod\bigl(\sfK\nabla^2f(x)\bigr)
 \notag\\
 &\qquad=
 \frac{
   |\sfa_{\varepsilon_j}(\sfK\nabla f(x))|
 }{|x-y|}
 \Biggl[
   (\gamma+1)
   \Biggl(
   \underbrace{
     1+(p-2)
     \frac{\sfK^2|\nabla f(x)|^2}
          {\varepsilon_j^2+\sfK^2|\nabla f(x)|^2}}_{\ge \min\{1,p-1\}}
   \Biggr)
   -(N-1)
 \Biggr]
 \notag\\
 &\qquad\ge
 \frac{
   |\sfa_{\varepsilon_j}(\sfK\nabla f(x))|
 }{|x-y|}
 \Bigl[
   (\gamma+1)\min\{1,p-1\}-(N-1)
 \Bigr]
 >0.
\end{align}
The strict positivity follows from the choice
of $\gamma$, together with $\sfK>0$ and
$|\nabla f(x)|>0$. To make the dependence on $\sfK$ explicit,
we use the definition of the flux to write
\begin{align*}
 \left|
   \sfa_{\varepsilon_j}
   \bigl(\sfK\nabla f(x)\bigr)
 \right|
 &=
 \sfK^{p-1}|\nabla f(x)|^{p-1}
 \left(
   1+\frac{\varepsilon_j^2}
           {\sfK^2|\nabla f(x)|^2}
 \right)^{\frac{p-2}{2}}.
\end{align*}
Since $|\nabla f|$ has a strictly positive
minimum on $\overline{B_R}$ and
$0<\varepsilon_j<1$, the last factor is
bounded below by a positive constant,
uniformly in $x$, $j$, and $\sfK\ge1$.
Together with the positive minimum of
$\frac1{|x-y|}$ and the choice of $\gamma$,
this shows that the right-hand side of
\eqref{supersol} is bounded below by a
positive constant times $\sfK^{p-1}$.
Since $p>1$, it can therefore be made
arbitrarily large by choosing $\sfK$
sufficiently large, uniformly in
$x\in\overline{B_R}$ and $j\in\N$.

For parameters $\delta\in(0,1)$ and
$\Lambda>0$, we define the barrier functions
\[
 w_\pm(x,t)
 :=
 u(x_{\mathrm b},t_{\mathrm b})
 \pm\Bigl[
   2\delta+\sfK f(x)
   +\Lambda(t-t_{\mathrm b})^2
 \Bigr],
 \qquad (x,t)\in Q_{R,S}.
\]
For fixed $\delta\in(0,1)$, the continuity
of $u$ and the uniform convergence
$g_j\to u$ on $\partial_{\mathrm{par}}Q_{R,S}$
yield a neighborhood $U$ of
$(x_{\mathrm b},t_{\mathrm b})$
and an index $j_\delta\in\N$ such that
\[
 |g_j(x,t)-u(x_{\mathrm b},t_{\mathrm b})|
 \le
 |g_j(x,t)-u(x,t)|
 +|u(x,t)-u(x_{\mathrm b},t_{\mathrm b})|
 <2\delta
\]
for every
$(x,t)\in U\cap\partial_{\mathrm{par}}Q_{R,S}$
and every $j\ge j_\delta$.
Since $f\ge0$ on $\overline{B_R}$ and
$\Lambda(t-t_{\mathrm b})^2\ge0$, it follows that
\[
 w_-\le g_j\le w_+
 \qquad\text{on }
 U\cap\partial_{\mathrm{par}}Q_{R,S},
 \qquad j\ge j_\delta.
\]
The neighborhood $U$ and the index $j_\delta$
can be chosen independently of $\sfK$
and $\Lambda$. Outside $U$, the function
$f(x)+(t-t_{\mathrm b})^2$ has a strictly
positive minimum on the closed parabolic
boundary. We can therefore choose $\Lambda$ and $\sfK$
sufficiently large in terms of $\sfM$ and
the positive minimum above, independently
of $j$, such that
\[
 w_-\le g_j=u_j\le w_+
 \qquad\text{on }
 \partial_{\mathrm{par}}Q_{R,S},
 \qquad j\ge j_\delta.
\]
Our next goal is to verify that $w_+$ is
a supersolution to \eqref{PDE2}.
Since $t,t_{\mathrm b}\in[-S,0]$, we have
\[
 \partial_t w_+
 =
 2\Lambda(t-t_{\mathrm b})
 \ge-2\Lambda S.
\]
Moreover, \(  \nabla w_+=\sfK\nabla f\) and
\( \nabla^2w_+=\sfK\nabla^2f\).
Keeping $\Lambda$ fixed, it therefore
suffices to choose $\sfK$ sufficiently
large that
\[
 -\frac{w_+^{2-p}}{p-1}
 D_\xi\sfa_{\varepsilon_j}
 \bigl(\sfK\nabla f(x)\bigr)
 \fprod\bigl(\sfK\nabla^2f(x)\bigr)
 \ge2\Lambda S
 \qquad\text{in }Q_{R,S}.
\]
This ensures that
\[
 \partial_t w_+
 -\frac{w_+^{2-p}}{p-1}
 D_\xi\sfa_{\varepsilon_j}(\nabla w_+)
 \fprod\nabla^2w_+
 \ge0.
\]
If $1<p\le2$, then
$w_+\ge u(x_{\mathrm b},t_{\mathrm b})\ge\sfm$,
and hence
\[
 \frac{w_+^{2-p}}{p-1}
 \ge
 \frac{\sfm^{2-p}}{p-1}>0.
\]
Together with \eqref{supersol} and the
uniform growth in $\sfK$ established above,
this allows us, after fixing $\Lambda$,
to choose $\sfK$ sufficiently large so that
the spatial contribution is at least
$2\Lambda S$ throughout $Q_{R,S}$,
uniformly in $j$.
Thus, in this case, $w_+$ is a supersolution
to \eqref{PDE2} throughout $Q_{R,S}$
for every $j\in\N$. If $p>2$, then at every point where
$w_+\le2\sfM$ we have
\[
 \frac{w_+^{2-p}}{p-1}
 \ge
 \frac{(2\sfM)^{2-p}}{p-1}>0.
\]
With $\Lambda$ fixed, the same argument
therefore shows that, for sufficiently
large $\sfK$, the supersolution inequality
holds at all such points, uniformly in $j$. Consequently, truncating from above at the
constant solution $2\sfM$, we obtain that
$\min\{w_+,2\sfM\}$ is a weak supersolution
to the regularized divergence equation
\eqref{PDE-reg} throughout $Q_{R,S}$,
for every $j\in\N$. Similarly, after increasing $\sfK$ if necessary,
$w_-$ satisfies the subsolution inequality
for \eqref{PDE2} wherever $w_->0$ if $p\ge2$,
and wherever $w_->\frac12\sfm$ if $1<p<2$.
For $p>2$, the positivity restriction is
needed because the coefficient contains
the negative power $w_-^{2-p}$.
For $p=2$, this coefficient equals one,
so the restriction can be omitted. Consequently, for every $p>1$, the truncated
function $\max\{w_-,\frac12\sfm\}$ is a weak
subsolution to \eqref{PDE-reg} throughout
$Q_{R,S}$, for every $j\in\N$.
We apply the comparison principle for
\eqref{PDE-reg} with
$\max\{w_-,\frac12\sfm\}$ as the subsolution.
As the supersolution, we use $w_+$ if
$1<p\le2$, and $\min\{w_+,2\sfM\}$ if $p>2$.
Both comparison functions are bounded below
by $\frac12\sfm>0$.
The previously established boundary ordering,
together with $\frac12\sfm\le u_j\le2\sfM$,
ensures that these comparison functions
bound $u_j$ on the parabolic boundary
for every $j\ge j_\delta$.
The comparison principle therefore yields
\[
 w_-\le u_j\le w_+
 \qquad\text{in }Q_{R,S},
 \qquad j\ge j_\delta.
\]
Keeping $\delta$, $\Lambda$, and $\sfK$ fixed,
we now pass to the locally uniform limit
as $j\to\infty$.
Since the barriers are independent of $j$,
this gives
\begin{align*}
 u(x_{\mathrm b},t_{\mathrm b})
 &-\Big[2\delta+\sfK f(x)+\Lambda(t-t_{\mathrm b})^2\Big]\\
 &\le w(x,t)\le
 u(x_{\mathrm b},t_{\mathrm b})
 +\Big[2\delta+\sfK f(x)+\Lambda(t-t_{\mathrm b})^2\Big]
\end{align*}
for every $(x,t)\in Q_{R,S}$.
Since $f(x_{\mathrm b})=0$, letting first
$(x,t)\to(x_{\mathrm b},t_{\mathrm b})$
and then $\delta\downarrow0$ gives
\[
 w(x,t)\to u(x_{\mathrm b},t_{\mathrm b})
 \qquad\text{as }
 Q_{R,S}\ni(x,t)\to
 (x_{\mathrm b},t_{\mathrm b}).
\]
Thus, $w$ attains the boundary values $u$
on $\partial B_R\times[-S,0]$.

Next, fix an initial boundary point
$(x_\ast,-S)$ with $x_\ast\in B_R$.
For parameters $\delta\in(0,1)$ and
$\sfK,\Lambda>0$, we define the barriers
\[
 w_\pm(x,t)
 :=
 u(x_\ast,-S)
 \pm\Bigl[
   2\delta+\sfK|x-x_\ast|^q
   +\Lambda(t+S)
 \Bigr],
 \qquad (x,t)\in Q_{R,S},
\]
where $q:=\max\{2,p'\}$ and
$p'=\frac p{p-1}$. We next verify that
\(
 D_\xi\sfa_{\varepsilon_j}(\nabla w_\pm)
 \fprod\nabla^2w_\pm
\)
is bounded in $Q_{R,S}$, uniformly in $j$. If $p\ge2$, then $q=2$, and hence
\[
 \nabla w_\pm=\pm2\sfK(x-x_\ast),
 \qquad
 \nabla^2w_\pm=\pm2\sfK I.
\]
Since $|x-x_\ast|\le2R$ and
$0<\varepsilon_j<1$, the formula for
$D_\xi\sfa_{\varepsilon_j}$ gives a bound
depending only on $N,p,R,\sfK$.
If $1<p<2$, then $q=p'>2$.
For $x\ne x_\ast$, the same formula and
$p-2<0$ give
\begin{align*}
 &\left|
   D_\xi\sfa_{\varepsilon_j}(\nabla w_\pm)
   \fprod\nabla^2w_\pm
 \right|\\
 &\qquad\le
 C|\nabla w_\pm|^{p-2}|\nabla^2w_\pm|
 \le
 C\sfK^{p-1}
 |x-x_\ast|^{(q-1)(p-2)+q-2}=C\sfK^{p-1},
\end{align*}
where $C$ depends only on $N,p$.
Here we used
\[
 (q-1)(p-2)+q-2
 =(q-1)(p-1)-1
 =0.
\]
The estimate extends to $x=x_\ast$ by
continuity, since $\varepsilon_j>0$
and $w_\pm$ is twice continuously
differentiable in space. Thus, in both cases, the bound depends only
on $N,p,R,\sfK$ and is independent of $j$
and $\Lambda$.

Fix $\sfK$ as above, and let
$C=C(N,p,R,\sfK)$ be the uniform bound
for the spatial contraction obtained above. If $p\ge2$, then $w_+\ge\sfm$ gives
\[
 \frac{w_+^{2-p}}{p-1}
 \le
 \frac{\sfm^{2-p}}{p-1}.
\]
Choosing $\Lambda>C\sfm^{2-p}/(p-1)$
therefore makes $w_+$ a supersolution
to \eqref{PDE2} in $Q_{R,S}$
for every $j\in\N$. If $1<p<2$, we have
\[
 w_+(x,t)
 \le
 \sfM+2+\sfK(2R)^q+\Lambda S.
\]
Since $0<2-p<1$, it follows that
\[
    \lim_{\Lambda\to\infty}
 \frac{
   C\bigl(\sfM+2+\sfK(2R)^q+\Lambda S\bigr)^{2-p}
 }{(p-1)\Lambda}
 =0.
\]
Choosing $\Lambda$ sufficiently large
therefore makes $w_+$ a supersolution
to \eqref{PDE2} in $Q_{R,S}$
for every $j\in\N$.

Since $\partial_t w_-=-\Lambda$, the bound
for the spatial term yields the subsolution
inequality for sufficiently large $\Lambda$
wherever $w_-^{2-p}$ is uniformly bounded.
For $p>2$, this factor may diverge as
$w_-\downarrow0$. We therefore work on
$\{w_->\frac12\sfm\}$, where
$\frac12\sfm<w_-\le\sfM$ provides a uniform
bound for every $p>1$. Truncating from below
at the constant solution $\frac12\sfm$
then gives the weak subsolution
$\max\{w_-,\frac12\sfm\}$.

To apply the comparison principle with these
barriers, we check that the previously
established boundary ordering remains valid.
Since $t+S\ge0$, increasing $\Lambda$ decreases
$w_-$ and increases $w_+$.
Moreover, the lower truncation cannot exceed
the boundary values, because
$g_j=u_j\ge\frac12\sfm$ there.
Consequently,
\[
 \max\{w_-,\tfrac12\sfm\}
 \le g_j=u_j\le w_+
 \qquad\text{on }
 \partial_{\mathrm{par}}Q_{R,S},
 \qquad j\ge j_\delta.
\]

The comparison principle for \eqref{PDE-reg}
gives $w_-\le u_j\le w_+$ in $Q_{R,S}$
for every $j\ge j_\delta$.
Keeping $\delta$, $\sfK$, and $\Lambda$ fixed
and passing to the locally uniform limit
as $j\to\infty$, we obtain
\[
 |w(x,t)-u(x_\ast,-S)|
 \le
 2\delta+\sfK|x-x_\ast|^q+\Lambda(t+S)
 \qquad\text{for every }(x,t)\in Q_{R,S}.
\]
Letting first $(x,t)\to(x_\ast,-S)$
within $Q_{R,S}$ and then $\delta\downarrow0$,
we conclude that
\[
 w(x,t)\to u(x_\ast,-S)
 \qquad\text{as }
 Q_{R,S}\ni(x,t)\to(x_\ast,-S).
\]
Thus, $w$ attains the prescribed initial
values on $B_R\times\{-S\}$.
Together with the lateral boundary limits,
this shows that $w$ extends continuously
to $\partial_{\mathrm{par}}Q_{R,S}$
with boundary values $u$, completing the
identification of the initial and boundary
values of the limit map.\\

\noindent\textbf{Step 4: Identification of the limit map.} It remains to show that $w$ is a weak solution
of~\eqref{PDE-reg} in $Q_{R,S}$ and to conclude $w=u$.
Fix a cylinder $Q_{r,s}\Subset Q_{R,S}$.
Since this cylinder is arbitrary, it suffices
to verify the weak formulation of~\eqref{PDE-reg}
in $Q_{r,s}$. By the gradient bound obtained in Step~2 from
Proposition~\ref{nn:prop:reg-gradient}, we have
\[
 \sup_{j\in\N}
 \|\nabla u_j\|_{L^\infty(Q_{r,s},\R^N)}
 <\infty.
\]
Since $0<\varepsilon_j<1$ and
\[
 \bigl|\sfa_{\varepsilon_j}(\nabla u_j)\bigr|
 \le (1+|\nabla u_j|)^{p-1},
\]
the preceding gradient bound also gives
\[
 \sup_{j\in\N}
 \bigl\|\sfa_{\varepsilon_j}(\nabla u_j)\bigr\|_{L^\infty(Q_{r,s},\R^N)}
 <\infty.
\]
By weak compactness on $Q_{r,s}$, we can pass
to a subsequence, still indexed by $j$,
such that
\begin{equation}\label{convergence-Duj}
\left\{
 \begin{array}{cl}
  \nabla u_j
  \wto \nabla w
  & \text{in }L^p(Q_{r,s},\R^N),\\[6pt]
  \sfa_{\varepsilon_j}(\nabla u_j)
  \wto \sfa_o
  &\text{in }L^{p'}(Q_{r,s},\R^N)
 \end{array}
\right.
\end{equation}
as $j\to\infty$, for some vector field
$\sfa_o\in L^{p'}(Q_{r,s},\R^N)$.
The first weak limit is $\nabla w$ because
the locally uniform convergence $u_j\to w$
implies that the spatial gradients converge
to $\nabla w$ in the sense of distributions. Since $u_j^{p-1}\to w^{p-1}$ locally uniformly in $Q_{r,s}$
and $\sfa_{\varepsilon_j}(\nabla u_j)\wto\sfa_o$
in $L^{p'}(Q_{r,s},\R^N)$ by \eqref{convergence-Duj},
we may pass to the limit in the weak formulation
of \eqref{PDE-reg} to obtain
\begin{equation}\label{PDE-w}
 \partial_t(w^{p-1})=\Div\sfa_o
 \qquad\text{weakly in }Q_{r,s}.
\end{equation}
To conclude that $w$ is a weak solution of \eqref{PDE}
in $Q_{r,s}$, it remains to show that
\[
 \sfa_o=|\nabla w|^{p-2}\nabla w
 \qquad\text{almost everywhere in }Q_{r,s}.
\]
Let $\varphi\in C^\infty_0(Q_{r,s})$ be nonnegative,
and let $t_-$ denote the lower endpoint of the time
interval of $Q_{r,s}$. For $h>0$, we define the
exponential time mollification of $w$ by
\[
 \lkk w\rkk_h(x,t)
 :=
 \frac1h\int_{t_-}^{t}
 e^\frac{\tau-t}{h}
 w(x,\tau)\,\d\tau,
 \qquad (x,t)\in Q_{r,s}.
\]
Differentiating with respect to $t$, we obtain
\begin{equation}\label{nn:eq:time-mollification-derivative}
 \partial_t\lkk w\rkk_h
 =
 \frac{w-\lkk w\rkk_h}{h}.
\end{equation}
For fixed $h>0$, this time derivative is bounded
because $w$ is bounded. Moreover,
$\nabla\lkk w\rkk_h=\lkk\nabla w\rkk_h$ belongs to
$L^p(Q_{r,s},\R^N)$. Since $\varphi$ has compact
support in $Q_{r,s}$, we can test \eqref{PDE-w} with
$\varphi\lkk w\rkk_h$. Expanding the time derivative and
the spatial gradient of this product gives
\begin{align*}
 &\iint_{Q_{r,s}}
   w^{p-1}\partial_t\lkk w\rkk_h\,\varphi\,\dx\,\dt
 +\iint_{Q_{r,s}}
   w^{p-1}\lkk w\rkk_h\,\partial_t\varphi\,\dx\,\dt\\
 &\qquad=
 \iint_{Q_{r,s}}
   \varphi\,\sfa_o\cdot\nabla\lkk w\rkk_h\,\dx\,\dt
 +\iint_{Q_{r,s}}
   \lkk w\rkk_h\,\sfa_o\cdot\nabla\varphi\,\dx\,\dt.
\end{align*}
In the first integral on the left-hand side, we add
and subtract $\lkk w\rkk_h^{p-1}$ and use
\eqref{nn:eq:time-mollification-derivative}. This gives
\begin{align*}
 &\iint_{Q_{r,s}}
   w^{p-1}\partial_t\lkk w\rkk_h\,\varphi\,\dx\,\dt\\
 &\quad=
 \iint_{Q_{r,s}}
   \lkk w\rkk_h^{p-1}
   \partial_t\lkk w\rkk_h\,\varphi\,\dx\,\dt
 +\frac1h\iint_{Q_{r,s}}
   \underbrace{
     \Bigl(w^{p-1}-\lkk w\rkk_h^{p-1}\Bigr)
     \bigl(w-\lkk w\rkk_h\bigr)
   }_{\ge0}
   \varphi\,\dx\,\dt\\
 &\quad\ge
 \iint_{Q_{r,s}}
   \lkk w\rkk_h^{p-1}
   \partial_t\lkk w\rkk_h\,\varphi\,\dx\,\dt=
 -\frac1p\iint_{Q_{r,s}}
   \lkk w\rkk_h^p\,\partial_t\varphi\,\dx\,\dt.
\end{align*}
The indicated product is nonnegative because
$s\mapsto s^{p-1}$ is increasing on $[0,\infty)$.
Together with $\varphi\ge0$, this gives the inequality.
In the last equality, we integrated by parts in time.
Substituting this lower bound into the preceding
identity and moving the term containing $\nabla\varphi$
to the right-hand side, we obtain
\begin{align*}
 \iint_{Q_{r,s}}&
   \varphi\,\sfa_o\cdot\nabla\lkk w\rkk_h\,\dx\,\dt\\
 &\ge
 \iint_{Q_{r,s}}
   \left(
     w^{p-1}\lkk w\rkk_h
     -\tfrac1p \lkk w\rkk_h^p
   \right)\partial_t\varphi\,\dx\,\dt
 -\iint_{Q_{r,s}}
   \lkk w\rkk_h\,\sfa_o\cdot\nabla\varphi\,\dx\,\dt.
\end{align*}
As $h\downarrow0$, we have
\[
 \lkk w\rkk_h\to w
 \quad\text{in }L^p(Q_{r,s}),
 \qquad
 \nabla\lkk w\rkk_h\to\nabla w
 \quad\text{in }L^p(Q_{r,s},\R^N),
\]
with strong convergence in both spaces. In particular,
the coefficient of $\partial_t\varphi$ converges in
$L^1(Q_{r,s})$ to $\frac{p-1}{p}w^p$.
Since $\sfa_o\in L^{p'}(Q_{r,s},\R^N)$, we may also
pass to the limit in the two integrals involving $\sfa_o$.
We conclude that
\begin{equation}\label{Bound-w}
 \iint_{Q_{r,s}}
   \varphi\,\sfa_o\cdot\nabla w\,\dx\,\dt
 \ge
 \frac{p-1}{p}\iint_{Q_{r,s}}
   w^p\partial_t\varphi\,\dx\,\dt
 -
 \iint_{Q_{r,s}}
   w\,\sfa_o\cdot\nabla\varphi\,\dx\,\dt
\end{equation}
for every nonnegative $\varphi\in C^\infty_0(Q_{r,s})$. Since $u_j$ is smooth, testing \eqref{PDE-reg} with
$\varphi u_j$ and integrating by parts gives
\begin{align}\label{Bound-uj}
 &\iint_{Q_{r,s}}
   \varphi\,\sfa_{\varepsilon_j}(\nabla u_j)
   \cdot\nabla u_j\,\dx\,\dt\notag\\
 &\qquad=
 \frac{p-1}{p}\iint_{Q_{r,s}}
   u_j^p\,\partial_t\varphi\,\dx\,\dt
 -
 \iint_{Q_{r,s}}
   u_j\,\sfa_{\varepsilon_j}(\nabla u_j)
   \cdot\nabla\varphi\,\dx\,\dt.
\end{align}
Since $Q_{r,s}\Subset Q_{R,S}$, the locally uniform convergence
$u_j\to w$ implies that $u_j^p\to w^p$ uniformly on $Q_{r,s}$.
Moreover, $u_j\nabla\varphi\to w\nabla\varphi$ strongly in
$L^p(Q_{r,s},\R^N)$. Together with the weak convergence of
$\sfa_{\varepsilon_j}(\nabla u_j)$ in \eqref{convergence-Duj},
this allows us to pass to the limit on the right-hand side
of \eqref{Bound-uj}. We obtain
\begin{align*}
 &\lim_{j\to\infty}\iint_{Q_{r,s}}
   \varphi\,\sfa_{\varepsilon_j}(\nabla u_j)
   \cdot\nabla u_j\,\dx\,\dt\\
 &\qquad=
 \frac{p-1}{p}\iint_{Q_{r,s}}
   w^p\,\partial_t\varphi\,\dx\,\dt
 -
 \iint_{Q_{r,s}}
   w\,\sfa_o\cdot\nabla\varphi\,\dx\,\dt\\
 &\qquad\le
 \iint_{Q_{r,s}}
   \varphi\,\sfa_o\cdot\nabla w\,\dx\,\dt,
\end{align*}
where the last inequality follows from \eqref{Bound-w}
since $\varphi\ge0$. Fix $\xi\in L^p(Q_{r,s},\R^N)$ and a nonnegative
$\varphi\in C^\infty_0(Q_{r,s})$.
Since $\sfa_{\varepsilon_j}(\xi)\to|\xi|^{p-2}\xi$
almost everywhere and
\( |\sfa_{\varepsilon_j}(\xi)|
 \le (1+|\xi|)^{p-1} \),
dominated convergence gives
\(  \sfa_{\varepsilon_j}(\xi)\to|\xi|^{p-2}\xi\)
strongly in \(L^{p'}(Q_{r,s},\R^N)\).
Combining this convergence with \eqref{convergence-Duj}
and the preceding estimate for the term containing
$\sfa_{\varepsilon_j}(\nabla u_j)\cdot\nabla u_j$,
we obtain
\begin{align*}
 &\iint_{Q_{r,s}}
   \varphi\bigl(\sfa_o-|\xi|^{p-2}\xi\bigr)
   \cdot(\nabla w-\xi)\,\dx\,\dt\\
 &\qquad\ge
 \lim_{j\to\infty}
 \iint_{Q_{r,s}}
   \varphi\bigl(
     \sfa_{\varepsilon_j}(\nabla u_j)
     -\sfa_{\varepsilon_j}(\xi)
   \bigr)
   \cdot(\nabla u_j-\xi)\,\dx\,\dt
 \ge 0.
\end{align*}
The last inequality follows from the monotonicity of
$\sfa_{\varepsilon_j}$ and the nonnegativity of $\varphi$. Let $e_\ell$ denote the $\ell$-th standard basis vector of
$\R^N$. For $\lambda>0$, we apply the preceding inequality
first with $\xi=\nabla w-\lambda e_\ell$ and then with
$\xi=\nabla w+\lambda e_\ell$. Dividing by $\lambda$ gives
\begin{align*}
 &\int_{Q_{r,s}}
   \varphi\,|\nabla w-\lambda e_\ell|^{p-2}
   (\nabla w-\lambda e_\ell)\cdot e_\ell\,\dx\,\dt\\
 &\qquad\le
 \int_{Q_{r,s}}
   \varphi\,\sfa_o\cdot e_\ell\,\dx\,\dt\\
 &\qquad\le
 \int_{Q_{r,s}}
   \varphi\,|\nabla w+\lambda e_\ell|^{p-2}
   (\nabla w+\lambda e_\ell)\cdot e_\ell\,\dx\,\dt.
\end{align*}
For $0<\lambda\le1$, the first and last integrands are
bounded in absolute value by $\varphi(1+|\nabla w|)^{p-1}$,
which is integrable on $Q_{r,s}$. By continuity of the map
$\zeta\mapsto|\zeta|^{p-2}\zeta$, we may let $\lambda\downarrow0$
and apply dominated convergence to obtain
\[
 \int_{Q_{r,s}}
   \varphi\bigl(\sfa_o-|\nabla w|^{p-2}\nabla w\bigr)
   \cdot e_\ell\,\dx\,\dt=0.
\]
Since this holds for every nonnegative
$\varphi\in C^\infty_0(Q_{r,s})$ and every
$\ell=1,\ldots,N$, we conclude that
\[
 \sfa_o=|\nabla w|^{p-2}\nabla w
 \qquad\text{almost everywhere in }Q_{r,s}.
\]
Substituting this identity into \eqref{PDE-w} shows that
$w$ is a weak solution of \eqref{PDE} in $Q_{r,s}$.
Since $Q_{r,s}\Subset Q_{R,S}$ was arbitrary, $w$ is
a weak solution in $Q_{R,S}$. 

It remains to show that $w=u$. Since the energy bounds
for $w$ are local, we apply the comparison principle
on the cylinders
\( Q_n
 :=
 B_{R-\frac1n}(0)
 \times\left(-S+\frac1n,-\frac1n\right) \),
where $n\in\N$ is sufficiently large that
$\frac1n<R$ and $\frac2n<S$. Then
$Q_n\Subset Q_{R,S}$ and $Q_n\uparrow Q_{R,S}$. By Step~3, $w-u$ extends continuously to
$\partial_{\mathrm{par}}Q_{R,S}$ and vanishes there.
Since this boundary is compact and every point of
$\partial_{\mathrm{par}}Q_n$ lies within distance $1/n$
of it, we obtain
\[
 \delta_n:=\sup_{\partial_{\mathrm{par}}Q_n}|w-u|
 \longrightarrow 0
 \qquad\text{as }n\to\infty.
\]
For $n$ sufficiently large, we have $\delta_n<\sfm$.
Since $u\ge\sfm$, the definition of $\delta_n$ gives
\[
 \left(1-\frac{\delta_n}{\sfm}\right)u
 \le w\le
 \left(1+\frac{\delta_n}{\sfm}\right)u
 \qquad\text{on }\partial_{\mathrm{par}}Q_n.
\]
Both factors are positive. By the homogeneity of
\eqref{PDE}, the two multiples of $u$ are therefore
weak solutions of \eqref{PDE}. Since $Q_n\Subset Q_{R,S}$,
the comparison principle from \cite{BS} applies on $Q_n$
and yields
\[
 \left(1-\frac{\delta_n}{\sfm}\right)u
 \le w\le
 \left(1+\frac{\delta_n}{\sfm}\right)u
 \qquad\text{in }Q_n.
\]

Fix $(x,t)\in Q_{R,S}$. Then $(x,t)\in Q_n$ for all
sufficiently large $n$. Letting $n\to\infty$ in the
preceding inequalities gives $w(x,t)=u(x,t)$.
Thus $w=u$ in $Q_{R,S}$, which completes the proof.
\end{proof}

\begin{proposition}[Positive cylinders]
\label{nn:prop:positive}
Let $r,h>0$, and let $u$ be a continuous weak
solution to \eqref{PDE} in a neighborhood of
$\overline{B_{2r}(x_o)}\times[s-2h,s]$.
Suppose that
\[
 0<\sfm\le u\le\sfM
 \qquad\text{in }B_{2r}(x_o)\times(s-2h,s].
\]
Then, for every $k\ge k_o$, with $k_o$ defined
in \eqref{nn:eq:k}, we have
\begin{equation}\label{nn:eq:positive}
 \big\|\nabla\bigl(u^{\frac1k}\bigr)\big\|_
 {L^\infty(B_r(x_o)\times(s-h,s))}
 \le
 C\sfM^{\frac1k}
 \left(\frac1r+\frac1{h^{\frac1p}}\right),
\end{equation}
where $C=C(N,p,k)$ is independent of
$\sfm,\sfM,r,h$.
\end{proposition}

\begin{proof}
By Lemma~\ref{lem:approximation}, there exist
$\varepsilon_j\downarrow0$ and smooth solutions
$u_j$ to \eqref{PDE-reg} in
$B_{2r}(x_o)\times(s-2h,s]$ such that \(  \tfrac12\sfm\le u_j\le2\sfM\)
and $u_j\to u$ locally uniformly in this cylinder.
Applying \eqref{nn:eq:reg-est} and absorbing
the fixed factors of $2$ into the constant,
we obtain
\[
 \big\|\nabla\bigl(u_j^{\frac1k}\bigr)\big\|_
 {L^\infty(B_r(x_o)\times(s-h,s))}
 \le
 C\sfM^{\frac1k}
 \left(\frac1r+\frac1{h^{\frac1p}}\right)
 +C\varepsilon_j\sfm^{-\frac{k-1}{k}}.
\]
Since $\sfm>0$ is fixed, the last term tends
to zero as $j\to\infty$.
Moreover, $u_j^{\frac1k}\to u^{\frac1k}$
locally uniformly in
$B_{2r}(x_o)\times(s-2h,s]$.
Thus, passing to the limit in the spatial
Lipschitz estimates gives
\[
 \big|
 u(x,t)^{\frac1k}-u(y,t)^{\frac1k}
 \big|
 \le
 C\sfM^{\frac1k}
 \left(\frac1r+\frac1{h^{\frac1p}}\right)
 |x-y|
\]
for every $x,y\in B_r(x_o)$ and
$t\in(s-h,s)$.
This implies \eqref{nn:eq:positive}.
\end{proof}

\begin{remark}[Approximation and order of limits]\upshape
\label{nn:rem:order}
For $p\ne2$, adding a positive constant to $u$
does not generally preserve \eqref{PDE}:
the spatial gradient remains unchanged,
but the time derivative of $u^{p-1}$ changes.
Therefore, we use solutions to the regularized
equation \eqref{PDE-reg}, rather than $u+\varepsilon$.
The regularized estimate contains the error term
\[
 \frac{C\varepsilon_j}{\sfm^{\frac{k-1}{k}}},
\]
which tends to zero as $\varepsilon_j\downarrow0$
for fixed $\sfm>0$, but is not uniform as
$\sfm\downarrow0$.
We first pass to the regularization limit on
each fixed positive cylinder.
The resulting estimate \eqref{nn:eq:positive}
is independent of $\sfm$ and can then be used
on cylinders whose positive lower bounds
approach zero.
\end{remark}

\section{Proof of Theorem~\ref{thm:Lipschitz}}
In this section, we first estimate how quickly
a nonnegative solution can become positive
inside a ball where it initially vanishes.
Comparison with an exponential barrier shows
that, for short times, the solution remains
exponentially small in a smaller concentric ball.
The precise result is stated in the following lemma.

We then prove local spatial Lipschitz continuity by combining this barrier estimate with the a priori gradient estimate on positive cylinders from Proposition~\ref{nn:prop:positive}. While the gradient estimate applies directly away from the zero set, the exponential smallness provided by the barrier estimate ensures that the gradient remains uniformly bounded as the solution becomes positive. In fact, we obtain a slightly stronger estimate showing that the gradient vanishes as the solution approaches the zero set; see Theorem~\ref{thm:lip-k}.

\subsection{Barrier estimates close to the zero time}\label{subsec:exp-barrier}

\begin{lemma}\label{lem:barrier}
There exist constants $c_o,\tau_o>0$, depending
only on $N,p$, with the following property.
Suppose that $u$ is a nonnegative weak solution
to \eqref{PDE} in
\( B_\rho(x_o)\times(t_\ast,t_\ast+\tau_o\rho^p]\),
satisfies $u\le\sfM$, and is continuous up to
the initial time, with
\[
 u(x,t_\ast)=0
 \qquad\text{for every }x\in B_\rho(x_o).
\]
Then
\begin{equation}\label{nn:eq:onset}
 u(x,t)
 \le
 \sfM\exp\left[
   -c_o\left(
     \frac{\rho^p}{t-t_\ast}
   \right)^{\frac1{p-1}}
 \right]
\end{equation}
for any 
\[
 (x,t)\in
 B_{\rho/2}(x_o)\times
 \left(t_\ast,t_\ast+\tfrac12 \tau_o\rho^p\right].
\]
If the solution is defined only up to an earlier
terminal time $T$, the same estimate holds
for $t<T$ in the indicated time interval.
It also holds at $t=T$ whenever $u$ extends
continuously to that time.
\end{lemma}

\begin{proof}
If $\sfM=0$, the assertion is immediate.
Therefore, we assume $\sfM>0$.
The change of variables
\[
 y=\frac{x-x_o}{\rho},
 \qquad
 \tau=\frac{t-t_\ast}{\rho^p}
\]
maps the cylinder
$B_\rho(x_o)\times(t_\ast,t_\ast+\tau_o\rho^p]$
onto \(  \widetilde Q:=B_1\times(0,\tau_o]\).
By Lemma~\ref{lem:scaling} the rescaled function
\[
 v(y,\tau)
 :=
 \frac{ u\bigl(x_o+\rho y,t_\ast+\rho^p\tau\bigr)}{\sfM}
\]
is a nonnegative weak solution to
\[
 \partial_\tau(v^{p-1})
 -
 \operatorname{div}
 \bigl(|\nabla v|^{p-2}\nabla v\bigr)
 =0
 \qquad\text{in }\widetilde Q.
\]
Moreover, $0\le v\le1$ in $\widetilde Q$,
and $v$ is continuous up to the initial time,
with $v(\cdot,0)=0$.
Thus, proving \eqref{nn:eq:onset} reduces to
proving
\begin{equation}\label{nn:eq:onset-normalized}
 v(y,\tau)
 \le
 \exp\left[-c_o\tau^{-\frac1{p-1}}\right]
 \qquad\text{in }
 B_{1/2}\times\left(0,
 \tfrac12 \tau_o\right].
\end{equation}
In the following, all spatial derivatives
are taken with respect to $y$.
Let $p'=\frac p{p-1}>1$ be the Hölder conjugate
of $p$. We define
\[
 \psi(y):=(1-|y|^{p'})^{p'}
 \qquad\text{for }y\in\overline{B_1}.
\]
For $0<|y|<1$, we have
\[
 \nabla\psi(y)
 =
 -(p')^2(1-|y|^{p'})^{p'-1}|y|^{p'-2}y.
\]
Since $(p'-1)p=p'$, we obtain
\begin{align}\label{nn:eq:psi-grad}
 |\nabla\psi|^p
 &=
 (p')^{2p}
 (1-|y|^{p'})^{(p'-1)p}
 |y|^{(p'-1)p}\notag\\
 &=
 (p')^{2p}|y|^{p'}\psi
 \le (p')^{2p}\psi.
\end{align}
Moreover, using $(p'-1)(p-1)=1$, we find
\begin{align}\label{nn:eq:psi-flux}
 |\nabla\psi|^{p-2}\nabla\psi
 &=
 -(p')^{2(p-1)}
 (1-|y|^{p'})^{(p'-1)(p-1)}
 |y|^{(p'-1)(p-1)-1}y\notag\\
 &=
 -(p')^{2(p-1)}(1-|y|^{p'})y.
\end{align}
Since $p'>1$, the last expression extends to
a $C^1$ vector field on $B_1$. Taking its divergence yields
\begin{align*}
 \Delta_p\psi
 &=
 -(p')^{2(p-1)}
 \diver\bigl((1-|y|^{p'})y\bigr)\\
 &=
 -(p')^{2(p-1)}
 \bigl(N(1-|y|^{p'})-p'|y|^{p'}\bigr)\\
 &=
 (p')^{2(p-1)}
 \bigl((N+p')|y|^{p'}-N\bigr)
 \qquad\text{in }B_1.
\end{align*}
To obtain a lower bound for $\Delta_p\psi$
in terms of $\psi$, we distinguish two cases.
If
\[
 |y|^{p'}\ge\frac N{N+p'},
\]
then $\Delta_p\psi\ge0$.
In the complementary region,
\[
 1-|y|^{p'}>\frac{p'}{N+p'},
 \qquad\text{and hence}\qquad
 \psi(y)>\left(\frac{p'}{N+p'}\right)^{p'}.
\]
Since $\Delta_p\psi\ge-N(p')^{2(p-1)}$, we obtain
\[
 \Delta_p\psi
 \ge -N(p')^{2(p-1)}
 \ge
 -N(p')^{2(p-1)}
 \left(\frac{N+p'}{p'}\right)^{p'}\psi.
\]
Combining the two regions yields
\begin{equation}\label{nn:eq:psi-lap}
 \Delta_p\psi\ge-C_o\psi
 \qquad\text{in }B_1,
 \qquad
 C_o:=
 N(p')^{2(p-1)}
 \left(\frac{N+p'}{p'}\right)^{p'}.
\end{equation}
The constant $C_o$ depends only on $N,p$.

Using the preceding estimates for $\psi$,
we now construct a barrier.
For $\lambda>0$ and $\delta\in(0,\tfrac12\tau_o]$, define
\[
 f(y,\tau)
 :=
 -\frac{\lambda\psi(y)}{(\tau+\delta)^{\frac1{p-1}}},
 \qquad
 V_\delta(y,\tau):=\exp(f(y,\tau)).
\]
The parameter $\lambda$ will be chosen below,
independently of $\delta$.
Since $\psi\in C^1(\overline{B_1})$,
both $f$ and $V_\delta$ are continuously
differentiable in space. Moreover,
\[
 |\nabla f|^{p-2}\nabla f
 =
 -\frac{\lambda^{p-1}|\nabla\psi|^{p-2}\nabla\psi}{\tau+\delta}
 .
\]
By \eqref{nn:eq:psi-flux}, this flux is
continuously differentiable in space.
Thus, the divergences below are well defined
pointwise, including at the origin $y=0$.
The chain rule gives
\[
 \partial_\tau(V_\delta^{p-1})
 =
 (p-1)\partial_\tau f\,e^{(p-1)f}.
\]
Similarly, since $\nabla V_\delta=e^f\nabla f$,
we have
\[
 |\nabla V_\delta|^{p-2}\nabla V_\delta
 =
 e^{(p-1)f}|\nabla f|^{p-2}\nabla f.
\]
Taking the divergence and applying the
product rule, we obtain
\begin{align*}
 \Delta_pV_\delta
 &=
 \nabla(e^{(p-1)f})\cdot
 \bigl(|\nabla f|^{p-2}\nabla f\bigr)
 +e^{(p-1)f}\Delta_pf\\
 &=
 e^{(p-1)f}
 \Bigl[
   (p-1)|\nabla f|^p+\Delta_pf
 \Bigr].
\end{align*}
Consequently,
\begin{equation}\label{nn:eq:barrier-identity}
 \partial_\tau(V_\delta^{p-1})
 -\Delta_pV_\delta
 =
 V_\delta^{p-1}\sfR,
\end{equation}
where
\[
 \sfR
 :=
 (p-1)\partial_\tau f
 -\Delta_pf
 -(p-1)|\nabla f|^p.
\]
By \eqref{nn:eq:barrier-identity},
$V_\delta$ is a supersolution whenever
$\sfR\ge0$.
We compute the three terms in $\sfR$
separately. Since $1+\frac1{p-1}=p'$, differentiation
with respect to time gives
\[
 (p-1)\partial_\tau f
 =
 \frac{\lambda\psi}
      {(\tau+\delta)^{\frac1{p-1}+1}}
 =
 \frac{\lambda\psi}{(\tau+\delta)^{p'}}.
\]
For the gradient term, we obtain
\[
 (p-1)|\nabla f|^p
 =
 \frac{(p-1)\lambda^p|\nabla\psi|^p}
      {(\tau+\delta)^{\frac p{p-1}}}
 =
 \frac{(p-1)\lambda^p|\nabla\psi|^p}
      {(\tau+\delta)^{p'}}.
\]
Finally, taking the spatial divergence of
\[
 |\nabla f|^{p-2}\nabla f
 =
 -\frac{\lambda^{p-1}}{\tau+\delta}
 |\nabla\psi|^{p-2}\nabla\psi
\]
yields
\[
 \Delta_pf
 =
 -\frac{\lambda^{p-1}}{\tau+\delta}
 \Delta_p\psi.
\]
Substituting these expressions into the
definition of $\sfR$, we find
\begin{align*}
 \sfR
 &=
 \frac{\lambda\psi}{(\tau+\delta)^{p'}}
 +\frac{\lambda^{p-1}}{\tau+\delta}
   \Delta_p\psi
 -\frac{(p-1)\lambda^p}{(\tau+\delta)^{p'}}
   |\nabla\psi|^p.
\end{align*}
Using \eqref{nn:eq:psi-grad} and
\eqref{nn:eq:psi-lap}, we obtain
\begin{equation}\label{nn:eq:barrier-residual}
 \sfR
 \ge
 \psi\Biggl[
   \bigl[1-(p-1)(p')^{2p}\lambda^{p-1}\bigr]
   \frac{\lambda}{(\tau+\delta)^{p'}}
   -\frac{C_o\lambda^{p-1}}{\tau+\delta}
 \Biggr].
\end{equation}
We now choose $\lambda$ and $\tau_o$ so that
the right-hand side of
\eqref{nn:eq:barrier-residual} is nonnegative.
First, we choose $\lambda$ small enough to make
the first term positive. We then restrict
$\tau+\delta$ by choosing $\tau_o$ sufficiently
small so that this positive term dominates
the negative one.
Both choices will depend only on $N,p$. The precise choices are as follows.
Specifically, let $\lambda>0$ satisfy
\[
 \lambda^{p-1}
 =
 \frac1{2(p-1)(p')^{2p}}.
\]
With this choice of $\lambda$, the expression
in the inner square brackets in
\eqref{nn:eq:barrier-residual} equals $\frac12$.
Using
\[
 \frac1{\tau+\delta}
 =
 \frac{(\tau+\delta)^{\frac1{p-1}}}
      {(\tau+\delta)^{p'}},
\]
we can therefore rewrite
\eqref{nn:eq:barrier-residual} as
\[
 \sfR
 \ge
 \frac{\psi}{(\tau+\delta)^{p'}}
 \left[
   \tfrac12 \lambda
   -C_o\lambda^{p-1}
    (\tau+\delta)^{\frac1{p-1}}
 \right].
\]
To make the expression in brackets nonnegative,
we choose $\tau_o>0$ such that
\[
 \tau_o^{\frac1{p-1}}
 =
 \frac{\lambda^{2-p}}{2C_o}.
\]
Indeed, whenever $0\le\tau$ and
$\tau+\delta\le\tau_o$, we have
\[
 C_o\lambda^{p-1}
 (\tau+\delta)^{\frac1{p-1}}
 \le
 C_o\lambda^{p-1}\tau_o^{\frac1{p-1}}
 =
 \tfrac12 \lambda.
\]
Thus, $\sfR\ge0$ for $y\in B_1$ and
$0\le\tau\le\frac12\tau_o$, 
provided that $0<\delta\le\frac12\tau_o$. 
Both $\lambda$ and $\tau_o$ depend only on $N,p$.
We conclude that
\begin{equation*}
 \partial_\tau(V_\delta^{p-1})
 -\Delta_pV_\delta
 \ge0
 \qquad\text{in }
 B_1\times(0,\tfrac12\tau_o].
\end{equation*}
Moreover, since $\psi\ge0$ on $\overline{B_1}$
and $\psi=0$ on $\partial B_1$, we have
\[
\mbox{
 \(0<V_\delta\le1\) 
 on 
\(\overline{B_1}\times[0,\tau_o)\)},
 \qquad
 \mbox{\(V_\delta=1\)
 on \(\partial B_1\times[0,\tau_o)\)}.
\]

We now use the comparison principle for
\eqref{PDE} (see \cite[Theorem~1.2]{BS})
to compare $v$ with the barrier.
Since no continuity of $v$ up to $\partial B_1$
is assumed, we first work in $B_{1-h}$,
where $0<h<1$, and subsequently let $h\downarrow0$. On the boundary of the smaller ball $B_{1-h}$,
the function $V_\delta$ is strictly less than one,
so it need not dominate $v$ there.
Therefore, we multiply $V_\delta$ by a suitable
positive constant so that the resulting function
is at least one on this boundary.
The homogeneity of the equation ensures that
this multiplication preserves the supersolution
property. The value of $\psi$ on $\partial B_{1-h}$ is
\[
 \psi_h:=\bigl(1-(1-h)^{p'}\bigr)^{p'}.
\]
For $h\in(0,1)$ and $\delta>0$ define
\[
 W(y,\tau)
 :=
 \exp\left(
   \frac{\lambda\psi_h}{\delta^{\frac1{p-1}}}
 \right)V_\delta(y,\tau).
\]
We have
\begin{align*}
 \partial_\tau(W^{p-1})-\Delta_pW
 &=
 \exp\left(
   \frac{(p-1)\lambda\psi_h}
        {\delta^{\frac1{p-1}}}
 \right)
 \Bigl[
   \partial_\tau(V_\delta^{p-1})
   -\Delta_pV_\delta
 \Bigr]
 \ge0
\end{align*}
in $B_1\times(0,\frac12\tau_o]$. 
Fix $0<\delta<\frac12\tau_o$ and $0<h<1$.
We compare $v$ and $W$ on $B_{1-h}\times\left(0,\tfrac12\tau_o\right)$. 
The closure of this cylinder lies inside
the region where $v$ is continuous and
$W$ is a supersolution.
On the lateral boundary, where $|y|=1-h$,
we have
\[
 W(y,\tau)
 =
 \exp\left[
   \lambda\psi_h
   \left(
     \frac1{\delta^{\frac1{p-1}}}
     -\frac1{(\tau+\delta)^{\frac1{p-1}}}
   \right)
 \right]
 \ge1\ge v(y,\tau).
\]
To obtain the ordering at the initial time, observe that
\[
v(\cdot,0)=0<W(\cdot,0)
\qquad\text{in } B_{1-h}.
\]
Together with the ordering on the lateral boundary,
this shows that $v\le W$ on the parabolic boundary of $B_{1-h}\times(0,\tfrac12\tau_o]$.
Moreover, $W$  satisfies
\[
W\ge \exp\left(
-\frac{\lambda}{\delta^{\frac1{p-1}}}
\right)>0
\]
throughout the comparison cylinder.
Since $0\le v\le1$, all assumptions of the
comparison principle \cite[Theorem~1.2]{BS}
are satisfied. Consequently,
\[
v\le W
\qquad\text{in }
B_{1-h}\times\left(0,\tfrac12\tau_o\right]
\]
whence
\[
v(y,\tau)
\le
\exp\left(
\frac{\lambda\psi_h}{\delta^{\frac1{p-1}}}
\right)V_\delta(y,\tau)
\qquad\text{in }
B_{1-h}\times\left(0,\tfrac12\tau_o\right].
\]

\begin{figure}[htbp]
\centering
\begin{tikzpicture}[
  x=1cm,
  y=1cm,
  >=Latex,
  every node/.style={font=\small}
]
\fill[blue!5] (-4,0) rectangle (4,5.3);
\fill[orange!10] (-2.6,0) rectangle (2.6,4);

\draw[blue!65!black,thick]
  (-4,0) rectangle (4,5.3);

\draw[orange!80!black,thick,dashed]
  (-2.6,4) -- (2.6,4);

\draw[orange!85!black,very thick]
  (-2.6,4) -- (-2.6,0) -- (2.6,0) -- (2.6,4);

\draw[->] (-4.7,-0.4) -- (-4.7,5.9)
  node[above] {$\tau$};
\draw[->] (-4.3,-0.55) -- (4.7,-0.55)
  node[right] {$y_1$};

\foreach \height/\labeltext in {
  0/{0},
  4/{\frac12\tau_o},
  5.3/{\tau_o-\delta}
}{
  \draw (-4.78,\height) -- (-4.62,\height);
  \node[left] at (-4.82,\height) {$\labeltext$};
}

\foreach \position/\labeltext in {
  -4/{-1},
  -2.6/{-(1-h)},
  0/{0},
  2.6/{1-h},
  4/{1}
}{
  \draw (\position,-0.63) -- (\position,-0.47);
  \node[below] at (\position,-0.67) {$\labeltext$};
}

\node[blue!65!black,align=center] at (0,4.7)
  {$B_1\times(0,\tau_o-\delta]$\\
   $V_\delta$ is a supersolution};

\node[align=center] at (0,2.5)
  {$B_{1-h}\times\left(0,\tfrac12\tau_o\right)$\\[3pt]
   Comparison principle\\[2pt]
   $\Longrightarrow\ v\le W$};

\node[
  anchor=west,
  text=orange!70!black
] (lateral-label) at (4.35,2.5)
  {$v\le1\le W$};

\draw[
  ->,
  thick,
  orange!70!black,
  shorten >=3pt
]
  (lateral-label.west) -- (2.6,2.5);

\node[
  text=orange!70!black
] (initial-label) at (0,0.7)
  {$v(\cdot,0)=0\le W(\cdot,0)$};

\draw[
  ->,
  thick,
  orange!70!black,
  shorten >=3pt
]
  (initial-label.south) -- (0,0);

\node[
  draw=blue!50!black,
  rounded corners,
  fill=blue!3,
  align=center,
  text width=10cm,
  inner sep=7pt,
  anchor=north
] at (0,-1.5)
  {Lift the barrier by a constant factor:\par\medskip
   $\displaystyle
     W=\exp\!\left(
       \frac{\lambda\,\psi_h}{\delta^{\frac1{p-1}}}
     \right)V_\delta,
     \qquad
     \psi_h=\bigl(1-(1-h)^{p'}\bigr)^{p'},
     \quad p'=\tfrac{p}{p-1}.$\par\medskip
   This ensures $W\ge1$ on the lateral boundary
   of the smaller cylinder.};

\end{tikzpicture}
\caption{A spatial cross-section of the comparison
cylinders. The solid orange boundary is the
parabolic boundary on which $v\le W$ is established.
The dashed upper time face is not part of the
parabolic boundary. The diagram is schematic,
with $0<\frac12\tau_o<\tau_o-\delta$.}
\label{fig:barrier-comparison}
\end{figure}

Fix $(y,\tau)\in B_1\times(0,\frac12\tau_o]$. 
For $0<h<1-|y|$, we have $y\in B_{1-h}$,
so the preceding estimate gives
\begin{align*}
 v(y,\tau)
 &\le
 \exp\left(
   \frac{\lambda\psi_h}{\delta^{\frac1{p-1}}}
 \right)V_\delta(y,\tau).
\end{align*}
Since $\psi_h\to0$ as $h\downarrow0$ and
$\delta$ remains fixed, the exponential factor
converges to one. Thus, $v(y,\tau)\le V_\delta(y,\tau)$.
As $(y,\tau)$ was arbitrary, this holds throughout
$B_1\times(0,\frac12\tau_o]$.
We now let $\delta\downarrow0$ and conclude that
\[
 v(y,\tau)
 \le
 \exp\left[
   -\frac{\lambda\psi(y)}
         {\tau^{\frac1{p-1}}}
 \right].
\]
For $y\in\overline{B_{1/2}}$, we have
\[
 \psi(y)
 \ge
 \left(1-\frac1{2^{p'}}\right)^{p'}.
\]
Thus, setting
\[
 c_o
 :=
 \lambda\left(1-\frac1{2^{p'}}\right)^{p'}>0,
\]
we obtain \eqref{nn:eq:onset-normalized}.
Since $\lambda$  depends only on $N$ and $p$,
the same is true for $c_o$.
Returning to the original variables, we obtain
\[
 u(x,t)
 \le
 \sfM\exp\left[
   -c_o
   \left(\frac{\rho^p}{t-t_\ast}\right)^{\frac1{p-1}}
 \right]
\]
for every
$(x,t)\in B_{\rho/2}(x_o)\times
(t_\ast,t_\ast+\frac12\tau_o\rho^p]$.

Finally, recall that $\tau_o$ was chosen solely
in the construction of the barrier and depends
only on $N$ and $p$.
Suppose now that $u$ is defined only on
$B_\rho(x_o)\times(t_\ast,T)$, where
\( t_\ast<T<t_\ast+\tau_o\rho^p \).
For any
\[
 0<\sigma<
 \min\left\{
   \frac{T-t_\ast}{\rho^p},
   \frac{\tau_o}{2}
 \right\},
\]
we can carry out the same comparison on
$B_{1-h}\times(0,\sigma]$.
The barrier and all constants remain unchanged.
Thus, the estimate holds whenever
\(  t_\ast<t<T\) and \( t\le t_\ast+\frac12 \tau_o\rho^p \),
with equality at the latter time obtained
by continuity.
If $T\le t_\ast+\frac12 \tau_o\rho^p$
and $u$ extends continuously to $t=T$,
passing to the limit as $t\uparrow T$
also gives the estimate at $T$.
\end{proof}

\begin{remark}
\upshape
Exponential-type comparison function have already been used before; see, for example, \cite{DBKV} in the context of singular porous medium equations, or \cite{SV} for doubly nonlinear parabolic equations.
\end{remark}

\subsection{Spatial Lipschitz continuity}

We now combine the barrier estimate of Lemma~\ref{lem:barrier} with the gradient estimate on positive cylinders from Proposition~\ref{nn:prop:positive} to establish the following result, which is stronger than Theorem~\ref{thm:Lipschitz} and implies it as a consequence. In fact, it quantifies the decay of the gradient as the solution approaches the zero set.

\begin{theorem}\label{thm:lip-k}
Let $u$ be a locally bounded,
nonnegative local weak solution of \eqref{PDE}
in $\Omega_T$.
Let $k_o$ be defined in \eqref{nn:eq:k}.
For every $k\ge k_o$, there exists a constant
$C=C(N,p,k)>0$ with the following property.
If $Q_{2R}(z_o)\Subset\Omega_T$ and
\[
 0\le u\le\sfM
 \qquad\text{in }Q_{2R}(z_o),
\]
then
\[
 \big\|
   \nabla\bigl(u^{\frac1k}\bigr)
 \big\|_{L^\infty(Q_R(z_o))}
 \le
 \frac{C\sfM^{\frac1k}}R.
\]
Moreover,
\[
 |\nabla u|
 \le
 \frac{C\sfM^{\frac1k}}R\,
 u^{1-\frac1k}
 \qquad\text{a.e. in }Q_R(z_o).
\]
\end{theorem}

\begin{proof}
If $\sfM=0$, the estimate is immediate.
Otherwise, by Lemma~\ref{lem:scaling}, we consider the rescaled function 
\[
 v(y,\tau)
 :=
 \frac1{\sfM}u(x_o+Ry,t_o+R^p\tau),
\]
which is a solution to~\eqref{PDE} on $Q_2:=B_2(0)\times(-2^p,0]$ and satisfies $0\le v\le1$.

Fix $k=k_o$, with $k_o$ defined in
\eqref{nn:eq:k}, so that $k$ depends only on $N,p$.
Let $t_\ast\in[-2^p,0]$ be the transition time
given by Proposition~\ref{prop:Harnack}
applied on $Q_2$, i.e. $v$ vanishes
throughout $B_2$ before $t_\ast$ and is
strictly positive after $t_\ast$. 
If $t_\ast=0$, then $v\equiv0$ in $Q_2$,
and the estimate is immediate.
Otherwise, fix $(\bar x,\bar t)\in Q_1$. Let $\tau_o=\tau_o(N,p)>0$ be the parameter
from Lemma~\ref{lem:barrier}, and choose
\[
 0<\tau_1<
 \min\big\{
   2^p-1,\,
   \tfrac12\,2^{-p}\tau_o
 \big\}.
\]
We distinguish three cases.
Suppose first that $\bar t\ge t_\ast+\tau_1$,
so that $\bar t$ is separated from the
transition time by at least $\tau_1$.
The closed cylinder
\[
 \overline{B_{\frac12}(\bar x)}
 \times\left[\bar t-\tfrac12\tau_1,\bar t\right]
 \subset B_2(0)\times(t_\ast,0]
\]
lies entirely in the region where $u$ is
strictly positive. By continuity and compactness, there exists
$\sfm>0$ such that
\[
 \sfm\le v\le1
 \qquad\text{on }
 \overline{B_{\frac12}(\bar x)}
 \times\left[\bar t-\tfrac12\tau_1,\bar t\right].
\]
Therefore, we can apply
Proposition~\ref{nn:prop:positive}
with $r=\frac14$ and $h=\frac14\tau_1$.
Since $v^{\frac1k}\le1$, this yields
\[
 \big\|\nabla\bigl(v^{\frac1k}\bigr)\big\|_{L^\infty(
   B_{\frac14}(\bar x)\times
   (\bar t-\frac14\tau_1,\bar t))}
 \le
 C\left[
   4+\left(\frac4{\tau_1}\right)^{\frac1p}
 \right]
 \le C(N,p,k).
\]
Here $\tau_1$ is fixed in terms of $N,p$,
and the bound is independent of $\sfm$.

Next, suppose that
$t_\ast<\bar t<t_\ast+\tau_1$.
Thus, $(\bar x,\bar t)$ lies in the positivity
set, but its time coordinate is close to
the transition time.
Since $\bar t>-1$ and $\tau_1<2^p-1$, we have
\[
 t_\ast>\bar t-\tau_1
 >-1-(2^p-1)=-2^p.
\]
In particular, the initial time $t_\ast$
lies inside the time interval of $Q_2$,
and $v(\cdot,t_\ast)=0$ in $B_2$.
We apply Lemma~\ref{lem:barrier} with
$\rho=\frac12$, center $\bar x$, and
upper bound $\sfM=1$, using its version
for shorter time intervals if necessary.
This gives
\[
 v(x,t)
 \le
 \exp\left[
   -c_o
   \left(
     \frac{2^{-p}}{t-t_\ast}
   \right)^{\frac1{p-1}}
 \right]
\]
for every $(x,t)$ satisfying
\[
 x\in B_{\frac14}(\bar x),
 \qquad
 t_\ast<t\le
 \min\left\{
   0,\,
   t_\ast+\tfrac12\,2^{-p}\tau_o
 \right\}.
\]
Since
$\bar t-t_\ast<\tau_1<\tau_o/2^{p+1}$,
this estimate applies, in particular,
for all $t\in(t_\ast,\bar t]$.
Here $c_o$ depends only on $N$ and $p$. By the choice of $\tau_1$, we have
\[
 t_\ast
 <
 \tfrac12(t_\ast+\bar t)
 <
 \bar t
 <
 t_\ast+\tau_1
 <
 t_\ast+\tfrac12\,2^{-p}\tau_o.
\]
Since the right-hand side of the preceding
estimate is increasing in $t$, it follows that
\[
 \sup_{
   B_{\frac14}(\bar x)\times
   (\frac12(t_\ast+\bar t),\bar t)
 }v
 \le
 \exp\left[
   -\frac{c_1}{(\bar t-t_\ast)^{\frac1{p-1}}}
 \right]
 =:\sfM_\ast,
 \qquad
 c_1:=c_o\,2^{-\frac p{p-1}}>0.
\]
By positivity and continuity, $u$ also has
a positive lower bound on
\( \overline{B_{\frac14}(\bar x)}
 \times\left[\tfrac12(t_\ast+\bar t),\bar t\right] \).
Therefore, we can apply
Proposition~\ref{nn:prop:positive}
with $(r,h,s,\sfM)$ replaced by \(\big( \frac18,  \frac14(\bar t-t_\ast), \bar t,\sfM_\ast\big)\). 
This yields
\begin{align*}
 &\big\|\nabla\bigl(v^{\frac1k}\bigr)\big\|_{L^\infty(
   B_{\frac18}(\bar x)\times
   (\bar t-\frac14(\bar t-t_\ast),\bar t))}\\
 &\qquad\le
 C\exp\left[
   -\frac{c_1}
         {k(\bar t-t_\ast)^{\frac1{p-1}}}
 \right]
 \left[
   8+\left(\frac4{\bar t-t_\ast}\right)^{\frac1p}
 \right]\le C(N,p,k).
\end{align*}
The last inequality follows because the function
\[
 \sigma\longmapsto
 \exp\left[
   -\frac{c_1}{k\sigma^{\frac1{p-1}}}
 \right]
 \left[1+\frac1{\sigma^{\frac1p}}\right]
\]
is bounded above by a constant $C=C(N,p,k)$
for all $\sigma>0$.
In particular, the resulting constant is
independent of $t_\ast$ and $\bar t$.

Finally, if $\bar t\le t_\ast$, then
$v(\cdot,\bar t)\equiv0$ in $B_2(0)$.
Consequently, $u^{\frac1k}(\cdot,\bar t)$
is constant, and its spatial gradient vanishes.

The first two cases yield a uniform bound
$C(N,p,k)$ for the spatial gradient of
$v^{\frac1k}$ in the positive region
$Q_1\cap\{t>t_\ast\}$.
Indeed, this region can be covered by cylinders
of the type considered above, and the bounds
are independent of their location.
For $t\le t_\ast$, the function
$v^{\frac1k}(\cdot,t)$ is identically zero,
so its spatial gradient vanishes.
Combining these estimates, we obtain
\begin{equation*}
     \big\|\nabla\bigl(v^{\frac1k}\bigr)\big\|_
 {L^\infty(Q_1)}
 \le C(N,p,k).
\end{equation*}
Rescaling back to the original solution $u$ on $Q_{2R}(z_o)$, we have shown that
\begin{equation*}
     \big\|\nabla\bigl(u^{\frac1k}\bigr)\big\|_
 {L^\infty(Q_R(z_o))}
 \le \frac{C(N,p,k)\sfM^{\frac1k}}{R},
\end{equation*}
which is the first assertion of the theorem. 
In particular, this shows $\nabla\bigl(u^{\frac1k}\bigr)\le C\sfM^{\frac1k}/R$ a.e. in $Q_R(z_o)$. By the spatial Sobolev chain rule, we have 
$\nabla u= k u^{1-\frac1k}\nabla\bigl(u^{\frac1k}\bigr)$
almost everywhere in \(Q_R(z_o)\), which in combination with the preceding inequality yields the second claim. 
\end{proof}

At this point, Theorem~\ref{thm:Lipschitz} is an immediate consequence of the preceding Theorem~\ref{thm:lip-k}. 

\begin{proof}[\textbf{\upshape Proof of Theorem~\ref{thm:Lipschitz}}]
Set $\sfM:=\|u\|_{L^\infty(Q_{2R}(z_o))}$ and choose $k=k_o$, where $k_o$ is the parameter defined in~\eqref{nn:eq:k}.
From the second inequality of Theorem~\ref{thm:lip-k}, we have 
\begin{equation*}
    \|\nabla u\|_{L^\infty(Q_{R}(z_o))}
    \le 
    \frac{C\sfM^{\frac1k}}R\,
    \big\|u^{1-\frac1k}\big\|_{L^\infty(Q_{R}(z_o))}
    \le 
    \frac{C}R\,
    \|u\|_{L^\infty(Q_{2R}(z_o))},
\end{equation*}
with a constant $C=C(N,p)$. This completes the proof of Theorem~\ref{thm:Lipschitz}.
\end{proof}

\section{\texorpdfstring{$C^{1,\alpha}$}{C1alpha}-estimates}

For solutions to \eqref{PDE} that are strictly positive and bounded, regularity results for equations of parabolic $p$-Laplace type yield H\"older estimates for the gradient. We recall the following result from \cite[Theorem~6.8]{p-Schauder}.
\begin{lemma}\label{nn:lem:normalized-schauder}
Suppose $U$ is a weak solution to \eqref{PDE} in $Q_4$, satisfying
\[
 \tfrac12\le U\le2\quad\text{a.e. in }Q_4.
\]
 Then $U$ has a continuous
spatial gradient on $Q_1$. Moreover,  for some $\beta\in(0,1)$ and $C>1$ depending only on $N$ and $p$, we have
\begin{equation*}
 \sup_{Q_1}|\D U|\le C,\qquad
 |\D U(z)-\D U(\zeta)|\le C \dist_{\mathrm{par}}(z,\zeta)^\beta
\end{equation*}
for any $ z,\zeta\in Q_1 $.
\end{lemma}

We are now in a position to give the 

\begin{proof}[\textbf{\upshape Proof of Theorem~\ref{thm:holder}}]

The proof is divided into three steps. In the first two steps, we establish the qualitative property of continuity of the gradient, first on the positivity set and then on the zero set. In the final step, we then prove the claimed H\"older continuity of the gradient.

\noindent\textbf{Step 1: Continuous differentiability on the positivity set.}
In this section, we aim to show that $\D u$ is continuous throughout $\{u>0\}\cap Q_1$.
We work henceforth with $0\le u\le1$ in $Q_2=B_2\times(-2^p,0]$.
The gradient bound from Theorem~\ref{thm:lip-k} and
Proposition~\ref{prop:holder} give constants $L,\,H$ and $\alpha,\eta\in(0,1)$, depending
only on $N$ and $p$, such that on $Q_{\frac74}$ we have that
\begin{align}
 |u^\alpha(x,t)-u^\alpha(y,t)|&\le L|x-y|,
 &|\D u|&\le L u^{1-\alpha}\quad\hbox{a.e.},
 \label{nn:eq:normalized-power}\\
 |u(z)-u(\zeta)|&\le H \dist_{\mathrm{par}}(z,\zeta)^\eta.
 &&\label{nn:eq:normalized-holder}
\end{align}
We used a.e. for the gradient since we do not know yet if it is continuous.
Pick $\zeta=(y,s)\in Q_1$ such that $u(\zeta)>0$; it suffices to show that for any $\varep>0$, there is $\dl(\varepsilon)>0$, such that
\[
|\D u(z)-\D u(\zeta)|\le \varep
\]
for any $z\in Q_\dl(\z)$.
To this end, we first fix $c\in(0,1)$, depending only on $N$ and $p$, so small that
\begin{equation}\label{nn:eq:c-choice}
 4c<3/4,\qquad
 (4c)^p<(7/4)^p-1,\qquad H(8c)^\eta\le 1/2.
\end{equation}
For another point $z=(x,t) \in Q_1$, we set
\[
 \delta=\dist_{\mathrm{par}}(z,\zeta),\qquad
 \mathsf{m}=\max\{u(z),u(\zeta)\}\in(0,1],\qquad
 \rho=c\,\mathsf{m}^{\frac1\eta}.
\]
We consider the case $\delta\le\rho/2$. Denote by $z_\mathsf{m}=(x_\mathsf{m},t_\mathsf{m})$ either $z$ or $\z$
where the maximum $\mathsf{m}$ is attained, and set $t_+=\max\{t,s\}$ to the later time.
Use the cylinder
\begin{equation*}
 \mathcal C_\rho=B_{4\rho}(x_\mathsf{m})
                 \times(t_+-(4\rho)^p,t_+].
\end{equation*}
The first two conditions in \eqref{nn:eq:c-choice} and the fact $\mathsf{m}\le1$ place this cylinder in
$Q_{\frac74}$. Observe that every $w\in\mathcal C_\rho$ satisfies
$\dist_{\mathrm{par}}(w,z_m)\le8\rho$: the spatial separation is at most $4\rho$, and the
time separation is at most $(4\rho)^p$, because
$0\le t_+-t_{\mathsf{m}}\le\delta^p\le(\rho/2)^p$.
Hence, by \eqref{nn:eq:normalized-holder} and \eqref{nn:eq:c-choice}$_3$ we get for all $w\in\mathcal{C}_\rho$ that
\begin{align*}
     |u(w)-u(z_\mathsf{m})| \le H \dist_{\mathrm{par}}(w,z_\mathsf{m})^\eta\le H (8\rho)^\eta= H (8c)^\eta \mathsf{m}\le\tfrac12 \mathsf{m}
\end{align*}
which in turn yields that
\begin{equation}\label{nn:eq:comparable}
 \tfrac12\mathsf{m}\le u(w)\le\tfrac32\mathsf{m},
 \qquad \mbox{for all } w\in\mathcal C_\rho.
\end{equation}
Both  points $z$ and $\z$ lie in
$\overline B_\rho(x_\mathsf{m})\times[t_+-\rho^p,t_+]$ since we have assumed their distance $\dl$ to be less than $\rho/2$.

Define on $Q_4$
\[
 U(Y,S):=\frac{u(x_\mathsf{m}+\rho Y,t_++\rho^pS)}{\mathsf{m}},
\]
which again solves \eqref{PDE}.
By \eqref{nn:eq:comparable}, we have $1/2\le U\le3/2$.
First apply Lemma~\ref{nn:lem:normalized-schauder} to $U$ and then scale back to get that
\begin{equation}\label{nn:eq:positive-pair}
 |\D u(z)-\D u(\zeta)|
 \le C\frac{\mathsf{m}}\rho\left(\frac{\delta}{\rho}\right)^\beta
 =C\delta^\beta \mathsf{m}^{\,1-\frac{1+\beta}{\eta}},
\end{equation}
for any $z\in Q_1$,
as long as the distance verifies that $\delta=\dist_{\mathrm{par}}(z,\zeta)\le\tfrac12\rho =\tfrac12 c \mathsf{m}^{\frac1\eta}$.
Note that the factor $c^{-1-\beta}$ has been absorbed in the constant $C$ of the above estimate.
The assumption $\delta\le\frac12\rho$ is in particular satisfied whenever $\delta\le\frac12 cu(\zeta)^\frac{1}{\eta}$. Therefore, estimate \eqref{nn:eq:positive-pair} implies that $\nabla u$ is continuous in any point $\zeta\in Q_1$ with $u(\zeta)>0$. 

\noindent\textbf{Step 2: Continuous differentiability on the zero set. }
Let $z_o=(x_o,t_o)\in Q_1$ with $u(z_o)=0$.
For a small step-size $h\in\R^N\setminus\{0\}$, estimate \eqref{nn:eq:normalized-power}$_1$ at the fixed time $t_o$
gives
\begin{equation}\label{nn:eq:zero-differentiability}
 0\le u(x_o+h,t_o)
 \le L^{\frac1\alpha}|h|^{\frac1\alpha}.
\end{equation}
Since $ \alpha<1$, we divide \eqref{nn:eq:zero-differentiability} by $|h|$ and let $h\to0$ to get Fr\'echet differentiability
in the spatial variables at $z_o$, with $\D u(z_o)=0$.

On $\{u>0\}\cap Q_1$ both $u$ and $\D u$ are continuous. The almost everywhere
weighted bound in \eqref{nn:eq:normalized-power} consequently holds
pointwise there. With the derivative just defined at zeros it becomes
\begin{equation}\label{nn:eq:pointwise-power}
 |\D u(z)|\le L u(z)^{1-\alpha}\quad\hbox{for every }z\in Q_1.
\end{equation}
Consider a sequence $Q_1\ni z_j\to z_o$; use \eqref{nn:eq:pointwise-power} with $z_j$, send $j\to\infty$, recall the continuity of $u$ on $Q_1$, and use $u(z_o)=0$ to conclude that the right-hand side tends to zero.
Thus $\D u$ is jointly continuous across the zero set.
We remark that an unweighted Lipschitz estimate would not by itself suffice for this step.

\noindent\textbf{Step 3: H\"older continuity of the gradient.} 
It remains to obtain a uniform H\"older modulus for the gradient. We continue in the normalized
setting $0\le u\le 1$ in $Q_2$ and write
\begin{equation}\label{choice-A-s}
 A=\frac{1+\beta}{\eta}-1>0,\qquad
 s_*=\frac{\beta}{A+1-\alpha}
     =\frac{\beta\eta}{1+\beta-\alpha\eta},
\end{equation}
with the parameter $\beta\in(0,1)$ from Lemma~\ref{nn:lem:normalized-schauder} and the parameters $\alpha,\eta\in(0,1)$ from~\eqref{nn:eq:normalized-power} and \eqref{nn:eq:normalized-holder}. 
All these parameters are now fixed by $N$ and $p$.
Because $\alpha\eta<1$, we have $0<s_*<\eta$. Fix $c$ from \eqref{nn:eq:c-choice} and choose  
$\delta_o=\delta_o(N,p)\in(0,1)$ such that
\begin{equation*} 
 \delta_o^{\,1-\frac{s_*}{\eta}}\le \tfrac12 c.
\end{equation*}
For two points $z,\zeta\in Q_1$, recall the notation $\delta=\dist_{\mathrm{par}}(z,\zeta)$ and
$\mathsf{m}=\max\{u(z),u(\zeta)\}$. The case $\delta=0$ is immediate.
If $\delta\ge\delta_o$, \eqref{nn:eq:pointwise-power} yields the desired
estimate~\eqref{claim-holder} for a sufficiently large constant $C$. Suppose $0<\delta<\delta_o$.

If $\mathsf{m}\le\delta^{s_*}$, the pointwise weighted bound gives
\begin{equation}\label{nn:eq:glue-small}
 |\D u(z)-\D u(\zeta)|
 \le2L \mathsf{m}^{1-\alpha}
 \le2L\delta^{s_*(1-\alpha)}.
\end{equation}
This case includes pairs for which both values vanish.

If $\mathsf{m}>\delta^{s_*}$, put $\rho=c \mathsf{m}^{\frac1\eta}$. Then
\[
 \frac{\delta}{\rho}
 <c^{-1}\delta^{\,1-\frac{s_*}{\eta}}\le\tfrac12,
\]
so  construction from Step 1 and \eqref{nn:eq:positive-pair} apply.
In particular, that construction also proves positivity at the other
point; it is not assumed in advance. From \eqref{nn:eq:positive-pair} and $A>0$ we deduce 
\begin{equation}\label{nn:eq:glue-large}
 |\D u(z)-\D u(\zeta)|
 \le C\delta^\beta \mathsf{m}^{-A}
 \le C\delta^{\beta-s_*A}.
\end{equation}
The choice of $s_\ast$ in~\eqref{choice-A-s} balances the exponents:
\[
 s_*(1-\alpha)=\beta-s_*A
 =\frac{\beta\eta(1-\alpha)}{1+\beta-\alpha\eta}
 =:\gamma>0.
\]
With the exponent $\gamma$ fixed above, equations \eqref{nn:eq:glue-small} and \eqref{nn:eq:glue-large} prove
\begin{equation}\label{claim-holder}
|\D u(z)-\D u(\zeta)|\le C \dist_{\mathrm{par}}(z,\zeta)^\gamma \qquad\mbox{on $Q_1$.}
\end{equation}

In the general case $u\colon Q_{2R}(z_o)\to[0,M]$, we consider the rescaled function
\[
U(Y,S)=M^{-1}u(x_o+RY,t_o+R^pS)
\qquad\mbox{for }(Y,S)\in Q_2.
\]
This is a solution to~\eqref{PDE} with $0\le U\le 1$, so that~\eqref{claim-holder} is applicable with $U$ in place of $u$.
Since the spatial gradient scales by $\D u=(M/R)\D U$ and the metric by
$\dist_{\mathrm{par}}(Y,S)=\dist_{\mathrm{par}}(z,\zeta)/R$, we obtain the asserted bound \eqref{nn:eq:holder-main}. This completes the proof of the theorem.
\end{proof}

{\small \noindent{\bf AI disclosure.} The development of this work and the preparation of the manuscript were partly assisted by AI tools.}

\end{document}